\documentclass[11pt]{article}

\usepackage[utf8]{inputenc}

\usepackage{amssymb,xcolor,xspace}
\usepackage{amsmath,graphicx}
\usepackage{gastex}

\usepackage{hyperref}
\usepackage{array}

\usepackage{geometry} 

\usepackage{fancyhdr} 
\usepackage{sectsty}
\allsectionsfont{\sffamily\mdseries\upshape} 

\usepackage[nottoc,notlof,notlot]{tocbibind} 
\usepackage[titles]{tocloft} 

\def\cqfd{\skip10=\parfillskip\parfillskip=0pt
\enspace\hfill\symbolecqfd\par\parfillskip=\skip10\par\medskip}
\def\symbolecqfd{\rlap{$\sqcap$}$\sqcup$}

\newtheorem{theorem}{Theorem}[section]
\newtheorem{proposition}[theorem]{Proposition}
\newtheorem{lemma}[theorem]{Lemma}
\newtheorem{corollary}[theorem]{Corollary}

\newtheorem{pro-fact}[theorem]{Fact}

\newtheorem{pro-example}[theorem]{Example}
\newenvironment{example}{\begin{pro-example}\rm}{\cqfd\end{pro-example}}

\newtheorem{pro-remark}[theorem]{Remark}
\newenvironment{remark}{\begin{pro-remark}\rm}{\cqfd\end{pro-remark}}

\newenvironment{preuve}{\rm \trivlist \item[\hskip \labelsep{\bf
Proof.}]}{\cqfd\endtrivlist}

\def\cqfd{\skip10=\parfillskip\parfillskip=0pt
\enspace\hfill\symbolecqfd\par\parfillskip=\skip10\par\medskip}
\def\symbolecqfd{\rlap{$\sqcap$}$\sqcup$}
\def\proof{\begin{preuve}}
\def\eop{\end{preuve}}

\let\phi\varphi

\def\findex{\textsf{fi}}
\def\crff{\textsf{cr-fr}}
\def\ff{\textsf{fr}}
\def\fffi{\textsf{fr-fi}}

\def\inv{^{-1}}
\let\epsilon\varepsilon
\def \calA {\mathcal{A}}
\def \calB {\mathcal{B}}

\def \O {\mathcal{O}}

\def\calS  {\mathcal{S}}

\def\E{\mathbb{E}}

\def\N{\mathbb{N}}
\def \P {\mathbb{P}}
\def\PSL{\textsf{PSL}_2(\Z)}

\def\Z{\mathbb{Z}}

\def\Schreier{\textsf{Schreier}}

\def\path{\textsf{path}}

\def\Gpr{G_\text{pr}}
\def\gpr{g_\text{pr}}
\def\ella{\ell_2}
\def\ellb{\ell_3}
\def\ka{k_2}
\def\kb{k_3}

\def\Gpra{\Gpr^{(a)}}
\def\Gprb{\Gpr^{(b)}}
\def\Gprz{\Gpr^{(0)}}
\def\Gpro{\Gpr^{(1)}}
\def\Gprv{\Gpr^{(v)}}

\def\gpra{\gpr^{(a)}}
\def\gprb{\gpr^{(b)}}
\def\gprz{\gpr^{(0)}}
\def\gpro{\gpr^{(1)}}
\def\gprv{\gpr^{(v)}}

\usepackage{todonotes}

\def\pw#1{\textcolor{purple}{PW -- #1 --}\xspace}

\title{Statistics of subgroups of the modular group}

\author{
    Fr\'ed\'erique Bassino, \small{\url{bassino@lipn.fr}}\\
    \small{Universit\'e Sorbonne Paris Nord, LIPN, CNRS UMR 7030, F-93430 Villetaneuse, France}%
    \and
    Cyril Nicaud, \small{\url{cyril.nicaud@u-pem.fr}}\\
    \small{LIGM, Univ Gustave Eiffel, CNRS, ESIEE Paris, F-77454, Marne-la-Vallée, France}%
    \and
    Pascal Weil, \small{\url{pascal.weil@labri.fr}}\\
    \small{Univ. Bordeaux, LaBRI, CNRS UMR 5800, F-33400 Talence, France}\thanks{%
    LaBRI, Univ. Bordeaux, 351 cours de la Lib\'eration, 33400 Talence, France.}\\
    \small{CNRS, ReLaX, UMI 2000, Siruseri, India}
    }

\begin{document}

\maketitle

\begin{abstract}
We count the finitely generated subgroups of the modular group $\PSL$. More precisely: each such subgroup $H$ can be represented by its Stallings graph $\Gamma(H)$, we consider the number of vertices of $\Gamma(H)$ to be the size of $H$ and we count the subgroups of size $n$. Since an index $n$ subgroup has size $n$, our results generalize the known results on the enumeration of the finite index subgroups of $\PSL$. We give asymptotic equivalents for the number of finitely generated subgroups of $\PSL$, as well as of the number of finite index subgroups, free subgroups and free finite index subgroups. We also give the expected value of the isomorphism type of a size $n$ subgroup and prove a large deviation statement concerning this value. Similar results are proved for finite index and for free subgroups. Finally, we show how to efficiently generate uniformly at random a size $n$ subgroup (resp. finite index subgroup, free subgroup) of $\PSL$.
\end{abstract}
\tableofcontents

\section{Introduction}

The modular group $\PSL$ is a central object in algebra, notably because of its relevance in geometry. It can be viewed as the quotient of the group $\textsf{SL}_2(\Z)$ of $(2,2)$-matrices with integer coefficients and determinant $1$ by its center (the 2-element subgroup $\{\pm \textsf{Id}_2\}$), as the group of integer coefficient matrices generated by
$$\left(\begin{array}{cc}0 & 1 \\-1 & 0\end{array}\right)\textrm{ and }\left(\begin{array}{cc}1 & 1 \\0 & 1\end{array}\right),$$
or as a set of isometries of the hyperbolic plane. The more relevant view in this paper is as the free product $\PSL = \Z_2 \ast \Z_3 = \langle a,b \mid a^2 = b^3 = 1\rangle$.

$\PSL$ is a standard example of a virtually free group, and of a group acting discretely on the hyperbolic plane. The latter property justifies the more than century-old interest in $\PSL$ and in the classification of its subgroups (\textit{e.g.}, Klein and Fricke \cite{1966:Klein-I,1966:Klein-II}). There is indeed a vast literature on the finite index subgroups of $\PSL$ and their enumeration, exact and asymptotic, later extended to Hecke groups (groups of the form $\Z_2 \ast \Z_q$), free products of finitely many cyclic groups (finite or infinite), and to Fuchsian groups. We will review this literature below.

The first objective of the present paper is to extend these results to all finitely generated subgroups of $\PSL$, for which we give exact enumeration formulas and asymptotic estimates.  We then exploit this first set of results to show how one can efficiently (in linear time) generate a subgroup of $\PSL$ uniformly at random. We also compute the expected value of the isomorphism type of a subgroup and show a concentration result around this expected type. All these results can be relativized to finite index subgroups (retrieving known formulas and asymptotic equivalents), to free subgroups and to free finite index subgroups.

As mentioned above, the enumeration and asymptotic estimation of the number of finite index subgroups of free products of cyclic groups has a long and rich history.
In 1949, Hall \cite{1949:Hall} gave an explicit formula for the number of finite index subgroups of a finite rank free group (a free product of a finite number of copies of $\Z$). 
In 1965, Dey \cite{1965:Dey} gave a recurrence relation for the number of index $n$ subgroups of $\PSL$, and more generally of a free product of finitely many finite groups. Hall's and Dey's results were then extended by Stothers in 1978 \cite{1978:Stothers-MathsComp}, who counted the finite index subgroups and the free finite index subgroups of a free product of finitely many cyclic groups. Shortly thereafter, 
Newman \cite{1976:Newman} refined Stothers's results and gave a new recurrence relation and an asymptotic equivalent for the number of finite index subgroups of a free product of cyclic groups. In the case of $\PSL$, this equivalent is
$$\left(12\pi\ e^{\frac12}\right)^{-\frac12}\ \exp\left(\frac 16 n\log n - \frac16n + n^{\frac12} + n^{\frac13} + \frac12\log n\right),$$
which we retrieve with our methods (Eq.~(\ref{eq: Hnfi}), Section~\ref{sec: finite index subgroups}).

In 1996, Müller \cite{1996:Muller} (see also \cite{2003:LubotzkySegal}) further generalized this result by giving the complete asymptotics of the number $s_n(G)$ of index $n$ subgroups of a free product $G$ of finitely many cyclic (finite and free) groups. This was extended by Mednykh in 1979 \cite{1979:Mednyh} and by Müller and Schlage-Puchta in 2002 \cite{2002:MullerPuchta} to the case where $G$ is a surface group, and in 2004 by Liebeck and Shalev \cite{2004:LiebeckShalev} to the case of Fuchsian groups (again, see \cite{2003:LubotzkySegal}). In 2019, Ciobanu and Kolpakov \cite{2017:CiobanuKolpakov} exploited the count of the free finite index subgroups of $\Z_2 \ast \Z_2 \ast \Z_2$ to enumerate 3-dimensional maps.

In 2004, Müller and Schlage-Puchta \cite{2004:MullerSchlage-Puchta} studied the isomorphism types of finite index subgroups of free products of cyclic groups (see also \cite{2010:MullerSchlage-Puchta}, and \cite{2007:MullerSchlage-Puchta} for an extension to Fuchsian groups). 
Let $G = C_{p_1}^{\ast e_1} \ast \cdots \ast C_{p_t}^{\ast e_t} \ast F_r$, where $C_q$ is the cyclic group of order $q$, the $p_i$ are pairwise distinct, the superscript ${}^{\ast e}$ means the free product of $e$ copies, and $F_r$ is the rank $r$ free group. Kurosh's theorem states that a finitely generated subgroup $H$ of $G$ is isomorphic to a free product of the form $C_{p_1}^{\ast \lambda_1} \ast \cdots \ast C_{p_t}^{\ast \lambda_t} \ast F_\mu$. The tuple $(\lambda_1,\ldots,\lambda_t,\mu)$ is called the \emph{isomorphism type} of $H$. Müller and Schlage-Puchta characterize the tuples that can arise as the isomorphism type of a finite index subgroup of $G$, and they describe their (asymptotic) distribution. In the case where $G = \PSL$, these isomorphism types are exactly the triples $(\ella,\ellb,r)$ such that $3\ella+4\ellb+6r > 7$, and the index is then $n = 3\ella + 4\ellb + 6(r -1)$.

As we said earlier, we extended these results to \emph{all} finitely generated subgroups of $\PSL$. This is made possible by a systematic usage of the so-called Stallings graph of a subgroup: a finite graph uniquely associated with the subgroup, which coincides with the Schreier graph (or coset graph) in the finite index case. Stallings graphs were first introduced by Stallings \cite{1983:Stallings} for finitely generated subgroups of free groups, and later extended to free products of finite groups (Markus-Epstein \cite{2007:Markus-Epstein}) and to quasi-convex subgroups of hyperbolic groups (Kharlampovich \textit{et al.} \cite{2017:KharlampovichMiasnikovWeil}). The consideration of Stallings graphs provides a unified approach to the (often efficient) solution of algorithmic problems, see among other \cite{1983:Stallings,1983:Gersten,2000:BirgetMargolisMeakin,2001:MargolisSapirWeil,2002:KapovichMyasnikov,2007:RoigVenturaWeil,2007:MiasnikovVenturaWeil,2008:BassinoNicaudWeil,2013:BassinoMartinoNicaud}.

In general, the Stallings graph of a subgroup is a connected graph with a designated base vertex, whose edges are labeled by the generators of the ambient graph. The number of vertices of this graph is taken to be a measure of the \emph{size} of the subgroup. In the specific case of the subgroups of $\PSL$, the edge labels are $a$ and $b$ (the order 2 and order 3 generators of $\PSL$), all vertices except perhaps the base vertex are adjacent to an $a$-labeled and to a $b$-labeled edge, the $a$-labeled edges represent an involution of the set of vertices to which they are adjacent, and the $b$-labeled edges represent a partial injection on the vertex set, whose square is a permutation of its domain, with orbits of size 1 or 3. The Stallings graph of a subgroup of $\PSL$ can be very efficiently computed given a set of generators for the subgroup, and it allows a simple characterization of the finite index property, freeness or the isomorphism type of a subgroup (see Sections~\ref{sec: PSL-reduced graphs} and~\ref{sec: properties}).

We compute precisely the number of finitely generated subgroups of $\PSL$ of a given size in Section~\ref{sec: counting}. The strategy for this consists in counting their Stallings graphs, which reduces to counting the $\PSL$-cyclically reduced graphs (those in which every vertex is adjacent to edges corresponding to the two generators $a$ and $b$ of $\PSL$). For this, we use the classical tools of analytic combinatorics (see \cite{2009:FlajoletSedgewick}), especially the notion of exponential generating series (EGS).

Counting is a stepping stone for the task of evaluating the asymptotic behavior of size $n$ subgroups or estimating the asymptotic probability of a property for the uniform random distribution over these subgroups, which is done in Section~\ref{sec: asymptotics}. 
The path to reach our results in this regard goes through two combinatorial digressions which yield results that are interesting in their own right. In the first, Section~\ref{sec: combinatorial digression 1}, we compute asymptotic equivalents and expected values, and we prove strong concentration results (in the form of large deviations theorems) for discrete objects that are disjoint unions of structures taken from a finite set: the sets of $a$-labeled edges and of $b$-labeled edges in the Stallings graph of a subgroup of $\PSL$ fall in this category, and the specific results for these structures are detailed in Section~\ref{sec: series T2 and T3}.

Next, we observe that a suitable set of $a$-labeled edges and a suitable set of $b$-labeled edges define a Stallings graph exactly if, together, they form a connected graph. Part of the computation of an asymptotic equivalent for the number of subgroups of $\PSL$ consists in showing that the probability that a set of $a$-edges and a set of $b$-edges on $n$ vertices define a connected graph tends to 1 when $n$ tends to infinity. This is done, in a more general setting, in our second combinatorial digression, Section~\ref{sec: asymptotics col}.

Another remark is that the isomorphism type of a subgroup $H$ is a computable triple $(\ella, \ellb,r)$, where $\ella$ (resp. $\ellb$) is the number of fixed points of the partial injection determined by the $a$-labeled (resp. $b$-labeled) edges of the Stallings graph of $H$ (see Proposition~\ref{prop: compute isomorphism type}).

This set of results allows us to compute asymptotic equivalents of the number of size $n$ subgroups (Section~\ref{sec: number of subgroups}) and of the expected isomorphism type of a size $n$ subgroup, as well as to prove a large deviations property for this type~(Section~\ref{sec: asymptotics expected values}).
In Section~\ref{sec: random generation}, our results are used to present a practical, low complexity random generation algorithm for finitely generated subgroups of $\PSL$ of a given size.

In Section~\ref{sec: finite and free}, we investigate in the same spirit particular classes of subgroups of $\PSL$. We obtain the same sort of enumeration and asymptotic properties as in the general case for finite index subgroups (here, most results were already known except for the large deviations property), for free subgroups or for free finite index subgroups. This is done by observing that freeness and finite index are characterized by purely combinatorial properties of the Stallings graphs, which are sufficiently smooth to fall into the classes covered by our combinatorial digressions, and can therefore be handled by the same strategy as the general case.

Appendix~\ref{app:coef} gives additional details about the methods used for exact enumeration computations, and shows the numbers of size $n$ subgroups (resp. index $n$ subgroups, size $n$ free subgroups) for small values of $n$.

\section{Combinatorial description of the subgroups of $\PSL$}\label{sec: subgroups}

Let $A = \{a,b\}$, and let $\tilde A = A \sqcup A\inv$. Every word on the alphabet $\tilde A$ represents an element of $\PSL = \langle a,b \mid a^2 = b^3 = 1\rangle \sim \Z_2 \ast \Z_3$, but different words can represent the same element. We say that a word $u$ is \emph{geodesic} (for $\PSL$) if it is a shortest word representing its value in $\PSL$: the geodesic words are exactly those words on $\tilde A$ which do not contain $aa\inv$, $a\inv a$, $bb\inv$, $b\inv b$ (such words are called \emph{freely reduced}) and also do not contain $a^2$, $a^{-2}$, $b^2$, $b^{-2}$. Equivalently, the geodesic words are those which consist of alternations of $a^{\pm1}$ and $b^{\pm1}$. Finding a geodesic representative for the element represented by a word $u$ is easily done: one first freely reduces $u$ (iteratively deleting every factor of the form $xx\inv$ with $x\in \tilde A$), then deletes every occurrence of $a^2$, $b^3$ and their inverses, and finally replace every occurrence of $b^2$ (resp. $b^{-2}$) by $b\inv$ (resp. $b$). Replacing every occurrence of $a\inv$ by $a$ in a geodesic word yields another geodesic representative for the same element, which is an alternation of $a$ and $b^{\pm 1}$. Such words are called \emph{shortlex geodesic}, and every element of $\PSL$ has a unique shortlex geodesic representative. Shortlex geodesic words provide a quick solution to the word problem in $\PSL$: two words represent the same element if and only if they have the same shortlex geodesic representative. 

\subsection{Stallings graphs}

It is well-known that $\PSL$ is hyperbolic, virtually free and that every one of its finitely generated subgroups is quasi-convex. It follows that every finitely generated subgroup $H$ of $\PSL$ is uniquely represented by a finite, labeled, rooted graph $(\Gamma(H),v)$, called its \emph{Stallings graph} \cite{2017:KharlampovichMiasnikovWeil} (see also \cite{2007:Markus-Epstein,2016:SilvaSoler-EscrivaVentura}), which can be efficiently computed, see below. The definition of the Stallings graph of a subgroup $H$ is as follows: we first consider the graph $\Schreier(\PSL,H)$, whose vertices are the cosets $Hg$ ($g\in \PSL$), with an $a$-labeled edge from $Hg$ to $Hga$ and a $b$-labeled edge from $Hg$ to $Hgb$ for every $g\in \PSL$. The notion of a labeled edge is extended to a labeled path in $\Schreier(\PSL,H)$: there is an empty path labeled by the empty word from $p$ to $p$ for every vertex $p$ of $\Schreier(\PSL,H)$; moreover, if $\pi$ is a path from vertex $p$ to vertex $q$ labeled by a word $u$, and $e$ is a $c$-labeled edge from $q$ to vertex $r$ ($c \in A$), then $\pi e$ is a path from $p$ to $r$ labeled $uc$. It is clear that a word $u \in \tilde A^*$ labels a loop at vertex $H$ in  $\Schreier(\PSL,H)$ if and only if $u$ represents an element of $H$.  The Stallings graph $\Gamma(H)$ is the (uniquely determined) subgraph of $\Schreier(\PSL,H)$, rooted at vertex $H$ and spanned by the loops at $H$, which are labeled by geodesic representatives of the elements of $H$. Because $H$ is quasi-convex, this graph is finite \cite{1991:GerstenShort}, see also \cite{2017:KharlampovichMiasnikovWeil} for an algorithmic point of view.

In \cite{2017:KharlampovichMiasnikovWeil} the general algorithm to compute the Stallings graph of a quasi-convex subgroup of a hyperbolic group (given a tuple of generators of $H$, in the form of words over the alphabet $\tilde A$) is rather complex, but it is much simpler in the case of subgroups of $\PSL$. We now give a short description of this algorithm. A crucial operation in this algorithm is that of \emph{edge folding} \cite{1983:Stallings}: if in a labeled graph $\Delta$, we have two edges with the same label $c \in A$, from vertex $p$ to vertices $q$ and $r$ (resp. from vertices $q$ and $r$ to vertex $p$), we identify the vertices $q$ and $r$ and keep only one $c$-labeled edge from $p$ to $q = r$ (resp. from $q = r$ to $p$).

Let $H$ be the subgroup of $\PSL$ generated by a tuple $\vec h = (h_1,\ldots, h_n)$, where each $h_i$ is a word on $\tilde A$. We may assume that the $h_i$ are non-empty and shortlex geodesic (or we can rewrite them in this form).

\paragraph*{Step 1.}
Compute $\Gamma_F(\vec h)$, the Stallings graph of the subgroup generated by $\vec h$ in the free group $F(A)$ \cite{1983:Stallings}. To this end, we first form a loop labeled by each $h_i$, around a common vertex $1$. The loop labeled $h_i$ has $|h_i|$ edges, labeled by the letters of $h_i$. We replace every edge from $p$ to $q$ labeled by  $b\inv$ by a $b$-labeled edge from $q$ to $p$. We then iteratively apply the edge folding operation, as long as there are foldable pairs of edges. The resulting graph is $\Gamma_F(\vec h)$ (rooted at vertex $1$).

\paragraph*{Step 2.}
For each $a$-labeled edge in $\Gamma_F(\vec h)$ from $p$ to $q$, add an $a$-labeled edge from $q$ to $p$. Similarly, for each $b$-labeled edge from $p$ to $q$, we add a new vertex $r$, and $b$-labeled edges from $q$ to $r$ and from $r$ to $p$.

We then iteratively apply the folding operation as long as it is possible; and we finally remove every vertex that is different from the root vertex, which has an incoming and an outgoing $b$-edge, and which is not adjacent to an $a$-edge. The resulting graph, rooted at 1, is (isomorphic to) $(\Gamma(H),H)$.

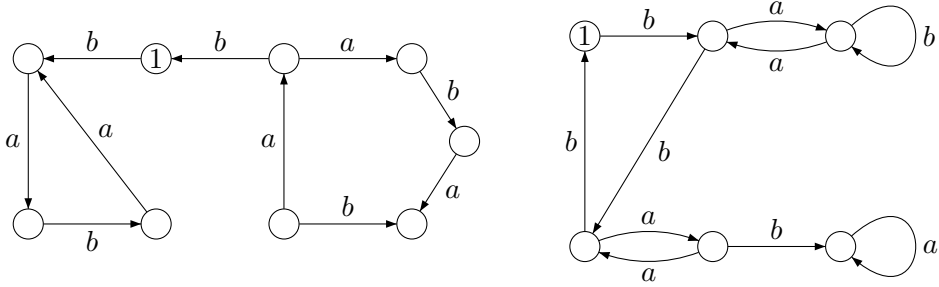
\begin{figure}[htbp]
\centering
\begin{picture}(88,35)(5,-37)
\gasset{Nw=4,Nh=4}

\node(n0)(3.0,-7.0){$1$}

\node(n1)(-14.0,-7.0){}

\node(n2)(-14.0,-29.0){}

\node(n3)(3,-29.0){}

\node(n4)(20,-7.0){}

\node(n5)(37.0,-7.0){}

\node(n6)(44.0,-18.0){}

\node(n7)(37.0,-29.0){}

\node(n8)(20.0,-29.0){}

\drawedge[ELside=r](n0,n1){$b$}

\drawedge[ELside=r](n1,n2){$a$}

\drawedge[ELside=r](n2,n3){$b$}

\drawedge[ELside=r](n3,n1){$a$}

\drawedge[ELside=r](n4,n0){$b$}

\drawedge(n4,n5){$a$}

\drawedge(n5,n6){$b$}

\drawedge(n6,n7){$a$}

\drawedge(n8,n7){$b$}

\drawedge(n8,n4){$a$}
\node(n0)(60.0,-4.0){$1$}

\node(n1)(77.0,-4.0){}

\node(n2)(60.0,-32.0){}

\node(n3)(94,-4.0){}

\node(n4)(77,-32.0){}

\node(n5)(94.0,-32.0){}

\drawedge[curvedepth=2.0](n1,n3){$a$}
\drawedge[curvedepth=2.0](n3,n1){$a$}

\drawedge[curvedepth=2.0](n2,n4){$a$}
\drawedge[curvedepth=2.0](n4,n2){$a$}

\drawedge(n0,n1){$b$}

\drawedge(n1,n2){$b$}

\drawedge(n2,n0){$b$}

\drawedge(n4,n5){$b$}

\drawloop[loopangle=0.0](n3){$b$}
\drawloop[loopangle=0.0](n5){$a$}

\end{picture}
\caption{\small The rooted labeled graphs $\Gamma_F(\vec h_2)$ and $\Gamma(H_2)$.}\label{fig: Stallings graphs H2}
\end{figure}

\begin{example}\label{ex: Stallings graph}
Let $H_1$ (resp. $H_2$) be the subgroups of $\PSL$ generated by the tuple $\vec h_1 = (abab\inv, babab)$ (resp.
$\vec h_2 = (babab\inv,b\inv abab\inv ab)$). The rooted labeled graphs $\Gamma_F(\vec h_1)$ and $\Gamma(H_1)$ (resp. $\Gamma_F(\vec h_2)$ and $\Gamma(H_2)$) are given in Figures~\ref{fig: Stallings graphs H1} and~\ref{fig: Stallings graphs H2}, respectively.
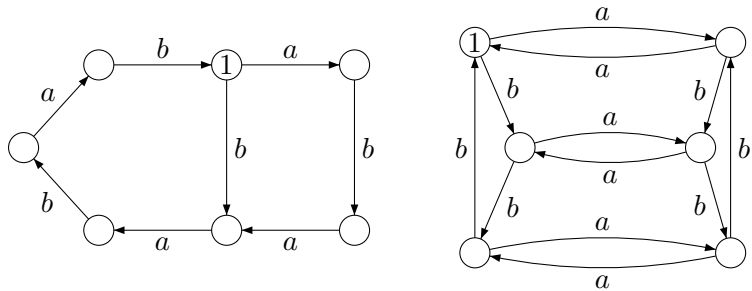
\begin{figure}[htbp]
\centering
\begin{picture}(88,37)(5,-37)
\gasset{Nw=4,Nh=4}

\node(n0)(27.0,-7.0){$1$}

\node(n1)(44.0,-7.0){}

\node(n2)(44.0,-29.0){}

\node(n3)(27,-29.0){}

\node(n4)(10,-29.0){}

\node(n5)(0.0,-18.0){}

\node(n6)(10.0,-7.0){}

\drawedge(n0,n1){$a$}

\drawedge(n0,n3){$b$}

\drawedge(n1,n2){$b$}

\drawedge(n2,n3){$a$}

\drawedge(n3,n4){$a$}

\drawedge(n4,n5){$b$}

\drawedge(n5,n6){$a$}

\drawedge(n6,n0){$b$}
\node(n0)(60.0,-4.0){$1$}

\node(n1)(66.0,-18.0){}

\node(n2)(60.0,-32.0){}

\node(n3)(94,-4.0){}

\node(n4)(90,-18.0){}

\node(n5)(94.0,-32.0){}

\drawedge(n0,n1){$b$}

\drawedge(n1,n2){$b$}

\drawedge(n2,n0){$b$}

\drawedge[ELside=r](n3,n4){$b$}

\drawedge[ELside=r](n4,n5){$b$}

\drawedge[ELside=r](n5,n3){$b$}

\drawedge[curvedepth=2.0](n0,n3){$a$}

\drawedge[curvedepth=2.0](n3,n0){$a$}

\drawedge[curvedepth=2.0](n1,n4){$a$}

\drawedge[curvedepth=2.0](n4,n1){$a$}

\drawedge[curvedepth=2.0](n2,n5){$a$}

\drawedge[curvedepth=2.0](n5,n2){$a$}
\end{picture}
\caption{\small The rooted labeled graphs $\Gamma_F(\vec h_1)$ and $\Gamma(H_1)$}\label{fig: Stallings graphs H1}
\end{figure}
\end{example}

In view of the algorithm to compute the Stallings graph of a subgroup of $\PSL$ outlined above, we note the following. If $H$ is a subgroup of $\PSL$ and $g \in \PSL$ has shortlex representative $u_g$, the Stallings graph of $H^g = g\inv Hg$ is obtained from the Stallings graph of $H$ as follows: add to $\Gamma(H)$ a $u_g$-labeled path, starting at the root vertex of $\Gamma(H)$ and otherwise consisting of $|u_g|$ new vertices; move the root to the end vertex of that path, and apply Step 2 above.

\subsection{$\PSL$-reduced and $\PSL$-cyclically reduced graphs}\label{sec: PSL-reduced graphs}

We now describe the graphs which may arise as Stallings graphs of finitely generated subgroups of $\PSL$.

An \emph{$A$-labeled graph} is a finite, directed graph, whose edges are labeled by elements of $A$. It is \emph{rooted} if it comes with a designated (base) vertex.  It is \emph{freely reduced} if no two edges out of (resp. into) the same vertex have the same label. A rooted $A$-labeled graph is $\PSL$\emph{-reduced} if it satisfies the following properties:
\begin{itemize}
\item[(1)] it is connected and freely reduced;

\item[(2)] for every $a$-labeled edge from $p$ to $q$, there is an $a$-edge from $q$ to $p$;

\item[(3)] if there are two consecutive $b$-edges, from $p$ to $q$ and from $q$ to $r$, then there is a third $b$-edge from $r$ to $p$, forming a \emph{$b$-triangle} (if $p, q, r$ are distinct);

\item[(4)] every vertex is adjacent to an $a$-edge and to a $b$-edge, except maybe the base vertex.
\end{itemize}

In the sequel, when drawing a $\PSL$-reduced graph, we will represent a pair of $a$-edges from $p$ to $q$ and from $q$ to $p$ ($p\ne q$) by a single undirected $a$-labeled edge between $p$ and $q$, which we will refer to as an \emph{isolated $a$-edge}. We will also refer to $b$-edges that are neither loops nor in a $b$-triangle as \emph{isolated $b$-edges}.

\begin{example}
The graph $\Gamma(H_2)$ in Figure~\ref{fig: Stallings graphs H2} exhibits an $a$- and a $b$-loop, two isolated $a$-edges, an isolated $b$-edge and a $b$-triangle.
\end{example}

An unrooted $A$-labeled graph is $\PSL$\emph{-cyclically reduced} if it is $\PSL$-reduced for every choice of a base vertex. Equivalently, an unrooted $A$-labeled graph is $\PSL$-cyclically reduced if it has a single vertex, or it satifies Properties (1), (2), (3) above, and the following property:
\begin{itemize}
\item[($4'$)] every vertex is adjacent to an $a$-edge and to a $b$-edge.
\end{itemize}

The following statement directly follows from the discussion we have carried out.

\begin{lemma}
An $A$-labeled graph is the Stallings graph of a finitely generated subgroup of $\PSL$ if and only if it is $\PSL$-reduced.

A class of conjugacy of finitely generated subgroups of $\PSL$ is uniquely characterized by a $\PSL$-cyclically reduced graph.
\end{lemma}

\begin{remark}
We will see in Section~\ref{sec: properties} how to effectively compute (a generating set of) the subgroup $H$ of $\PSL$ of which a given $\PSL$-reduced graph $(\Gamma,1)$ is the Stallings graph. The $\PSL$-cyclically reduced graph $\Delta$ which characterizes the conjugacy class of $H$ is obtained from $\Gamma$ by forgetting the base vertex, and iteratively removing every vertex that is not adjacent to an $a$-edge and a $b$-edge.
\end{remark}

We will say that a $\PSL$-cyclically reduced (resp. $\PSL$-reduced) graph is \emph{proper} if it satisfies Properties (1), (2), (3) and ($4'$). 

\begin{example}\label{ex: size 1 subgroups}
There are  three non-proper $\PSL$-cyclically reduced graphs, see Figure~\ref{fig: 1 vertex non-proper}. Each has a single vertex, is $\PSL$-reduced (when rooted at its only vertex), and they correspond to the subgroups $H = 1, \langle a\rangle \textrm{ and } \langle b\rangle$.
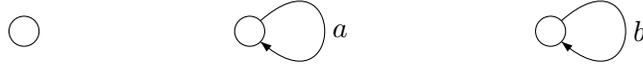
\begin{figure}[htbp]
\centering
\begin{picture}(76,10)(5,-13)
\gasset{Nw=4,Nh=4}

\node(n0)(3.0,-10.0){}

\node(n1)(33.0,-10.0){}
\drawloop[loopangle=0.0](n1){$a$}

\node(n2)(73.0,-10.0){}
\drawloop[loopangle=0.0](n2){$b$}

\end{picture}
\caption{\small The non-proper $\PSL$-cyclically reduced graphs.}\label{fig: 1 vertex non-proper}
\end{figure}
In addition, there is one proper 1-vertex $\PSL$-cyclically reduced graph: it has both an $a$- and a $b$-loop, and it corresponds to the subgroup $\PSL$. In particular, $\PSL$ has exactly 4 size 1 subgroups.
\end{example}

Finally we define the \emph{combinatorial type} of a $\PSL$-cyclically reduced or a $\PSL$-reduced graph $\Gamma$ to be the tuple $(n,\ka,\kb,\ella,\ellb,m)$, where:
\begin{itemize}
\item $n$ is the number of vertices (called the \emph{size} of $\Gamma$);
\item $\ka$, the number of isolated $a$-edges;
\item $\kb$, the number of isolated $b$-edges;
\item $\ella$, the number of $a$-loops;
\item $\ellb$, the number of $b$-loops;
\item and $m$, the number of $b$-triangles.
\end{itemize}
%

\subsection{Reading properties of a subgroup off its Stallings graph}\label{sec: properties}

Let $H$ be a finitely generated subgroup of $\PSL$. We show how such properties of $H$ as finite index, isomorphism type or freeness can be read directly off the combinatorial type of the Stallings graph of $H$.

\subsubsection{Combinatorial type}
The next result determines which tuples of integers are the combinatorial type of a proper $\PSL$-cyclically reduced graph. It is analogous in spirit to early results of Millington \cite{1969:Millington-2,1969:Millington-1} and Stothers \cite{1974:Stothers} on so-called impossible specifications of finite index subgroups of $\PSL$ (see Section~\ref{sec: finite index Stallings} for the connection with finite index).

\begin{proposition}\label{prop: realizability combinatorial type}
Let $\gamma = (n,\ka,\kb,\ella,\ellb,m)$ be a tuple of natural integers. There exists a proper $\PSL$-cyclically reduced graph with combinatorial type $\gamma$ if and only if $n = 2\ka + \ella = 2\kb + \ellb + 3m$ and $m-\ella-\ellb$ is even and greater than or equal to $-2$.
\end{proposition}

\proof
It is easily verified that, if $\Gamma$ is proper $\PSL$-cyclically reduced, we have
$$n = 2\ka + \ella = 2\kb + \ellb + 3m.$$
It follows that $n = \ka + \kb + \frac12(3m + \ella + \ellb)$, and therefore $m - \ella - \ellb$ is even.

In addition, let $T$ be a spanning tree of $\Gamma$ containing two edges of every $b$-triangle (it is not difficult
to verify that every forest in $\Gamma$ can be completed to a spanning tree). Then $T$ has $n-1$ edges, namely $2m$ from the $b$-triangles, some of the $\ka$ isolated $a$-edges and some of the $\kb$ isolated $b$-edges. We have therefore
\begin{align*}
&2m + \ka +\kb \enspace\ge\enspace \ka + \kb + \frac12(3m + \ella + \ellb) - 1, \\
&\textrm{and hence }m-\ella-\ellb \enspace\ge\enspace -2.
\end{align*}

To prove that the condition is sufficient, we consider a tuple $\gamma = (n,\ka,\kb,\ella,\ellb,m)$ satisfying the conditions in the statement. We distinguish three cases, depending whether $m-\ella-\ellb$ is equal to $-2$, is equal to 0 or is greater than or equal to 2.

We observe that $2\ka + \ella = 2\kb + \ellb + 3(m-\ella-\ellb) + 3(\ella + \ellb)$, which implies that $\ka = \kb + 2\ellb + \ella + \frac32(m-\ella-\ellb)$.

\paragraph*{Case $m-\ella-\ellb = -2$.}
Suppose that $m \ge 1$. We consider graphs $\Delta_1, \ldots, \Delta_m$, consisting each of a $b$-triangle, at one vertex of which we attach either an $a$-loop, or an isolated $a$-edge carrying a $b$-loop. This is done in such a way that the maximum possible number of $b$-loops are used (all of them unless $\ellb > m$, that is, $\ella \le 1$). If $\Delta_1$ carries an $a$-loop, we modify it by inserting a path labeled $(ab)^{\kb}$ between the triangle and the $a$-loop. If instead $\Delta_1$ carries a $b$-loop, we replace the isolated $a$-edge by a path labeled $(ab)^{\kb}a$. We then construct $\Delta$ by connecting the $\Delta_i$ in a line, with isolated $a$-edges between them. Note that $\Delta$ has two vertices, at each extremity, which are on a $b$-triangle but not adjacent to an $a$-edge. We call these \emph{open} vertices. 

If $\ellb \le m$, we ensured that $\Delta$ contains $m$ $b$-triangles, $\ellb$ $b$-loops, $\kb$ isolated $b$-edges, $\ella-2$ $a$-loops, and $\kb + \ellb + m-1 = \kb + 2\ellb + \ella - 3 = \ka$ isolated $a$-edges. We then attach to each of the open vertices an $a$-loop: the resulting graph $\Gamma$ has combinatorial type $\gamma$.

If $m \ge 1$ and $\ellb = m+1$ or $m+2$ (so that $\ella \le 1$), then $\Delta$ has $m$ $b$-triangles, $m$ $b$-loops, $\kb$ isolated $b$-edges, no $a$-loop, and $2m-1+\kb = \kb + 2\ella + 2\ellb -5 = \ka + \ella -2$ isolated $a$-edges. If $\ella = 1$, we attach to one of the open vertices of $\Delta$ an $a$-loop, and to the other an isolated $a$-edge followed by a $b$-loop. If instead $\ella = 0$, we attach to each open vertex an isolated $a$-edge with a $b$-loop. In either case, this adds $2 - \ella$ isolated $a$-edges and $\ella$ $a$-loops, and the resulting graph $\Gamma$ has combinatorial type $\gamma$.

Finally, if $m = 0$ and $\ella+\ellb = 2$, we let $\Gamma$ be a $\PSL$-reduced path using $\ka$ isolated $a$-edges and $\kb$ isolated $b$-edges, with loops at both ends. Again, $\Gamma$ has combinatorial type $\gamma$.

\begin{figure}[htbp]
\centering
\begin{picture}(60,42)(5,-40)
\gasset{Nw=2,Nh=2}

\node(n10)(0,-2){}
\node(n11)(10.0,-14.0){}
\node(n12)(14,-5.0){}
\node(n13)(18.0,-14.0){}
\drawedge(n11,n12){}
\drawedge(n12,n13){}
\drawedge(n13,n11){}
\drawbpedge[AHnb=0](n10,200,-15,n12,25,-20){}
\drawloop[loopdiam=4,loopangle=150,AHnb=1](n10){}
\node(n21)(26.0,-14.0){}
\node(n22)(30,-5.0){}
\node(n23)(34.0,-14.0){}
\drawedge(n21,n22){}
\drawedge(n22,n23){}
\drawedge(n23,n21){}
\drawloop[loopdiam=4,loopangle=90,AHnb=1](n22){}
\node(n31)(55.0,-14.0){}
\node(n32)(59,-5.0){}
\node(n33)(63.0,-14.0){}
\drawedge(n31,n32){}
\drawedge(n32,n33){}
\drawedge(n33,n31){}
\drawloop[loopdiam=4,loopangle=90,AHnb=0](n32){}
\drawedge[AHnb=0](n13,n21){}
\drawline[AHnb=0](35,-14)(41,-14)
\put(42,-10){$\dots$}
\drawline[AHnb=0](48,-14)(54,-14)
\drawloop[loopdiam=4,loopangle=180](n11){}
\drawloop[loopdiam=4,loopangle=00,AHnb=0](n33){}

\node(n10)(0,-24){}
\node(n11)(10.0,-36.0){}
\node(n12)(14,-27.0){}
\node(n13)(18.0,-36.0){}
\drawedge(n11,n12){}
\drawedge(n12,n13){}
\drawedge(n13,n11){}
\drawbpedge[AHnb=0](n10,200,-15,n12,25,-20){}
\drawloop[loopdiam=4,loopangle=150,AHnb=1](n10){}
\node(n21)(26.0,-36.0){}
\node(n22)(30,-27.0){}
\node(n23)(34.0,-36.0){}
\drawedge(n21,n22){}
\drawedge(n22,n23){}
\drawedge(n23,n21){}
\drawloop[loopdiam=4,loopangle=90,AHnb=1](n22){}
\node(n31)(55.0,-36.0){}
\node(n32)(59,-27.0){}
\node(n33)(63.0,-36.0){}
\drawedge(n31,n32){}
\drawedge(n32,n33){}
\drawedge(n33,n31){}
\drawloop[loopdiam=4,loopangle=90,AHnb=0](n32){}
\drawedge[AHnb=0](n13,n21){}
\drawline[AHnb=0](35,-36)(41,-36)
\put(40,-32){$\dots$}
\drawline[AHnb=0](48,-36)(54,-36)
\drawedge[curvedepth=-4,AHnb=0](n11,n33){}
\end{picture}
\caption{\small The graph $\Gamma$ in typical situations when $m-\ella-\ellb = -2$ or $0$}\label{fig: case -2 and 0}
\end{figure}
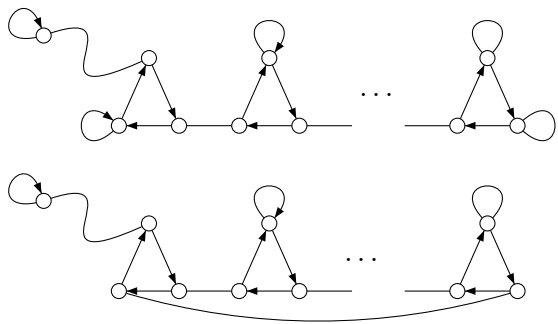

\paragraph*{Case $m-\ella-\ellb = 0$.}
Suppose that $m \ge 1$ and consider the graphs $\Delta_1, \ldots, \Delta_m$ as in the previous case. We then arrange the $\Delta_i$ in a cycle, separating them by isolated $a$-edges. The resulting graph $\Gamma$ has $\kb + \ellb + m = \kb + \ella + 2\ellb = \ka$ isolated $a$-edges, and combinatorial type $\gamma$.

If instead $m = \ella = \ellb = 0$, then we can pick $\Gamma$ to be a cycle labeled by $(ab)^{\kb}$ to achieve combinatorial type $\gamma$.

\paragraph{Case $m-\ella-\ellb \ge 2$.}
Then $m \ge 2 + (\ella + \ellb)$.

Again we consider graphs $\Delta_1, \ldots, \Delta_{\ella+\ellb}$ as before, which use all the $a$- and $b$-loops and all the isolated $b$-edges available, and we connect them in a line with isolated $a$-edges, to produce a graph with 2 open vertices, one at each end of this line. We then connect these 2 open vertices to 2 vertices of a new $b$-triangle, using isolated $a$-edges. The resulting graph $\Delta$ has $\kb + \ellb + (\ella + \ellb - 1) + 2 = \kb + 2\ellb + \ella + 1$ isolated $a$-edges and $\ella + \ellb + 1$ $b$-triangles.

Note that $m - (\ella+\ellb+1)$ is odd and positive. We consider a binary rooted tree T with $m - \ella - \ellb - 1$ vertices (each node has either 2 children or is a leaf): such a tree has necessarily $\frac12(m-\ella-\ellb)$ leaves. We replace each node of $T$ with a $b$-triangle, and each edge in $T$ by an isolated $a$-edge. In particular, the root $b$-triangle is adjacent to 2 $a$-edges, each other internal $b$-triangle of $T$ is adjacent to 3 $a$-edges, and each leaf $b$-triangle is adjacent to a single $a$-edge. On each leaf $b$-triangle, we add an isolated $a$-edge to connect its 2 open vertices. The resulting graph $\Delta'$ uses $m - \ella - \ellb - 1$ $b$-triangles. It also uses $m - \ella - \ellb - 2$ isolated $a$-edges to connect these $b$-triangles, and an additional $\frac12(m-\ella-\ellb)$ $a$-edges on the leaf $b$-triangles. Therefore $\Delta'$ uses $\frac32(m-\ella-\ellb) - 2$ isolated $a$-edges.

Finally we construct $\Gamma$ by connecting the open vertices of $\Delta$ and $\Delta'$ (one each) with an $a$-edge. Then $\Gamma$ has $(\kb + 2\ellb + \ella + 1) + (\frac32(m-\ella-\ellb) - 2) + 1 = \ka$ isolated $a$-edges, and hence $\Gamma$ has combinatorial type $\gamma$, which concludes the proof.
\eop
\begin{figure}[htbp]
\centering
\begin{picture}(60,68)(5,-68)
\gasset{Nw=2,Nh=2}

\node(n10)(0,-2){}
\node(n11)(10.0,-14.0){}
\node(n12)(14,-5.0){}
\node(n13)(18.0,-14.0){}
\drawedge(n11,n12){}
\drawedge(n12,n13){}
\drawedge(n13,n11){}
\drawbpedge[AHnb=0](n10,200,-15,n12,25,-20){}
\drawloop[loopdiam=4,loopangle=150,AHnb=1](n10){}
\node(n21)(26.0,-14.0){}
\node(n22)(30,-5.0){}
\node(n23)(34.0,-14.0){}
\drawedge(n21,n22){}
\drawedge(n22,n23){}
\drawedge(n23,n21){}
\drawloop[loopdiam=4,loopangle=90,AHnb=1](n22){}
\node(n31)(55.0,-14.0){}
\node(n32)(59,-5.0){}
\node(n33)(63.0,-14.0){}
\drawedge(n31,n32){}
\drawedge(n32,n33){}
\drawedge(n33,n31){}
\drawloop[loopdiam=4,loopangle=90,AHnb=0](n32){}
\drawedge[AHnb=0](n13,n21){}
\drawline[AHnb=0](35,-14)(41,-14)
\put(42,-10){$\dots$}
\drawline[AHnb=0](48,-14)(54,-14)
\node(n41)(36.5,-29.0){}
\node(n42)(32.5,-22.0){}
\node(n43)(40.5,-22.0){}
\drawedge(n41,n42){}
\drawedge(n42,n43){}
\drawedge(n43,n41){}
\drawedge[AHnb=0](n11,n42){}
\drawedge[AHnb=0](n33,n43){}
\node(n51)(32.5,-43.0){}
\node(n52)(36.5,-36.0){}
\node(n53)(40.5,-43.0){}
\drawedge(n51,n52){}
\drawedge(n52,n53){}
\drawedge(n53,n51){}
\drawedge[AHnb=0](n41,n52){}
\node(n61)(24.5,-56.0){}
\node(n62)(28.5,-49.0){}
\node(n63)(32.5,-56.0){}
\drawedge(n61,n62){}
\drawedge(n62,n63){}
\drawedge(n63,n61){}
\drawedge[AHnb=0](n51,n62){}
\node(n71)(40.5,-56.0){}
\node(n72)(44.5,-49.0){}
\node(n73)(48.5,-56.0){}
\drawedge(n71,n72){}
\drawedge(n72,n73){}
\drawedge(n73,n71){}
\drawedge[curvedepth=-2,AHnb=0](n71,n73){}
\drawedge[AHnb=0](n53,n72){}
\node(n81)(16.5,-68.0){}
\node(n82)(20.5,-62.0){}
\node(n83)(24.5,-68.0){}
\drawedge(n81,n82){}
\drawedge(n82,n83){}
\drawedge(n83,n81){}
\drawedge[curvedepth=-2,AHnb=0](n81,n83){}
\drawedge[AHnb=0](n61,n82){}
\node(n91)(32.5,-68.0){} 
\node(n92)(36.5,-62.0){} 
\node(n93)(40.5,-68.0){} 
\drawedge(n91,n92){}
\drawedge(n92,n93){}
\drawedge(n93,n91){}
\drawedge[curvedepth=-2,AHnb=0](n91,n93){}
\drawedge[AHnb=0](n63,n92){}
%
{\gasset{Nw=0,Nh=0}\node(n1)(75,0){}\node(n2)(75,-30){}}
{\gasset{ELdistC=y,ELdist=0}\drawedge[ATnb=1,eyo=-1,dash={1}0](n1,n2){\colorbox{white}{$\Delta$}}}
{\gasset{Nw=0,Nh=0}\node(n3)(75,-35){}\node(n4)(75,-69){}}
{\gasset{ELdistC=y,ELdist=0}\drawedge[ATnb=1,eyo=-1,dash={1}0](n3,n4){\colorbox{white}{$\Delta'$}}}
\end{picture}
\caption{\small The graph $\Gamma$ in a typical situation when $m-\ella-\ellb > 0$}\label{fig: case greater than 0}
\end{figure}

\subsubsection{Finite index subgroups}\label{sec: finite index Stallings}

We instantiate a more general result to establish the following statement about finite index subgroups.

\begin{proposition}\label{prop: kb=0}
Let $H$ be a finitely generated subgroup of $\PSL$, whose Stallings graph $\Gamma(H)$ has combinatorial type $(n,\ka,\kb,\ella,\ellb,m)$. Then $H$ has finite index if and only if $\Gamma(H)$ is proper $\PSL$-cyclically reduced and $\kb = 0$. In that case, $H$ has index $n$.
\end{proposition}

\proof
Let $L$ be the language $L$ of shortlex geodesics of $\PSL$. \cite[Thm 6.12]{2017:KharlampovichMiasnikovWeil} states that, if $L$ satisfies the so-called extendability property (see below), then $H$ has finite index if and only if every word of $L$ labels a path in $\Gamma(H)$ starting at the base vertex (which we denote by 1).

We first verify that $L$ satisfies the extendability property. Namely, we show that for every $u\in L$, there exists a sequence of words $(v_n)_n$ such that, for each $n$, $uv_n \in L$ and for almost all $m$, $u$ is a prefix of the shortlex geodesic representative of $uv_nv_m\inv u\inv$. This is indeed the case: if $u$ is of the form $u = u_0a$, then we can let $v_n = b\inv(ab)^n$ (and consider all $m > n$); if $u = u_0b$ or $u_0b\inv$, we let $v_n = ab\inv(ab)^n$.

We now verify that every word of $L$ labels a path in $\Gamma(H)$ starting at vertex 1, if and only if every vertex has an incoming and an outgoing $b$-edge, and is adjacent to an $a$-edge. If the latter property holds, the letters $a$ and $b$ both label permutations of the vertex set of $\Gamma(H)$ and every word can be read along a path from any vertex. Conversely, suppose that every word of $L$ labels a path from 1, and let $p$ be a vertex of $\Gamma(H)$. If $p = 1$, the fact that $a$, $b$ and $b\inv$ are shortlex geodesic guarantees that $p$ is adjacent to an $a$-edge and to incoming and outgoing $b$-edges. If $p \ne 1$, then by definition of $\Gamma(H)$, there exists a word $u_1u_2 \in L$, representing an element of $H$, such that there is a $u_1$-labeled path from 1 to $p$ and a $u_2$-labeled path from $p$ to 1. Either $u_1$ ends with $a$ and $u_2$ starts with $b$ or $b\inv$, or $u_1$ ends with $b^{\pm1}$ and $u_2$ starts with $a$. Without loss of generality, we may assume that $u_1$ ends with $a$ and $u_2$ starts with $b$. Then $p$ is adjacent to an $a$-edge and to an outgoing $b$-edge. Moreover $u_1b\inv$ is also in $L$, so it labels a path in $\Gamma(H)$ starting at 1, and it follows that $p$ also has an incoming $b$-edge.

Thus $H$ has finite index if and only if the letters $a$ and $b$ label permutations of the vertex set of $\Gamma(H)$ and the result follows directly.
\eop

\subsubsection{Isomorphism type and freeness}\label{sec: isomorphism type}

We now turn to the isomorphism type of $H$.  According to Kurosh's classical theorem (\emph{e.g.} \cite[Thm 11.55]{1995:Rotman}), $H$ is isomorphic to a free product of $r_2$ copies of $\Z_2$, $r_3$ copies of $\Z_3$ and a free group of rank $r$. The triple $(r_2,r_3,r)$, called the \emph{isomorphism type} of $H$, characterizes $H$ up to isomorphism.

We first show how to construct from $\Gamma(H)$ a so-called \emph{independent generating set} $\calB$ of $H$, witnessing this particular decomposition: that is, $\calB$ is a generating set of $H$ consists of $r_2$ elements of order 2 generating a subgroup $H_2$ isomorphic to the free product of $r_2$ copies of $\Z_2$, $r_3$ elements of order $3$ generating a subgroup $H_3$ isomorphic to free product of $r_3$ copies of $\Z_3$ and $r$ elements freely generating a subgroup $H_\infty$ of rank $r$, such that $H = H_2 \ast H_3 \ast H_\infty$.

Let $\tilde H$ be the subgroup of $F(\{a,b\})$ whose Stallings graph is $\Gamma(H)$ and let $T$ be a spanning tree of $\Gamma(H)$ containing two edges of every $b$-triangle.
Then $T$ specifies a basis $\tilde\calB$ of $\tilde H$, namely the words of the form $u_p c u\inv_q$, for every edge from vertex $p$ to vertex $q$ with label $c$ which is not in $T$ --- where the word $u_p$ labels a shortest path in $T$ from 1 to $p$. The elements of $\tilde\calB$ (in reduced form) can be partitioned into six subsets $(\tilde\calB_i)_{1\leq i \leq 6}$ depending on their type:
\begin{itemize}
\item[(1)] $\ella$ words of the form $uau\inv$;

\item[(2)] $\ellb$ words of the form $ubu\inv$;

\item[(3)] for every isolated $b$-edge not in $T$, a word of the form $ubv\inv$ with $u\ne v$;

\item[(4)] for every non-loop $a$-edge that is not in $T$ and such that its twin edge is not in $T$ either, two words of the form $uav\inv$ and $vau\inv$, with $u\ne v$;

\item[(5)] for every $a$-edge in $T$, a word of the form $ua^2u\inv$;

\item[(6)] $m$ words of the form $ub^3u\inv$ (where $m$ is the number of $b$-triangles in $\Gamma(H)$).
\end{itemize}
We let $\calB_2$ (resp. $\calB_3$) be the projections in $\PSL$ of the elements of $\tilde\calB$ of type (1) (resp. (2)): each has order 2 (resp. order 3). The projections in $\PSL$ of the elements of $\tilde\calB$ of type (5) and (6) are trivial. The projections of the pairs of elements of $\tilde\calB$ of type (4) are mutually inverse, and we let $\calB_{1,2}$ consist of the projections of an arbitrary choice of one element in each of these pairs. We also let $\calB_{1,3}$ consist of the projections of the elements of $\tilde\calB$ of type (3). Finally we let $\calB = \calB_{1,2} \cup \calB_{1,3} \cup \calB_2 \cup \calB_3$.

\begin{lemma}\label{lm: independent generating set}
The set $\calB$ generates $H$ with the following property: every element of $H$ can be uniquely written as a word on the elements of $\calB$ and their inverses, which is freely reduced and avoids the elements of $\calB_2\inv$ and $\calB_{1,2}\inv$ and the squares of elements of $\calB_2$ or $\calB_3$.

The isomorphism type of $H$ is $(\ella, \ellb, r)$, where $r = |\calB_{1,2}| + |\calB_{1,3}|$.
\end{lemma}

\proof
The generation statement is immediate: $H$ is the projection of $\tilde H$ and, if $x\in \calB_2$ (resp. $x\in \calB_{1,2}$, $x\in \calB_3$), then $x\inv = x$ and $x^2 = 1$ (resp. $x\inv = x$, $x^2 = x\inv$).

To establish the uniqueness statement, first consider a word $u$ on the elements of $\calB$ and their inverses, which is freely reduced and avoids the elements of $\calB_2\inv$ and $\calB_{1,2}\inv$ and the squares of elements of $\calB_2$ or $\calB_3$. By construction, the free reduction $u'$ of the resulting word in $\{a,b,a\inv,b\inv\}^*$ labels a loop at the base vertex in $\Gamma(H)$. By construction also, it does not contain $a^2$, $b^3$ or their inverses as factors.

Let $e$ be the last letter of $u$. If $e \in \calB$, it is derived from an edge of type (1)-(4) and we let $q$ be its end vertex $c\in \{a,b\}$ be its label. If $e\in \calB\inv$, then $e\inv$ is derived from an edge of type (2)-(3), whose start vertex and label we denote by $q$ and $c\inv$ (in this case, $c\inv = b$). In either case, $u' = u''cu_q\inv$. 

Let $v$ be another word on $\calB \cup \calB\inv$, satisfying the same constraints, and representing the same element of $H$ as $u$. Then the free reduction $v'$ of the corresponding word in $\{a,b,a\inv,b\inv\}^*$ factors as $v' = v''du_r\inv$, for some letter $d\in \{a,b,b\inv\}$ and vertex $r$ determined by the the last letter of $v$.

We want to show that $u = v$, by induction on $|u|+|v|$. Then the free reduction $w'$ of $u'{v'}\inv$ represents the trivial element of $\PSL$. In particular, if it is not empty, then it contains $a^2$, $b^3$ or their inverses as a factor, whereas these words do not appear in $u'$ or $v'$.

If $q\ne r$, then the free reduction of $z = u_q\inv u_r$ labels the geodesic path in $T$ from $q$ to $r$, and hence it is non trivial. It follows that $w' = u''czd\inv {v''}\inv \ne 1$. This leads to a contradiction as this word, by construction, does not contain $a^2$, $b^3$ or their inverses.

If $q = r$ and $c\ne d$, then $w' = u''cd\inv{v''}\inv \ne 1$. The only possibility for it to contain $a^2$ or its inverse is if $c = d\inv = a$, a contradiction. If instead it contains $b^3$ or $b^{-3}$, then $c = d\inv \in \{b,b\inv\}$. This implies that $T$ avoids two edges from the same $b$-triangle, again a contradiction.

Therefore $q = r$ and $c = d$, so that $u$ and $v$ have the same last letter (in $\calB\cup\calB\inv$) and we conclude by induction.
\eop

%
%
%

\begin{proposition}\label{prop: compute isomorphism type}
Let $H$ be a finitely generated subgroup of $\PSL$ with combinatorial type $(n,\ka,\kb,\ella,\ellb,m)$, such that $\Gamma(H)$ is $\PSL$-cyclically reduced and $n\ge 2$. Then
\begin{itemize}
\item the isomorphism type of $H$ is $(\ella, \ellb, \frac{n - 2\kb - 3\ella - 4 \ellb}6 + 1)$, and $(\ella, \ellb, \frac{n - 3\ella - 4 \ellb}6 + 1)$ if $H$ has finite index;

\item $H$ is free if and only if $\ella = \ellb = 0$; in that case, $H$ has rank $\frac{n-2\kb}6 +1$, and $\frac n6 + 1$ if $H$ has finite index.
\end{itemize}
If $\Gamma(H)$ is not $\PSL$-cyclically reduced (and $n\ge 2$), then
\begin{itemize}
\item the isomorphism type of $H$ is $(\ella, \ellb, \frac{n - 2\kb - 3\ella - 4 \ellb}6 + \frac13)$ if the base vertex of $\Gamma(H)$ is adjacent to an $a$-edge, and $(\ella, \ellb, \frac{n - 2\kb - 3\ella - 4 \ellb}6 + \frac12)$ if it is adjacent to a $b$-edge;

\item $H$ is free if and only if $\ella = \ellb = 0$, and in that case, $H$ has rank $\frac16 (n-2\kb+2)$ if the base vertex of $\Gamma(H)$ is adjacent to an $a$-edge, and $\frac16 (n-2\kb+3)$ if it is adjacent to a $b$-edge.
\end{itemize}

\end{proposition}

\proof
We first assume that $\Gamma(H)$ is a proper $\PSL$-cyclically reduced labeled graph. Let $T$ be the spanning tree of $\Gamma(H)$ and $\calB$ be the generating set of $H$ constructed above. Let $k'_2$ (resp. $k'_3$) be the number of $a$-edges (resp. $b$-edges outside the $b$-triangles) in $T$. Then $T$ has $k'_2 + k'_3 + 2m$ edges and, since it spans an $n$-vertex graph, we have $k'_2 + k'_3 + 2m = n - 1$.

By Lemma~\ref{lm: independent generating set},  the isomorphism type of $H$ is $(\ella,\ellb,r)$ with
$$r = |\calB_{1,2}| + |\calB_{1,3}| = \ka - \ka' + \kb - \kb' = \ka + \kb - (n-2m-1).$$
%
%
Using the fact that $n = 2\ka+\ella = 2\kb + \ellb + 3m$, we get
\begin{align*}
r &= \ka + \kb + 2m - (2\ka+\ella) + 1 = \kb + 2m -\ka - \ella +1\textrm{ and}\\
r &= \ka + \kb + 2m - (2\kb+\ellb+3m) + 1 = \ka -\kb - \ellb - m +1.
\end{align*}
Adding these equations, we find that
\begin{align*}
2(r-1) &= m-\ella-\ellb = \frac{n-2\kb-\ellb}3-\ella-\ellb\textrm{, and hence}\\
6(r-1) &= n - 2\kb - 3\ella - 4 \ellb,
\end{align*}
as announced.

The statements on the case where $H$ is free follow immediately since freeness corresponds to an isomorphism type of the form $(0,0,r)$. The statements on the case where $H$ has finite index follow from Proposition~\ref{prop: kb=0}.

We now consider the case where $H$ is not $\PSL$-cyclically reduced. Then either the base vertex $v_0$ of $\Gamma(H)$ is adjacent to an $a$-edge but not to a $b$-edge, or it is adjacent to a $b$-edge but not to an $a$-edge. In the first case, we let $\Delta$ be the graph obtained from $\Gamma(H)$ by adding a $b$-loop at $v_0$. Then $\Delta$ is $\PSL$-cyclically reduced, and it is the Stallings graph of a subgroup $K$ of $\PSL$. By construction, the combinatorial type of $K$ is $(n,\ka,\kb,\ella,\ellb+1,m)$ and therefore its isomorphism type is $(\ella,\ellb+1, r$) with $6(r-1) = n - 2\kb - 3\ella - 4 (\ellb+1)$. In particular, $r = \frac{n - 2\kb - 3\ella - 4 \ellb}6 + \frac13$. Moreover, $K$ is isomorphic to $\Z_3 \ast H$, so the isomorphism type of $H$ is $(\ella,\ellb,r)$ and the proposition follows.

If instead the base vertex $v_0$ of $\Gamma(H)$ is adjacent to a $b$-edge but not to an $a$-edge, we let $\Delta$ be the graph obtained from $\Gamma(H)$ by adding an $a$-loop at $v_0$. Again $\Delta$ is $\PSL$-cyclically reduced, and it is the Stallings graph of a subgroup $K$ of $\PSL$. By construction, the combinatorial type of $K$ is $(n,\ka,\kb,\ella+1,\ellb,m)$ and therefore its isomorphism type is $(\ella+1,\ellb, r$) with $6(r-1) = n - 2\kb - 3(\ella+1) - 4 \ellb$. In particular, $r = \frac{n - 2\kb - 3\ella - 4 \ellb}6 + \frac12$. Since $K$ is isomorphic to $\Z_2 \ast H$, the isomorphism type of $H$ is $(\ella,\ellb,r)$ and the proposition follows.
\eop

\begin{remark}
Muller and Schlage-Puchta \cite[Corol. 7]{2004:MullerSchlage-Puchta} showed that a tuple $(\ella,\ellb,r)$ is the isomorphism type of a finite index subgroup of $\PSL$ if and only if $3\ella + 4\ellb + 6(r-1) >0$: Proposition~\ref{prop: compute isomorphism type} above confirms that the value $3\ella + 4\ellb + 6(r-1)$ is then the index of the subgroup. The only impossible isomorphism types, for a finite index subgroup of $\PSL$, are therefore $\Z_2$, $\Z_3$, $\Z_2 \ast \Z_2$ and $\Z$. These are of course possible if we allow arbitrary index subgroups: \textit{e.g.} $\langle a\rangle$, $\langle b\rangle$, $\langle a, bab\inv\rangle$ and $\langle ab\rangle$, whose Stallings graphs have 1 or 2 vertices.
\end{remark}

\section{Counting strategy for finitely generated subgroups of $\PSL$}\label{sec: counting}
We want to count the finitely generated subgroups of $\PSL$, by size. That is: we want to compute for each $n \ge 1$ the number $H_n$ of size $n$ subgroups or, equivalently, the number of size $n$ $\PSL$-reduced rooted graphs.

For this, we use the machinery of analytic combinatorics \cite{2009:FlajoletSedgewick}, which is particularly well suited to count labeled structures.

\subsection{From finitely generated subgroups to $\PSL$-cyclically reduced graphs}\label{sec: counting strategy}

In our situation, a \emph{labeled graph} is one which is equipped with a bijection between its vertex set and a set of the form $\{1,2,\dots,n\}$. 

\begin{figure}[htbp]
\centering
\begin{picture}(30,5)(5,-7)
\gasset{Nw=4,Nh=4}

\node(n0)(7.0,-5.0){1}
\node(n1)(27.0,-5.0){2}

\drawedge[AHnb=0](n0,n1){$a$}
\drawloop[loopangle=180](n0){$b$}
\drawloop[loopangle=0](n1){$b$}
\end{picture}
\caption{\small This graph has only one labeling}\label{fig: unique labeling}
\end{figure}
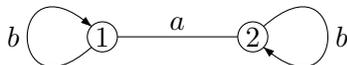

\begin{example}\label{ex: labeled rooted}
The graph $\Gamma$ in Figure~\ref{fig: unique labeling} has only one labeling, but if $v$ is any vertex of $\Gamma$, then the rooted graph $(\Gamma,v)$ has two distinct labelings.
\end{example}

Generalizing Example~\ref{ex: labeled rooted}, we note that an $n$-vertex $\PSL$-cyclically reduced graph may have less than $n!$ distinct labelings (depending on its symmetries, that is, on the size of its automorphism group), but that a rooted $\PSL$-reduced graph $(\Gamma,v)$ always has $n!$ distinct labelings. It follows that, to count $n$-vertex rooted $\PSL$-reduced graphs, it suffices to count the labeled such rooted graphs, and then to divide by $n!$. Similarly, to draw a size $n$ rooted $\PSL$-reduced graph uniformly at random (see Section~\ref{sec: random generation}), it suffices to draw a labeled one, and then to forget the labeling.

The traditional tool when studying the enumeration of discrete structures is the so-called \emph{exponential generating series}, or EGS: if there are $a_n$ size $n$ labeled structures of a certain kind, the corresponding EGS is the formal power series $\sum_n \frac{a_n}{n!} z^n$. If $S(z)$ is a formal power series, we denote by $[z^n]S(z)$ the coefficient of $z^n$ in $S$.

To compute the number $H_n$ of size $n$ subgroups ($n \ge 1$), we will thus compute the EGS $L(z)$ of labeled $\PSL$-reduced rooted graphs
$$L(z) = \sum_n \frac{L_n}{n!}\ z^n,$$
where $L_n$ is the number of labeled $\PSL$-reduced rooted graphs of size $n$. As discussed above, the number of size $n$ subgroups of $\PSL$ is
\begin{equation}\label{eq: Hn}
H_n = \frac{L_n}{n!}.
\end{equation}

Let $(\Gamma,v_0)$ be a labeled rooted $\PSL$-reduced graph. Then we are in exactly one of the following situations.
\begin{itemize}
\item[(1)] $\Gamma$ is $\PSL$-cyclically reduced.
\item[(2)] $\Gamma$ is not $\PSL$-cyclically reduced, and in that case $\Gamma$ has at least 2 vertices, and vertex $v_0$ is adjacent to an $a$-labeled (resp. $b$-labeled) edge and is not adjacent to a $b$-labeled (resp. $a$-labeled) edge.
\end{itemize}
We note that, in situation (2), adding a $b$-loop (resp. an $a$-loop) at $v_0$ yields a $\PSL$-cyclically reduced graph. As a result, the set of structures $(\Gamma, v_0)$ in this situation is in bijection with the structures $(\Delta, f)$ formed by a labeled proper $\PSL$-cyclically reduced graph with a selected loop $f$, where $\Delta$ has the same number of vertices and one more loop than $\Gamma$.

In view of this analysis, we need to study the bivariate EGS of labeled $\PSL$-cyclically reduced graphs,
$$G(z,u) = \sum_{n, \ell} \frac{g_{n,\ell}}{n!}z^n u^\ell,$$
where $g_{n,\ell}$ is the number of labeled proper $\PSL$-cyclically reduced graphs with $n$ vertices and $\ell$ ($a$- or $b$-labeled) loops. Then the number of $n$-vertex labeled rooted $\PSL$-reduced graphs in situation (1) is $n\sum_\ell g_{n,\ell}$, and the number of $n$-vertex labeled rooted $\PSL$-reduced graphs in situation (2) is $\sum_\ell \ell g_{n,\ell}$.

Thus we have (see Example~\ref{ex: size 1 subgroups})
\begin{equation}\label{eq: Ln}
L_1 = 4 \text{ and, for $n\ge 2$, }L_n = \sum_{\ell = 0}^n(n+\ell)g_{n,\ell}.
\end{equation}

\subsection{Counting labeled $\PSL$-cyclically reduced graphs}\label{sec: count cyclically reduced}

We now discuss the EGS $G(z,u)$ of labeled $\PSL$-cyclically reduced graphs. It is actually more convenient to concentrate on the bivariate EGS $\Gpr(z,u)$ of labeled \emph{proper} $\PSL$-cyclically reduced graphs in which, again, the variable $u$ is used to count the number of loops. The graphs we miss with this EGS (those that are not proper) are the 1-vertex graphs with 0 or 1 loop, counted by $z(1+2u)$, see Section~\ref{sec: PSL-reduced graphs} and Figure~\ref{fig: 1 vertex non-proper}, that is,
\begin{equation}\label{eq: G vs Gpr}
G(z,u) = z + 2zu + \Gpr(z,u).
\end{equation}

In a proper $\PSL$-cyclically reduced graph with $n$ vertices ($n\ge 1$), every vertex is either on an $a$-loop, or on an isolated $a$-edge: that is, the $a$-edges determine a partition of the vertex set into 1- and 2-element sets, which we call a $\tau_2$-structure. Similarly every vertex is on a $b$-loop, an isolated $b$-edge or a $b$-triangle, and the $b$-edges determine a partition of the vertex set into 1-, 2- and 3-element sets, together with an orientation for each of the 2- and 3-element sets, which we call a $\tau_3$-structure. Conversely, a $\tau_2$-structure and a $\tau_3$-structure on an $n$-element set determine a graph, which is a proper $\PSL$-cyclically reduced graph if and only if it is connected.

Therefore we introduce now the EGSs  $T_2(z,u)=\sum_{n,\ell} \frac{t^{(2)}_{n,\ell}}{n!}z^n u^\ell$ where $t^{(2)}_{n,\ell}$ is the number of labeled $\tau_2$-structures with $n$ vertices and  $\ell$ $a$-loops and $T_3(z,u, v)\sum_{n,\ell, k} \frac{t^{(3)}_{n,\ell,k}}{n!}z^n u^\ell v^k$ where $t^{(3)}_{n,\ell,k}$ is the number of labeled $\tau_3$-structures with $n$ vertices, $\ell$ $b$-loops and $k$ isolated $b$-edges. We compute these series using the symbolic method on labeled combinatorial structures (see \cite[Thms II.1 and III.2, Example II.13]{2009:FlajoletSedgewick}).

In particular, since a $\tau_2$-structure is a set of 1-element and 2-element sets, corresponding respectively to vertices with an $a$-loop and to pairs of distinct vertices linked with an $a$-edge. The symbolic method takes into account the symmetry on the 2-element sets and yields the following formula:
\begin{equation}\label{eq: T2}
  T_2(z,u) = \exp\left(zu + \frac{z^2}2\right)
\end{equation}
(we consider, for the sake of the series, that there is a single size 0 $\tau_2$-structure).
\begin{remark}\label{rk: tau_2 = involution}
A $\tau_2$-structure on an $n$-element set is, combinatorially, the same thing as an involution on that set.
\end{remark}

Similarly, a $\tau_3$-structure is a set of 1-element sets, pairs and 3-cycles, corresponding respectively to vertices with a $b$-loop, pairs of vertices linked by a directed $b$-edge and directed $b$-triangles. 
Taking into account the symmetry in the triangles, the corresponding EGS is
\begin{equation}\label{eq: T3}
T_3(z,u,v) = \exp\left(zu + z^2v + \frac{z^3}3\right).
\end{equation}
The third variable $v$ counting the isolated $b$-edges  will be used in Sections~\ref{sec: asymptotics expected values} and~\ref{sec: random generation}. Setting $v=1$ amounts to ignoring that count.

Then the bivariate EGS $\widetilde \Gpr(z,u)$ of non-empty labeled graphs whose set of $a$-edges (resp. $b$-edges) is determined by a $\tau_2$-structure (resp. a $\tau_3$-structure) satisfies, for $n \ge 1$ and $0 \le \ell \le n$,
\begin{equation}\label{eq: tilde G}
[z^n u^\ell]\widetilde \Gpr = n! \sum_{\ell_1+\ell_2 = \ell} [z^n u^{\ell_1}]T_2(z,u)\ [z^n u^{\ell_2}]T_3(z,u,1).
\end{equation}
The graphs counted by $\widetilde \Gpr(z,u)$, that may not be connected,  are disjoint unions of proper $\PSL$-cyclically reduced graphs, so that we have $1 + \widetilde \Gpr(z,u) = \exp(\Gpr(z,u))$ or, in other words,
\begin{equation}\label{eq: G}
\Gpr(z,u) = \log(1 + \widetilde \Gpr(z,u)).
\end{equation}
This suffices to compute the coefficients of $G$ and $L$, as well as the numbers $H_n$ of size $n$ subgroups of $\PSL$, see Section~\ref{sec: exact counting and random generation} and Appendix~\ref{app:coef} for more details.

We now turn to the study of the asymptotic behavior of $H_n$.

\section{Combinatorial digression: sets of structures taken from a finite collection}\label{sec: combinatorial digression 1}

The counting strategy outlined in Section~\ref{sec: counting} requires a study of the asymptotic behavior of the coefficients of the series $T_2$, $T_3$, $G$, etc. The EGSs $T_2$ and $T_3$ are particular cases of a general situation, which is the focus of this section. We will return to the specific cases of $T_2$ and $T_3$ in Section~\ref{sec: series T2 and T3}, and we will consider other specific cases in Section~\ref{sec: finite and free}, when we investigate finite index and free subgroups of $\PSL$.

\medskip

Let $\calS$ be a finite combinatorial class with no element of size $0$, and let $\calA$ be the class of labeled structures that are sets of structures in $\calS$. 
The discussion in Section~\ref{sec: count cyclically reduced} makes it clear that if $\calS$ is the class of (labeled) 1- or 2-element sets, then $\calA$ is the class of labeled $\tau_2$-structures. Labeled $\tau_3$-structures can be described in a similar fashion.

Let $A(z)$ and $S(z)$ be the EGSs of $\calA$- and $\calS$-structures. By definition $S(z) = \sum_{i=1}^d s_i z^i$ is a polynomial with non-negative coefficients, and such that $S(0)=0$ and $s_d > 0$. We also assume that $S(z)$ is \emph{aperiodic}, in the sense that it is not of the form $S(z) = T(z^p)$ for some polynomial $T$ and integer $p \ge 2$ (equivalently: $\gcd\{j:s_j\neq 0\}=1$). Using the calculus on labeled combinatorial structures (see \cite[Thms II.1]{2009:FlajoletSedgewick}), we have $A(z)=\exp (S(z))$.

\smallskip

In this section, we investigate the asymptotic behavior of the class $\calA$: we find asymptotic equivalents for the coefficients of $A(z)$, we compute the expected value of the number $\chi_n^{(t)}$ of size $t$ $\calS$-components in a size $n$ $\calA$-structure, we show a large deviation result for $\chi_n^{(t)}$, which is indeed very concentrated around its expected value, and we discuss an efficient algorithm to randomly generate an $\calA$-structure of a given size.

\subsection{Asymptotics of $\calA$-structures}\label{sec: asymptotics A-structures}

Since the series $A(z)$ is the exponential of a polynomial with non-negative coefficients, it is \emph{$H$-admissible} (a technical condition which we do not need to define here, see \cite[Sec. VIII]{2009:FlajoletSedgewick}) and hence amenable to an analysis using saddle point asymptotics. This is done in Proposition~\ref{prop: col polynom} below \cite[Cor VIII.2, p. 568]{2009:FlajoletSedgewick}.

\begin{proposition}\label{prop: col polynom}
Let $S(z) = \sum_{i=1}^d s_i z^i$ be an aperiodic polynomial with non-negative coefficients and let $A(z) = \exp(S(z))$. Then
\[ [z^n] A(z) \sim \frac{e^{S(c_n)}}{\sqrt{2\pi \lambda(c_n)}\, c_n^n},\ \text{ where } \lambda(z) = z^2\frac{d^2}{dz^2} S(z) + z\frac{d}{dz} S(z)\]
where $c_n$ is the least positive solution of the equation $z\frac{d}{dz}S(z)=n$.
\end{proposition}

To use Proposition~\ref{prop: col polynom}, we need to evaluate the asymptotic behaviors of $c_n$ and $c_n^n$. This is done in the Lemma~\ref{lm:dev cn} and Corollary~\ref{cor:cnn} below.
\begin{lemma}\label{lm:dev cn}
Let $S(z)=\sum_{i=1}^d s_i z^i$ be a polynomial of degree $d\geq 1$ with non-negative coefficients. Then the least positive solution $c_n$ of the equation  $z \frac{d}{dz}S(z)=n$ satisfies
\[
c_n = \left(\frac n{d s_d}\right)^{\frac1d} \,\left(1 + R(n^{-\frac1d}) +o\left(n^{-1}\right)\right),
\]
where $R$ is a polynomial of degree at most $d$, depending on $S$ only and such that $R(0) = 0$.
\end{lemma}

\proof
By definition of $c_n$, we have $\sum_{i=1}^d i s_i c_n^i = n$. Since the $s_i$ are non-negative, it follows that
\begin{equation}\label{eq: cn^d}
c_n^d = \frac1{ds_d}\left(n-\sum_{i=1}^{d-1}is_ic_n^i\right) \quad\text{and hence}\quad c_n \leq \left(\frac{n}{d s_d}\right)^\frac1d.
\end{equation}
In particular $c_n = \O(n^{1/d})$, so $\sum_{i=1}^{d-1}\frac{is_i}{ds_d}c_n^i = \O(n^{\frac{d-1}d})$, from which we get the first term in the asymptotic development of $c_n$:
\begin{equation}\label{eq:init cn}
c_n = \left(\frac{n}{ds_d}\right)^{\frac1d}\,\left(1+\O(n^{-\frac{1}d})\right).
\end{equation}
We now prove by induction on $k\geq 1$ that there exists a polynomial $R_k$  of degree at most $k-1$, depending only on $S$ and satisfying $R_k(0) = 1$, such that Property $(I_k)$ below holds:
\begin{equation}\label{eq:Ik}
c_n = \left(\frac{n}{ds_d}\right)^{\frac1d}\,\left(R_k(n^{-\frac1d})+\O(n^{-\frac{k}d})\right)\tag{$I_k$}
\end{equation}
The idea is to bootstrap the development of $c_n$ at rank $k-1$ into the saddle-point equation, gaining one term in the development as the error term goes from $\O(n^{-k/d})$ to $\O(n^{-(k+1)/d})$. More precisely, by Equation~\eqref{eq:init cn}, the property holds for $k=1$. Now assume that it holds for $k\geq 1$. From Equation~\eqref{eq: cn^d} we deduce
\[
c_n = \left(\frac{n}{ds_d}\right)^{\frac1d}\left(1 - \sum_{i=1}^{d-1} is_i \frac{c^i_n}{n} \right)^{\frac1d}.
\]
By induction hypothesis, for $1 \le i \le d-1$, we have
\begin{align*}
c^i_n & =  \left(\frac{n}{ds_d}\right)^{\frac{i}{d}}\,\left(R_k(n^{-\frac1d}) + \O(n^{-\frac{k}{d}})\right)^i \\
& =  \left(\frac{n}{ds_d}\right)^{\frac{i}{d}}\,\left(R^{(i)}_k(n^{-\frac1d}) + \O(n^{-\frac{k}{d}})\right)
\end{align*}
where $R^{(i)}_k$ is a polynomial of degree at most $k-1$ depending only on $S$, satisfying $R^{(i)}_k(0) = 1$. Since $i\leq d-1$, this yields
\[
ds_d\,\frac{c^i_n}{n} = \left(\frac{n}{ds_d}\right)^{\frac{i-d}d}\, R^{(i)}_k(n^{-\frac1d}) + \O(n^{-\frac{k+1}{d}}).
\]
In particular, there exists a polynomial $R^{(0)}_k(z)$ of degree at most $k$, with $R^{(0)}_k(0) = 0$, such that
\begin{equation}\label{eq:sum cj}
\sum_{i=1}^{d-1} is_i\,\frac{c^i_n}{n} = R^{(0)}_k(n^{-\frac1d}) + \O(n^{-\frac{k+1}{d}}).
\end{equation}
The Taylor development of  $x\mapsto (1-x)^{1/d}$ at order $k$ is
\[
(1-x)^{\frac1d} = 1 + \sum_{j=1}^k \frac1{j!}\frac1d\left(\frac1d-1\right)\cdots\left(\frac1d-j+1\right)
(-x)^j +\O(x^{k+1}).\]
Substituting $\sum_{i=1}^{d-1}is_i\frac{c^i_n}{n}$ for $x$ in this formula and applying Equation~\eqref{eq:sum cj}, we find that
\[
\left(1-\sum_{i=1}^{d-1}is_i\,\frac{c^j_n}{n}\right)^{\frac1d} = 1 + R_{k+1}(n^{-\frac1d}) +
\O(n^{-\frac{k+1}{n}}),
\]
for a polynomial $R_{k+1}(z)$ of degree at most $k$ satisfying $R_{k+1}(0)=0$. This concludes the induction and the proof of the lemma.
\eop

\begin{corollary}\label{cor:cnn}
With the same notation as in Lemma~\ref{lm:dev cn}, we have
\[
c_n^n \enspace\sim\enspace  \left(\frac n{d s_d}\right)^{\frac{n}{d}}\, \exp\left( Q(n^{\frac1d})\right),
\]
where $Q$ is a polynomial of degree at most $d-1$, which depends on $S$ only. 
\end{corollary}

\proof
By Lemma~\ref{lm:dev cn}, there is a polynomial $R$ of degree at most $d$, satisfying $R(0) = 0$, such that
\begin{align*}
c_n &=  \left(\frac n{d s_d}\right)^{\frac1d} \, \left(1 + R(n^{-\frac1d})+o\left(n^{-1}\right)\right),\\
c_n^n &= \left(\frac n{d s_d}\right)^{\frac{n}{d}} \, \left(1 + R(n^{-\frac1d})+o\left(n^{-1}\right)\right)^n.
\end{align*}
We have
\[
\left( 1 + R(n^{-\frac1d}) + o\left(n^{-1}\right)\right)^n
= \exp\left(n\log\left( 1 + R(n^{-\frac1d})+ o\left(n^{-1}\right)\right)\right).
\]
Since $\log(1+x)$ admits a Taylor development at $x=0$ at any order, starting with $\log(1+x) = x + \ldots$, there exists a polynomial $R_1$ of degree at most $d$, satisfying $R_1(0) = 0$, such that
\[
\log\left(1 + R(n^{-\frac1d})+o\left(n^{-1}\right)\right)
= R_1(n^{-\frac1d}) + o\left(n^{-1}\right).
\]
If $R_1(z) = \sum_{j=1}^d \beta_jz^j$, we have
\[
\left( 1 + R(n^{-\frac1d}) + o\left(n^{-1}\right)\right)^n
= \exp\left(\sum_{j=1}^d  \beta_j n^{\frac{d-j}d}\right)\left(1+o(1)\right).
\]
Letting $Q(z) = \sum_{j=0}^{d-1} \beta_{d-j}z^j$ yields the announced result.
\eop

We can now apply Proposition~\ref{prop: col polynom} to compute an asymptotic equivalent of $[z^n]A(z)$.

\begin{proposition}\label{pro:equiv exp polynom}
Let $S(z) = \sum_{i=1}^d s_i z^i$ be an aperiodic polynomial with non-negative coefficients and let $A(z) = \exp(S(z))$. Then
\[ [z^n] A(z) \sim \frac1{\sqrt{2\pi d\,n}}\,n^{-\frac{n}{d}}\exp\left(\frac{n}d\, \left(1+\log(ds_d)\right) + T\left(n^{\frac1d}\right)\right),\]
for some polynomial $T(z)$ of degree at most $d-1$.
\end{proposition}

\proof
We use the notation of Proposition~\ref{prop: col polynom}. In particular, $\lambda(z)$ is a polynomial of degree $d$,
with leading coefficient $d(d-1)s_d+ds_d=d^2 s_d$. Hence, by Equation~\eqref{eq:init cn}
\begin{equation}\label{eq: lambda(cn)}
\sqrt{\lambda(c_n)}\sim \sqrt{d^2 s_d  \left((ds_d)^{-\frac1d}\ n^{\frac1d}\right)^{d}} = \sqrt{dn}.
\end{equation}
Note that
\[
\exp\left(S(c_n)\right) = \exp\left(\sum_{i=1}^d s_i c_n^i\right)
= \exp\left(s_d c_n^d + \sum_{i=1}^{d-1} s_i c_n^i\right).
\]
Using Lemma~\ref{lm:dev cn}, we can develop the expression inside the exponential into the sum of a degree $d$ polynomial in $n^{1/d}$ with leading coefficient $s_d (ds_d)^{-d/d} = \frac1d$, and a polynomial in $n^{-1/d}$. In particular, isolating the leading term, there exists a polynomial $P(z)$ of degree at most $d-1$ such that
\[
\exp\left(S(c_n)\right) \sim  \exp\left(\frac{n}d  + P\left(n^{\frac1d}\right)\right).
\]
By Corollary~\ref{cor:cnn}, $c_n^n \sim \left(\frac n{d s_d}\right)^{n/d} \,\exp\left( Q(n^{1/d})\right)$. Now Proposition~\ref{prop: col polynom} states that

$$[z^n]A(z) \sim \frac1{\sqrt{2\pi d}}\,n^{-\frac{n}{d}-\frac1{2}}\, \exp\left(\frac{n}d\, \left(1+\log(ds_d)\right) + P\left(n^{\frac1d}\right) - Q\left(n^{\frac1d}\right)\right),$$
as announced.
\eop

\begin{remark}\label{rk: computable polynomials}
The proofs of Lemma~\ref{lm:dev cn}, Corollary~\ref{cor:cnn} and Proposition~\ref{pro:equiv exp polynom} are constructive, which means that we can explicitly compute the polynomials $R$, $Q$ and $T$ in these statements. This is used in Section~\ref{sec: series T2 and T3} for the EGSs $T_2$ (resp. $T_3$) where  the polynomial $S(z)$ is $z + \frac12z^2$ (resp. $z + z^2 + \frac13z^3$) as well as in Section~\ref{sec: finite and free}.
\end{remark}

\subsection{Expected number of size $t$ components}

In this section we consider the random variable $\chi_n^{(t)}$ ($\chi_n$ if there is no ambiguity on $t$) which counts the number of size $t$ $\calS$-components in a size $n$ $\calA$-structure taken uniformly at random (amongst $\calA$-structures of size $n$). In particular, here and in the next section, we assume that $s_t > 0$, that is, that there exists at least one $\calS$-structure of size $t$.  The bivariate series associated with the random variable $\chi_n^{(t)}$ is $A^{(t)}(z,u)$, or $A(z,u)$ which is, by definition,
\[
A(z,u) = \sum_{\alpha\in\calA}\frac{1}{|\alpha|\,!}z^{|\alpha|} u^{\chi(\alpha)}.
\]
In our case, using the marking technique of the symbolic method (see \cite[Sec. III]{2009:FlajoletSedgewick}), we obtain
\begin{equation}\label{eq: Azu}
A(z,u) = \exp\left(S(z) + s_t\,(u-1) z^t\right).
\end{equation}

\begin{proposition}\label{pro:expected number}
Let $\calS$ be a finite class of combinatorial structures of size at least 1, with EGS $S(z) = \sum_{i=1}^d s_iz^i$, such that $S(z)$ is aperiodic. Let $t > 1$ be an integer such that there is at least one size $t$ $\calS$-structure. Let $\calA$ be the class of combinatorial structures which are sets of $\calS$-structures and let $\chi_n$ be the random variable which counts the number of size $t$ $\calS$-components in a size $n$ $\calA$-structure (for the uniform distribution). Then the expected value of $\chi_n$ satisfies
\[
\E[\chi_n] \sim s_t\,(ds_d)^{-\frac{t}{d}}\,n^{\frac{t}{d}}.
\]
\end{proposition}

\proof
It is well known (see \cite[Prop. III.2]{2009:FlajoletSedgewick}) that the expectation of $\chi_n$ can be computed directly from Equation~(\ref{eq: Azu}), namely from the formula $A(z,u) = \exp\left(S(z) + (u-1)s_t z^t\right)$. More precisely,
\[
\E[\chi_n] = \frac{[z^n] \frac{\partial}{\partial u}A(z,u)\Big|_{u=1}}{[z^n]A(z,1)}
\]
As $A(z,1)=A(z)$, Proposition~\ref{pro:equiv exp polynom} provides an asymptotic equivalent of the denominator: $[z^n]A(z,1)\sim \frac1{\sqrt{2\pi d}}\,n^{-n/d-1/2}\exp\left(\frac nd (1+\log(ds_d)) + T(n^{1/d})\right)$,
where $T(z)$ is a polynomial of degree at most $d-1$.

Moreover, differentiating along variable $u$, we get $\frac{\partial}{\partial u}A(z,u) = s_t z^t A(z,u)$, so that
$\frac{\partial}{\partial u}A(z,u)\Big|_{u=1} = s_t z^t A(z)$ and hence $[z^n]\frac{\partial}{\partial u}A(z,u)\Big|_{u=1} = s_t
[z^{n-t}]A(z)$.  Therefore
\begin{align*}
\E[\chi_n] &\sim s_t\ \frac{(n-t)^{-\frac{n-t}{d}}}{n^{-\frac{n}{d}}} \left(\frac n{n-t}\right)^{\frac12} 
\exp\left(\frac{1+\log(ds_d)}d ((n-t)-n)\right)\ \exp\left(T((n-t)^{\frac1d}) - T(n^{\frac1d})\right) \\
&\sim s_t\ \frac{(n-t)^{-\frac{n-t}{d}}}{n^{-\frac{n}{d}}}
\exp\left(\frac{-t}d\,(1+\log(ds_d))\right)\ \exp\left(T((n-t)^{\frac1d}) - T(n^{\frac1d})\right)
\end{align*}
Note that
\[
(n-t)^{-\frac{n-t}{d}} = n^{-\frac{n-t}{d}}\exp\left(-\frac{n-t}{d} \log\left(1-\frac{t}{n} \right)\right)
\sim n^{-\frac{n-t}{d}} \exp\left(\frac{t}{d}\right).
\]
Also, observe that, if $0 \le \gamma < d$, 
\[
(n-t)^{\frac{\gamma}{d}} - n^{\frac{\gamma}{d}} = n^{\frac{\gamma}{d}}\left(\left(1-\frac{t}{n}\right)^{\frac{\gamma}{d}}-1\right) \sim
-\frac{t\gamma}d n^{\frac{\gamma-d}d},
\]
which vanishes when $n$ tends to infinity. Since $T$ is a polynomial of degree at most $d-1$, it follows that $\lim T\left((n-t)^{1/d}\right)-T(n^{1/d}) = 0$.

We now have
$$\E[\chi_n] \sim s_t\ n^{\frac{t}{d}} \exp\left(-\frac td\log(ds_d)\right) =  s_t\ (ds_d)^{-\frac{t}{d}}\ n^{\frac{t}{d}},$$
as announced.
\eop

\subsection{Large deviations}

In this section we establish a result of large deviations for $\chi_n$. We first establish the following two-variable variant of the saddle-point inequality.

\begin{lemma}\label{lem: large deviation}
Let $X$ be a random variable on a set of discrete structures and let $F(z,u)$ be the corresponding bivariate EGS, $F(z,u) = \sum_{n,k}\frac{f_{n,k}}{n!}\ z^n u^k$, where $f_{n,k}$ is the number of size $n$ structures for which $X = k$. Suppose that $F$ converges in an open domain $D$ containing $(0,0)$, let $(\rho,r) \in D$, with $r,\rho > 0$, and let $n,k\ge 1$. Then
\begin{align*}
\P(X = k) &\enspace\le\enspace \frac{F(\rho,r) \rho^{-n} r^{-k}}{[z^n]F(z,1)}.
\end{align*}
\end{lemma}

\proof
The numbers of size $n$  structures and of size $n$ structures where $X = k$ are, respectively,
$$\sum_{j = 0}^\infty f_{n,j} = n!\ [z^n]F(z,1)\quad\textrm{and}\quad f_{n,k} = n!\ [z^n u^k]F(z,u).$$
Since the coefficients of the series $F(z,u)$ are all non-negative, we have
$$[z^nu^k] F(z,u)\rho^nr^k \le F(\rho,r),$$
for every $(\rho,r) \in D$ with $\rho, r > 0$, and this establishes the announced inequality.
\eop

We can now give the announced large deviations result.

\begin{proposition}\label{pro:large dev}
Let $\calS$, $\calA$, $S(z) = \sum_{i=1}^d s_iz^i$, $t$ and $\chi_n$ be as in Proposition~\ref{pro:expected number}. There exist real numbers $0 < \lambda_0 \le 1 \le \mu_0$ satisfying the following: for every $\lambda, \mu$ such that $0 < \lambda < \lambda_0$ and $\mu > \mu_0$, there exists $0 < \gamma < 1$ such that
\[
\P(\chi_n \leq \lambda\,n^{\frac{t}{d}}) = \O\left(\gamma^{n^{\frac{t}{d}}}\right)\quad\text{and}\quad
\P(\chi_n \geq \mu\,n^{\frac{t}{d}}) = \O\left(\gamma^{n^{\frac{t}{d}}}\right).
\]
\end{proposition}

\proof
We start with the first inequality (the left tail). We apply Lemma~\ref{lem: large deviation} to $A(z,u)$, which converges everywhere, choosing $\rho=c_n$, $0 < r < 1$ to be fixed later, and  $k=r_n-m$ for $r_n=\lfloor r n^{t/d}\rfloor$ and some integer $m\geq 0$. Then
\[
\P\left(\chi_n = r_n-m\right) \enspace\leq\enspace \frac{A(c_n,r)\, c_n^{-n}\, r^{-r_n+m}}{[z^n]A(z)}.
\]
Combining Proposition~\ref{prop: col polynom} and Equations~(\ref{eq: lambda(cn)}) and~(\ref{eq: Azu}) we find that
\[
\frac{A(c_n,r)\, c_n^{-n}}{[z^n]A(z)} \sim \sqrt{2\pi d}\, \sqrt{n}\,\exp\left((r-1)s_t c_n^t\right).
\]
Hence, for $\kappa > \sqrt{2\pi d}$ and $n$ sufficiently large, we have
\[
 \P(\chi_n = r_n -m) \enspace\le\enspace \frac{A(c_n,r) c_n^{-n} r^{-r_n+m}}{[z^n]A(z)}
 \leq\kappa\, \sqrt{n} \,r^m \exp\left((r-1)s_t c_n^t - r n^{\frac{t}{d}} \log r\right).\]
Now, by Lemma~\ref{lm:dev cn}, $(r-1)s_t c_n^t - r n^{t/d} \log r \sim \left((r-1)\,s_t\,(ds_d)^{-t/d} - r\log r\right) n^{t/d}$ (for any fixed value of $r$). Let $f(r) = (r-1)\,s_t\,(ds_d)^{-t/d} - r\log r$. Since $\lim_{r\to 0^+}f(r) = -s_t\,(ds_d)^{-t/d} < 0$, we can choose $0 < \lambda_0 \le 1$ such that $f(\lambda) < 0$ for every $0 < \lambda < \lambda_0$.

At this point, we can fix $r$: for $0 < \lambda < \lambda_0$, we let $\gamma$ be such that $\exp(f(\lambda)) < \gamma < 1$ and we let $r = \lambda$. Then
\begin{align*}
\P(\chi_n = \lfloor \lambda\,n^{\frac{t}{d}}\rfloor-m) &\enspace=\enspace \P(\chi_n = \lfloor rn^{\frac{t}{d}}\rfloor-m) \\
&\enspace\le\enspace \kappa'\, \sqrt n\, r^m\, \exp\left(f(\lambda)\,{n^{\frac{t}{d}}}\right) \enspace\le\enspace \kappa' \lambda^m \gamma^{n^{\frac{t}{d}}}
\end{align*}
for some $\kappa' > 0$ and for $n$ large enough (uniformly in $m$). Summing over $m$ yields
\[
\P(\chi_n \leq  \lambda n^{\frac{t}{d}} ) \leq  \kappa'\,\gamma^{n^{\frac{t}{d}}} \sum_{m=0}^\infty \lambda^m
= \kappa'\,\frac1{1-\lambda}  \gamma^{n^{\frac{t}{d}}}.
\]
This concludes the proof for the left tail.

We proceed in the same way for the right tail, applying Lemma~\ref{lem: large deviation} to $A(z,u)$ with $\rho=c_n$, $r$ to be fixed later, $r_n=\lceil r n^{t/d}\rceil$ and  $k = r_n+m$. With the same computation as before, we find that
\[
\P(\chi_n = \lceil  r n^{\frac{t}{d}}\rceil + m) \leq \kappa'\sqrt{n} \,r^{-m}  \exp\left(f(r)\,n^{\frac{t}{d}}\right).
\]
As $\lim_{r\to +\infty}f(r) = -\infty$, there exists $\mu_0 \ge 1$ such that $f(\mu) < 0$ for every $\mu > \mu_0$. For a fixed $\mu > \mu_0$, we fix $r = \mu$ and we choose $\gamma$ such that $\exp (f(\mu)) < \gamma < 1$. Then, for $n$ sufficiently large,
\[
\P(\chi_n = \lceil  \mu\, n^{\frac{t}{d}}\rceil + m) \leq \kappa'\,\mu^{-m}\,\gamma^{n^{\frac{t}{d}}}.
\]
Again, summing on $m\geq 0$ yields the announced result.
\eop

\begin{remark}\label{rk: optimize bounds}
We can optimize the values of $\lambda_0$ and $\mu_0$ in Proposition~\ref{pro:large dev} by solving the equation $f(r) = 0$, in the intervals $(0,1]$ and $[1,+\infty)$.
\end{remark}

\subsection{Exact counting and random generation of $\calA$-structures}\label{sec: exact counting and random generation}

Let $a_n = n![z^n]A(z)$ be the number of size $n$ elements of $\calA$. Recall that $A(z) = \exp(S(z))$. Differentiating this equation, we get
\begin{equation}\label{eq: pointing A}
\frac{d}{dz} A(z) = \frac{d}{dz} S(z)\ A(z),
\end{equation}
and it follows that the $a_n$ ($n\ge d$) satisfy the recurrence relation below \eqref{eq: rec ai} (say, by considering the coefficients of $z^{n-1}$). We use the \emph{falling factorial} notation $x^{\underline{k}}=x(x-1)\cdots(x-k+1)$:
\begin{equation}\label{eq: rec ai}
a_n  = \sum_{i=1}^d (n-1)^{\underline{i-1}}\ is_i \,a_{n-i}.
\end{equation}
Equation~\eqref{eq: rec ai} can clearly be used to compute the values of $a_n$, once the initial values are obtained by hand. In order to evaluate the complexity of this computation, we consider two models of computation:
\begin{itemize}
\item the \emph{unit-cost model}, where each elementary arithmetic operation (addition, multiplication, \dots) and the uniform random generation of a number in $\{1,\ldots,n\}$ for some $n$, can be performed in constant time;
\item the \emph{bit-cost model}, where we take the number of bits $\|n\|$ of $n$ into account for the complexities. In particular if $m\leq n$, adding $n$ and $m$ costs time $\|n\|$, as does generating a random number in $\{1,\dots,n\}$. This model is more realistic when dealing with numbers such as the $a_n$, which grow exponentially fast. Note that $\|n\|$ is within one unit of $\log_2(n)$, and is therefore equivalent to $\log n$.
\end{itemize}

\begin{proposition}\label{pro:counting}
The first $n$ values of $a_n$ can be computed in $\O(n)$ time in the unit-cost model and $\O(n^2(\log n)^2)$ time in the bit-cost model.
\end{proposition}

\proof
We use Equation~\eqref{eq: rec ai} iteratively. As $d$ is fixed, each step costs constant time in the unit-cost model. 

For the bit-cost model, first observe that since $a_n=n! [z^n]A(z)$, we have
$\log a_n = \Theta(n\log n)$ by Proposition~\ref{pro:equiv exp polynom}. Also, since $(n-1)^{\underline{i-1}}\leq n^{d-1}$, its logarithm is $\Theta(\log n)$, and it can be computed in $\O(n\log n)$ time. So, using the naive multiplication algorithm, we see that computing one term of the sum in Equation~\eqref{eq: rec ai} can be done in time $\O(n(\log n)^2)$. This concludes the proof as there are $n$ values to compute.
\eop

We now turn to the random generation of a size $n$ $\calA$-structure. We rely on the so-called \emph{recursive method}, using the pointing technique: if we point (distinguish) one atom in a non-empty element of $\calA$,  we also distinguish the $\calS$-component it belongs to. Combinatorially, there are $n$ ways of pointing a labelled atom, which translates into the operator $z\frac{d}{dz}$ on the EGS and amounts to Equation~\eqref{eq: pointing A}. Recall that this equation led to the recurrence relation~\eqref{eq: rec ai}.

This in turn translates into the following combinatorial interpretation: a pointed $\calA$-structure can be thought of as a pointed $\calS$-structure together with an $\calA$-structure (this is what Equation~\eqref{eq: pointing A} says). Then the recurrence relation~\eqref{eq: rec ai} can be thought of as partitioning the pointed $\calA$-structures according to the size of the pointed $\calS$-component. More precisely, the probability $p_i$ that this component has size $i$ is 
\[
p_i = \P\left(\text{pointed $\calS$-component has size }i\right)=\frac{(n-1)^{\underline{i-1}}\ is_i \,a_{n-i}}{a_n}.
\]

These considerations lead to the following algorithm: to draw a size $n$ $\calA$-structure uniformly at random, one first draws a random value for $i$ according to the distribution specified by the $p_i$ ($1\le i \le d$). For this value of $i$, one chooses uniformly at random a size $i$ $\calS$-structure (there are exactly $i!s_i$ of them). Then one repeats the procedure to draw a size $n-i$ $\calA$-structure.

The output of this procedure is a sequence $\sigma$ of $\calS$-structures, whose total size is $n$. We then draw uniformly at random a permutation $\tau$ of $\{1,\ldots, n\}$ and label the atoms of $\sigma$ using $\tau$. 

Note that this algorithm supposes a pre-processing phase, namely the computation of the $a_i$ ($i \le n$).

\begin{proposition}\label{pro: complexity random}
The complexity of the random generation algorithm described above is as follows.
\begin{itemize}
\item The pre-processing phase can be performed in $\O(n)$ time in the unit-cost model and $\O(n^2(\log n)^2)$ time in the bit-cost model.
\item Once done, each random generation can be performed in $\O(n)$ time in the unit-cost model and $\O(n^2(\log n)^2)$ time in the bit-cost model.
\end{itemize}
\end{proposition}

\proof
The pre-processing phase consists in computing the $a_i$ ($i \le n$). Its complexity was studied in Proposition~\ref{pro:counting}.

Turning to the random generation phase of the algorithm, we note that each $p_i a_n = (n-1)^{\underline{i-1}}\ is_i \,a_{n-i}$ ($1\le i \le d$) is computed in constant time in the unit-cost model, and in $\O(n\log n)$ in the bit-cost model.

Next, choosing $i\in \{1,\dots,n\}$ according to the distribution given by the $p_j$ can be done by first choosing $m$ uniformly at random in $\{1,\ldots,a_n\}$, then iteratively subtracting $p_ja_n$, for $j=1,\ldots,d$ until this is not possible any more. This is done in constant time in the unit-cost model, and in $\O(n(\log n)^2)$ in the bit-cost model.

Choosing a size $i$ $\calS$-structure uniformly at random is done in constant time in any model (since there are finitely many of them). One then iterates the procedure to draw a size $n-i$ $\calA$-structure: there are at most $n$ iterations.

The last step of the algorithm requires drawing a size $n$ permutation uniformly at random. The Fisher-Yates shuffle algorithm \cite{1969:Knuth} does that in $\O(n)$ time in the unit-cost model and $\O(n\log n)$ time in the bit-cost model.

This concludes the proof.
\eop

\section{Asymptotic behavior of the coefficients of the series $T_2$, $T_3$ and $\widetilde\Gpr$}\label{sec: series T2 and T3}

We now return to the strategy outlined in Section~\ref{sec: counting} to evaluate the asymptotic behavior of the number $H_n$ of size $n$ subgroups of $\PSL$. The first step is to explore the asymptotic behavior of the coefficients of the EGSs $T_2(z,1)$ and $T_3(z,1,1)$ of $\tau_2$- and $\tau_3$-structures. The results in Section~\ref{sec: combinatorial digression 1} apply directly to these series, since $T_2(z,1) = \exp(S_2(z))$ and $T_3(z,1,1) = \exp(S_3(z))$ with $S_2(z) = \frac12 z^2 + z$ and $S_3(z) = \frac13 z^3 + z^2 + z$.

As observed in Remark~\ref{rk: computable polynomials}, the polynomials $R(z)$, $Q(z)$ and $T(z)$ from Lemma~\ref{lm:dev cn}, Corollary~\ref{cor:cnn} and Proposition~\ref{pro:equiv exp polynom} can be computed explicitly. We find the polynomials in Figure~\ref{table: RQT}, which leads to Proposition~\ref{prop: asymptotics T2 T3 Gpr} below.
\bgroup
\def\arraystretch{1.5}
\begin{figure}[htp]
\begin{center}
\begin{tabular}{c|c|c|}
& $S_2 = \frac12\,z^2 + z$ & $S_3 = \frac13z^3+z^2+z$ \\
\hline
$R(z)$ & $\frac18 z^2 - \frac12 z$ & $\frac2{81}z^3 + \frac19 z^2 -\frac23 z$ \\
\hline
$Q(z)$ & $ - \frac12 z$ & $-\frac23 z^2 - \frac19 z$ \\
\hline
$T(z)$ & $z - \frac14$ & $z^2 + \frac13 z - \frac29$\\
\hline
\end{tabular}
\caption{Intermediary polynomials in the computation of asymptotics for $T_2$ and $T_3$}\label{table: RQT}
\end{center}
\end{figure}
\egroup%

\begin{proposition}\label{prop: asymptotics T2 T3 Gpr}
Asymptotic equivalents for $[z^n]T_2$, $[z^n]T_3$ and $[z^n]\widetilde \Gpr$ are as follows
\begin{align}
[z^n]T_2(z,1) &\sim \frac{e^{-\frac14}}{2\sqrt{\pi n}}\ n^{-\frac n2} \exp\left(\frac 12n+n^{\frac12}\right) \label{eq: asymptotics T2} \\
[z^n]T_3(z,1,1) &\sim \frac{e^{-\frac29}}{\sqrt{6\pi n}}\ n^{-\frac n3} \exp\left(\frac 13n + n^{\frac23} + \frac13 n^{\frac13}\right). \label{eq: asymptotics T3} \\
[z^n]\widetilde \Gpr(z,1) &\sim \frac{e^{-\frac{17}{36}}}{\sqrt{12\pi n}}\ n^{\frac n6} \exp\left(-\frac n6 + n^{\frac23} + n^{\frac12} + \frac13 n^{\frac13}\right). \label{eq: widetilde Gpr}
\end{align}
\end{proposition}

\proof
The asymptotic equivalents \eqref{eq: asymptotics T2} and \eqref{eq: asymptotics T3} are direct applications of Proposition~\ref{pro:equiv exp polynom}.
To establish \eqref{eq: widetilde Gpr}, we note that
$$[z^n]\widetilde \Gpr(z,1)= n!\ [z^n]T_2(z,1)\ [z^n]T_3(z,1,1)$$
(see Equation~\eqref{eq: tilde G}) and Stirling's equivalent of $n!$, namely $n! \sim \sqrt{2\pi n}\ e^{-n} n^n$.
\eop

\begin{remark}
The series $T_2(z,1)$, which is also the EGS of involutions (see Remark~\ref{rk: tau_2 = involution}), was studied before, see \cite{1951:ChowlaHersteinMoore,1955:MoserWyman}, and Eq.~(\ref{eq: asymptotics T2}) can already be found in \cite[Prop. VIII.2]{2009:FlajoletSedgewick}. 
\end{remark}

Now let $\ella$ (resp. $\ellb$) be the random variable which counts the number of loops in a $\tau_2$- (resp. $\tau_3$-) structure --- namely the number of size 1 components in a $\tau_2$- (resp. $\tau_3$-) structure, in the vocabulary of Section~\ref{sec: combinatorial digression 1}. Let also $\kb$ be the random variable which counts the number of isolated $b$-edges (that is: of size 2 components) in a $\tau_3$-structure. Then Proposition~\ref{pro:expected number} directly yields the following evaluation of the expected value of these random variables.

\begin{proposition}\label{prop: expected isomorphism type}
The expected value of $\ella$ on the set of size $n$ $\tau_2$-structures is equivalent to $n^{\frac12}$.

The expected value of $\ellb$ (resp. $\kb$) on the set of size $n$ $\tau_3$-structures is equivalent to $n^{\frac13}$ (resp. $n^{\frac23}$).
\end{proposition}

\begin{remark}
The expected value and standard deviation for $\ella$ can also be found in the literature, \textit{e.g.} \cite[Corollary 12]{2004:MullerSchlage-Puchta} and \cite[Proposition IX.19]{2009:FlajoletSedgewick}.
\end{remark}

Proposition~\ref{pro:large dev} leads to a result of strong concentration of these random variables around their expected value.

\begin{theorem}\label{thm: large deviation ella}
For every real numbers $0 < \lambda < 1$ and $1 < \mu$, the following holds:
\begin{itemize}
\item there exists $0 < \gamma < 1$ such that
\[
\P(\ella \leq \lambda\,n^{\frac1{2}}) = \O\left(\gamma^{n^{\frac1{2}}}\right)\quad\text{and}\quad
\P(\ella \geq \mu\,n^{\frac1{2}}) = \O\left(\gamma^{n^{\frac1{2}}}\right);
\]
\item there exists $0 < \gamma < 1$ such that
\[
\P(\ellb \leq \lambda\,n^{\frac1{3}}) = \O\left(\gamma^{n^{\frac1{3}}}\right)\quad\text{and}\quad
\P(\ellb \geq \mu\,n^{\frac1{3}}) = \O\left(\gamma^{n^{\frac1{3}}}\right);
\]
\item there exists $0 < \gamma < 1$ such that
\[
\P(\kb \leq \lambda\,n^{\frac{2}{3}}) = \O\left(\gamma^{n^{\frac{2}{3}}}\right)\quad\text{and}\quad
\P(\kb \geq \mu\,n^{\frac{2}{3}}) = \O\left(\gamma^{n^{\frac{2}{3}}}\right).
\]
\end{itemize}
\end{theorem}

\proof
With the notation in the proof of Proposition~\ref{pro:large dev}, and in the spirit of Remark~\ref{rk: optimize bounds}, we note that, in the particular situation of $\ella, \ellb$, where $t = ds_d = 1$ and $\kb$ where $t = 2$, we have $s_t = 1$ and the function $f(r)$ is equal to $f(r) = r - 1 -r\log r$. It is easily verified that $f$ is increasing on $(0,1)$ and decreasing on $(1,+\infty)$, and that $f(1) = 0$, so that $f(r) < 0$ for all positive $r \ne 1$. Therefore, for all three random variables, Proposition~\ref{pro:large dev} holds for $\lambda_0 = \mu_0 = 1$, which is exactly what our statement expresses.
\eop

We conclude this section by noting that Section~\ref{sec: exact counting and random generation} above gives an algorithm to randomly generate a $\tau_2$- or a $\tau_3$-structure of size $n$, in time $\O(n)$ in the unit-cost model, $\O(n^2(\log n)^2)$ in the bit-cost model.

\section{Combinatorial digression: estimating connectedness}\label{sec: asymptotics col}

Pursuing the counting strategy in Section~\ref{sec: counting}, we now need to handle the question of connectedness (see Section~\ref{sec: count cyclically reduced}, and especially Equation~\eqref{eq: G}). More precisely, given a pair of $\tau_2$- and a $\tau_3$-structure, both of size $n$, we want to estimate how frequently these structures together determine a connected graph over $n$ vertices.

For this purpose, we establish a general technical result, Proposition~\ref{prop: bender revisited} below, which is of independent interest and will be used repeatedly in the sequel of this paper. It is an instantiation of a theorem of Bender \cite{1974:Bender-a} and \cite[Theorem 3]{1974:Bender-b}, see also \cite[p. 162]{1955:MoserWyman} for an earlier, weaker statement. Its application to the counting of $\PSL$-cyclically reduced graphs will be considered in Section~\ref{sec: asymptotics}. 

\begin{proposition}\label{prop: bender revisited}
Let $A(z) = 1 + \sum_{n\ge 1} A_nz^n$ be a formal power series such that
$$A_n \enspace\sim\enspace \alpha n^\beta e^{Q(n)} n^{\epsilon n},$$
with $Q(n) = \sum_{i = 1}^d \gamma_i n^{\delta_i}$, for some non-zero constants $d, \alpha, \beta, \gamma_i, \delta_i, \epsilon$ satisfying $d \ge 1$ and $0 < \delta_d < \cdots < \delta_2 < \delta_1 = 1$.
Then $\frac{A_{n-1}}{A_n} \sim e^{-\gamma_1-\epsilon}n^{-\epsilon}$.

If, in addition, $\alpha,\epsilon > 0$ and $B(z)$ is such that $B(z) = \log(A(z))$, then for every $s \ge 1$, we have
$$[z^n]B(z) =  \sum_{k = 0}^{s-1} b_kA_{n-k} + \O(A_{n-s}),$$
where $b_0 = 1$ and for each $k\ge 1$
$b_k = \sum_{\ell=1}^k (-1)^\ell \sum_{j_1+\ldots+j_\ell =k}A_{j_1}\cdots A_{j_\ell}.$
\end{proposition}

\begin{remark}\label{rk: compute the bk}
The first coefficients $b_k$ in Proposition~\ref{prop: bender revisited} are $b_1 = -A_1$, $b_2 = A_1^2 - A_2$ and $b_3 = -A_1^3 + 2A_1A_2 - A_3$.
\end{remark}

\proof
Bender's theorem \cite{1974:Bender-a} states that, if $s_0 \ge 1$ and the two following conditions
\begin{itemize}
\item[\textit{(i)}] $A_{n-1} = o(A_n)$;

\item[\textit{(ii)}] for every $s \geq s_0$, $\sum_{k = s}^{n-s}A_k\ A_{n-k} = \O(A_{n-s})$
\end{itemize}
hold, then, for every $s \ge s_0$,
$$[z^n]B(z) = \sum_{k = 0}^{s-1} b_k\ A_{n-k} + \O(A_{n-s}),$$
where $b_k = [z^k]\frac1{A(z)}$, and therefore satisfies the formula  given in the statement of the proposition. 

Let us start with Condition \textit{(i)}. We have
$$\frac{A_{n-1}}{A_n} = \frac{(n-1)^\beta}{n^\beta} \ \exp(Q(n-1) - Q(n))\ \frac{(n-1)^{\epsilon(n-1)}}{n^{\epsilon n}}.$$
Since $Q(n-1) - Q(n) = \sum_{i = 1}^d \gamma_i \left((n - 1)^{\delta_i} - n^{\delta_i}\right) 
= \sum_{i = 1}^d \gamma_i n^{\delta_i} \left(\left(1 - \frac1n\right)^{\delta_i} - 1\right)$
and  $\delta_1 = 1$, this is equal to $-\gamma_1 + \O(n^{\delta_2-1})$ (or simply $-\gamma_1$ if $d = 1$).
Moreover
$$\frac{(n-1)^{\epsilon(n-1)}}{n^{\epsilon n}} = \frac1{(n-1)^\epsilon}\ \left(1 - \frac1n\right)^{\epsilon n} \sim e^{-\epsilon} n^{-\epsilon}.$$
Therefore $\frac{A_{n-1}}{A_n} \sim e^{-\gamma_1-\epsilon}n^{-\epsilon}$ tends to 0, and Condition \textit{(i)} is satisfied.

We now turn to Condition \textit{(ii)}. Let $\bar A_k = \alpha k^\beta e^{Q(k)} k^{\epsilon k}$.
For integers $n,s,k$ such that $n > 2s$ and $s \le k \le n-s$, we let
$$S_k = \frac{A_kA_{n-k}}{A_{n-s}}\quad\text{and}\quad \bar S_k = \frac{\bar A_k\bar A_{n-k}}{\bar A_{n-s}}.$$

Since $A_k \sim \bar A_k$, there exists a value $s_0$ such that, for $k \ge s_0$, we have $\frac12 \bar A_k \le A_k \le 2 \bar A_k$. Therefore, for $s \geq s_0$, $n > 2s$ and $s \le k \le n-s$, we have $S_k \le 8 \bar S_k$. In particular, in order to show that $\sum_{k=s}^{n-s} S_k$ is bounded (this is Condition \textit{(ii)}), we only need to show that $\sum_{k=s}^{n-s} \bar S_k$ is bounded.

By symmetry, it suffices to show that $\sum_{k=s}^{\frac n2} \bar S_k$ is bounded. We first show (part (a) below) that, for a well-chosen $\zeta \in (0,\frac12)$, the sequence $(\bar S_k)_{k = s}^{\zeta n}$ is bounded above by a geometric series with ratio less than 1, so its sum is bounded, and then (part (b)) we show that $\lim_n\sum_{k = \zeta n}^{\frac n2} \bar S_k = 0$.

\paragraph*{part (a).}
By definition,
\begin{align*}
\bar S_k& = \alpha \underbrace{\left(\frac{k(n-k)}{n-s}\right)^\beta}_{P_1} \ \underbrace{\exp(Q(k) + Q(n-k) - Q(n-s))}_{P_2}\ \underbrace{\Bigg(\frac{k^{k} (n-k)^{(n-k)}}{(n-s)^{(n-s)}}\Bigg)^{\epsilon}}_{P_3},
\end{align*}
so we have
\begin{align*}
\frac{\bar S_{k+1}}{\bar S_k} &= \underbrace{\left(\frac{(k+1)(n-k-1)}{k(n-k)}\right)^\beta}_{P_1'}  \underbrace{\exp(Q(k+1)-Q(k) + Q(n-k-1) - Q(n-k))}_{P_2'} \\
&\hskip 8cm \underbrace{\left(\frac{(k+1)^{k+1} (n-k-1)^{n-k-1}}{k^{k}(n-k)^{n-k}}\right)^{\epsilon}}_{P_3'}.
\end{align*}

For $s\leq k \leq \frac{n}2$, it holds
$1 \leq \left(1+ \frac{1}{k}\right)\ \left(1- \frac{1}{n-k}\right) \leq 1+\frac1s$. Therefore $P_1' \leq 1$ if $\beta \leq 0$ and $P_1'\leq (1+\frac1s)^\beta$ if $\beta >0$.

Further, we have
\begin{align*}
\log P_2'  &= \sum_{i=1}^d \gamma_i \left((k+1)^{\delta_i}- k^{\delta_i} + (n-k-1)^{\delta_i}- (n-k)^{\delta_i} \right) \\
&= \sum_{i=2}^d \gamma_i \left((k+1)^{\delta_i}- k^{\delta_i} + (n-k-1)^{\delta_i}- (n-k)^{\delta_i} \right) \text{ since $d_1 = 1$}\\
&= \sum_{i=2}^d \gamma_i \left(k^{\delta_i}\left(\left(1 + \frac1k\right)^{\delta_i} - 1\right) + (n-k)^{\delta_i}\left(\left(1 - \frac1{n-k}\right)^{\delta_i} - 1\right)\right).
\end{align*}
For each $i \ge 2$, we have $0 < \delta_i < 1$, so
$$\left(1 + \frac1k\right)^{\delta_i} \le 1 + \delta_i\ k\inv \quad\text{ and}\quad\left(1 - \frac1{n-k}\right)^{\delta_i} \le 1 - \delta_i\ (n-k)\inv.$$
It follows that the coefficient of $\gamma_i$ in $\log P'_2$ is bounded above by
$$C_i \enspace=\enspace \delta_i\ (k^{\delta_i-1} - (n-k)^{\delta_i-1}) \enspace\le\enspace \delta_i,$$
where the inequality follows from the fact that $C_i$  is a decreasing function of $k$ and  $1\leq s\leq k$. It follows that each summand in $\log P'_2$ is bounded above by $|\gamma_i|\delta_i$, and hence $P'_2 < C$ with $C=\exp (\sum_{i=2}^d|\gamma_i|\delta_i)$.

Turning to $P'_3$, we find that
\begin{align*}
\frac{\log P'_3}\epsilon &= (k+1)\log(k+1) - k\log k + (n-k-1)\log(n-k-1) - (n-k)\log(n-k) \\
&= k\log\left(1+\frac1k\right) + (n-k)\log\left(1-\frac1{n-k}\right) + \log\left(\frac{k+1}{n-k-1}\right).
\end{align*}
Recall that $\log(1+x) \le x$ for all $x > -1$, so the sum of the first two summands above is negative. Now fix $0 < \zeta < \frac12$ and suppose that $k \le \zeta n$. Then
$\frac{k+1}{n-k-1} \le \frac\zeta{1-\zeta} +\O(n\inv)$ so that
\begin{align*}
\frac{\log P'_3}\epsilon &\enspace\le\enspace \log\left(\frac\zeta{1-\zeta} +\O(n\inv)\right) 
\enspace\le\enspace \log\left(\frac\zeta{1-\zeta}\right) +\O(n\inv)
\end{align*}
Hence, $ P'_3 \le  \left(\frac\zeta{1-\zeta}\right)^\epsilon\ \left(1 + \O(n\inv)\right)$.
Therefore $\frac{\bar S_{k+1}}{\bar S_k} \le C\ \left(1+\frac1s\right)^{|\beta|} \left(\frac\zeta{1-\zeta}\right)^\epsilon\ \left(1 + \O(n\inv)\right)$.
On the interval $(0,\frac12)$, the function $\zeta\mapsto \frac\zeta{1-\zeta}$ increases from $0$ to $1$, so the parameter $\zeta$ can be chosen in such a way that $C \left(1+\frac1s\right)^{|\beta|} \left(\frac\zeta{1-\zeta}\right)^\epsilon \le \frac12$. With this value of $\zeta$, we have $\frac{\bar S_{k+1}}{\bar S_k} < 1$ for $n$ large enough, as announced.

\paragraph*{part (b).}
We now want to show that $\lim_n \sum_{k=\zeta n}^{\frac n2} \bar S_k = 0$. We first bound each $\bar S_k$, when $\zeta n \leq k  \leq \frac{n}{2}$.

Since $\zeta(1-\zeta)n^2 \le k(n-k) \le \frac{n^2}4$ and $s < \frac n2$, we have $0 < \frac{k(n-k)}{n-s} \le \frac1{2n} \le \frac12$, and hence $P_1$ is bounded above by $2^{-\beta}$.

Since $d_1=1$ and $k \ge s$, we have
$$P_2= \sum_{i=1}^d \gamma_i  \left(k^{\delta_i} + (n-k)^{\delta_i} - (n-s)^{\delta_i}\right) \leq \gamma_1 s+ \sum_{i=2}^d |\gamma_i|  k^{\delta_i} \leq |\gamma_1| s+ \sum_{i=2}^d |\gamma_i|  \left(\frac{n}{2}\right)^{\delta_i},$$
that is, $P_2$ is bounded by a constant if $d = 1$, and is $\O(n^{\delta_2})$ if $d\ge 2$.
Finally,
\begin{align*}
P_3 &= \left(\frac{k^k\ (n-k)^{n-k}}{(n-s)^{n-s}}\right)^\epsilon = \left(\frac{\left(\frac kn\right)^{\frac kn}\ \left(1-\frac kn\right)^{1-\frac kn}}{\left(1-\frac sn\right)^{1-\frac sn}}\ n^{-s}\right)^{\epsilon n} \\
\frac{\log P_3}{\epsilon n} &= \frac kn \log\left(\frac kn\right) + \left(1-\frac kn\right)\log\left(1 - \frac kn\right) - \left(1-\frac sn\right)\log\left(1 - \frac sn\right)  - s \log n.
\end{align*}
Note that the function $x \mapsto x\log x + (1-x)\log(1-x)$ decreases on interval $(0,\frac12)$, from $0$ to $-\log 2$. Therefore
$$\frac{\log P_3}{\epsilon n} \le \zeta\log \zeta + (1-\zeta)\log(1-\zeta) -s \log n + \O(n\inv).$$
In particular, for $n$ large enough and for some $s' > 0$, we have $\log P_3 \le -\epsilon s' n \log n$ and $P_3 \le n^{-\epsilon s'n}$. It follows that, for some $0 < \eta < 1$, $n$ large enough and $\zeta n \le k \le \frac n2$, we have $\bar S_k = P_1P_2P_3 \le \eta^n$. This ensures that
$\lim_n \sum_{k=\zeta n}^{\frac n2} \bar S_k = 0$
as expected.

Summarizing, we have shown that there exists $s_0$ such that  for $s\geq s_0, \sum_{k = s}^{n-s}A_k\ A_{n-k} = \O(A_{n-s}).$
It remains to extend this property to any $s\geq 1$. For $1\leq s < s_0$,
$$\sum_{k = s}^{n-s}A_k\ A_{n-k} = \O(A_{n-s_0}) + 2\sum_{k = s}^{s_0-1}A_k\ A_{n-k} =\O(A_{n-s}),$$ since it is the sum of a finite number of terms which are all $\O(A_{n-s})$ in view of Condition~\textit{(i)}.

 This concludes the verification that Condition~\textit{(ii)} holds for any $s\geq 1$, and hence the proof of Proposition~\ref{prop: bender revisited}.
\eop

\section{Asymptotics of finitely generated subgroups of $\PSL$}\label{sec: asymptotics}

We again return to finitely generated subgroups of $\PSL$ and, first, to labeled $\PSL$-cyclically reduced graphs, as explained in Section~\ref{sec: counting}.

\subsection{Asymptotic behavior of the number of finitely generated subgroups of $\PSL$}\label{sec: number of subgroups}

With the notation of Section~\ref{sec: counting}, the EGS $\Gpr$ of proper labeled $\PSL$-cyclically reduced graphs satisfies $1+\widetilde\Gpr = \exp(\Gpr)$ (Equation~\eqref{eq: G}). Recall also that we found an asymptotic equivalent of the coefficients of $\widetilde\Gpr$ (Equation~\eqref{eq: widetilde Gpr}), namely
$$[z^n]\widetilde\Gpr \enspace\sim\enspace \frac{e^{-\frac{17}{36}}}{\sqrt{12\pi n}}\ n^{\frac n6} \exp\left(-\frac n6 + n^{\frac23} + n^{\frac12} + \frac13 n^{\frac13}\right).$$
Proposition~\ref{prop: bender revisited} then leads to the following result.

\begin{theorem}\label{thm: probability of being connected}
The probability $p_n$ that an $n$-vertex labeled graph whose set of $a$-edges (resp. $b$-edges) is determined by a $\tau_2$-structure (resp. a $\tau_3$-structure) is a proper $\PSL$-cyclically reduced graph, satisfies
$$p_n = \frac{[z^n]\Gpr(z,1)}{[z_n]\widetilde\Gpr(z,1)} = 1 - n^{-\frac16} + \O(n^{-\frac13}).$$
\end{theorem}

\proof
As mentioned above, we have $1+\widetilde\Gpr = \exp(\Gpr)$. Moreover $\widetilde\Gpr$ satisfies the hypothesis of Proposition~\ref{prop: bender revisited} with $\epsilon = -\gamma_1 = \frac16$, $b_1 = -[z]\widetilde\Gpr(z,1) = -1$ (see Remark~\ref{rk: compute the bk}). Then, by Proposition~\ref{prop: bender revisited}, we obtain
\begin{align*}
[z^n]\Gpr(z,1) &= [z^n]\widetilde\Gpr(z,1) - [z^{n-1}]\widetilde\Gpr(z,1) + \O\left([z^{n-2}]\widetilde\Gpr(z,1)\right) \\
p_n &= 1 - \frac{[z^{n-1}]\widetilde\Gpr(z,1)}{[z^n]\widetilde\Gpr(z,1)} + \O\left(\frac{[z^{n-2}]\widetilde\Gpr(z,1)}{[z^n]\widetilde\Gpr(z,1)}\right).
\end{align*}
By Proposition~\ref{prop: bender revisited} again, $\frac{[z^{n-1}]\widetilde\Gpr(z,1)}{[z^n]\widetilde\Gpr(z,1)} \sim e^{-\gamma_1-\epsilon}n^{-\epsilon} = n^{-\frac16}$.

Moreover $\frac{[z^{n-2}]\widetilde\Gpr(z,1)}{[z^n]\widetilde\Gpr(z,1)} = \frac{[z^{n-2}]\widetilde\Gpr(z,1)}{[z^{n-1}]\widetilde\Gpr(z,1)}\ \frac{[z^{n-1}]\widetilde\Gpr(z,1)}{[z^n]\widetilde\Gpr(z,1)} \sim n^{-\frac13}$, and this concludes the proof.
\eop

We can now derive information about the asymptotic behavior of the coefficients of $\Gpr$, $G$ and $L$, and of the $H_n$.

\begin{corollary}\label{cor: asymptotics for G}
The EGSs $G$ and $\Gpr$, respectively of labeled $\PSL$-cyclically reduced graphs and labeled proper $\PSL$-cyclically reduced graphs, satisfy the following:
$$[z^n]G(z,1) \sim [z^n]\Gpr(z,1) \sim \frac{e^{-\frac{17}{36}}}{\sqrt{12\pi n}}\ n^{\frac n6} \exp\left(-\frac n6 + n^{\frac23} + n^{\frac12} + \frac13 n^{\frac13}\right).$$
\end{corollary}

\proof
Theorem~\ref{thm: probability of being connected} implies that $[z^n]\Gpr(z,1) \sim [z^n]\widetilde\Gpr(z,1)$. Moreover $G$ and $\Gpr$ differ only in their $z$ and $zu$ coefficients. The result then follows from the computation of an asymptotic equivalent for $[z^n]\widetilde\Gpr(z,1)$, Equation~\eqref{eq: widetilde Gpr}. 
\eop

Finally, we find an asymptotic estimate of the number of size $n$ subgroups of $\PSL$.

\begin{theorem}\label{thm: asymptotics for H}
The number $H_n$ of size $n$ subgroups of $\PSL$ satisfies the following:
$$H_n = \Theta\left(n^{\frac n6+\frac12} \exp\left(-\frac n6 + n^{\frac23} + n^{\frac12} + \frac13 n^{\frac13}\right)\right).$$
\end{theorem}

\proof
In view of Eq.s~\eqref{eq: Hn} and~\eqref{eq: Ln}, we have $V_n < H_n < 2V_n$, where $V_n = n\ [z^n]G(z,1)$. The result then follows from Corollary~\ref{cor: asymptotics for G}.
\eop

\subsection{Asymptotic behavior of the isomorphism type}\label{sec: asymptotics expected values}

Let $(\ella, \ellb, r)$ be the isomorphism type of a size $n$ $\PSL$-cyclically reduced subgroup (see Section~\ref{sec: isomorphism type}), seen as a triple of random variables.

\begin{theorem}\label{thm: expected isomorphism type}
On the set of size $n$ $\PSL$-cyclically reduced subgroups, the expected value of $\ella$ (resp. $\ellb$, $r$) is $n^{\frac12} + o(n^{\frac12})$ (resp. $n^{\frac13} + o(n^{\frac13})$, $\frac n6 - \frac13n^{\frac23} + o(n^{\frac23})$).
\end{theorem}

\proof
We already computed the asymptotic behavior of the expected value of the random variables $\ella$, $\ellb$, $\kb$, on the sets of size $n$ $\tau_2$- and $\tau_3$-structures, respectively (Section~\ref{sec: series T2 and T3}).

On the set of size $n$ $\PSL$-cyclically reduced subgroups, the expected value of $\ella$ is the sum $\sum_{k= 0}^n k\,p_k$, where $p_k$ is the probability that a size $n$ $\PSL$-cyclically reduced subgroup satisfies $\ella = k$. Let $B_n$ be the number of $\PSL$-cyclically reduced subgroups of size $n$ and let $A_{n,k}$ be the number of those subgroups such that $\ella = k$. Let also $\tilde B_n$ be the number of pairs of a $\tau_2$-structure and a $\tau_3$-structure, both of size $n$, and $\tilde A_{n,k}$ be the number of such pairs such that $\ella = k$. Then
$$p_k = \frac{A_{n,k}}{B_n} \le \frac{\tilde A_{n,k}}{B_n} = \frac{\tilde A_{n,k}}{\tilde B_n}\ \frac{\tilde B_n}{B_n}.$$
Since $\frac{\tilde B_n}{B_n} \sim 1$ (Theorem~\ref{thm: probability of being connected}) and $\frac{\tilde A_{n,k}}{\tilde B_n} \sim n^{\frac12}$ (since this ratio is the expected value of $\ella$ on $\tau_2$-structures), this upper bound is asymptotically equivalent to $n^{\frac12}$.

Now let $C_n = \tilde B_n - B_n$ be the number of pairs of a $\tau_2$- and a $\tau_3$-structure of size $n$ that do not define a Stallings graph (irrespective of the value of $\ella$). Then
$$p_k = \frac{A_{n,k}}{B_n} \ge \frac{\tilde A_{n,k}-C_n}{B_n} = \frac{\tilde A_{n,k}}{\tilde B_n}\ \frac{\tilde B_n}{B_n} - \frac{C_n}{\tilde B_n}\ \frac{\tilde B_n}{B_n}.$$
We just saw that the first term in this difference is equivalent to $n^{\frac12}$. We also know (Theorem~\ref{thm: probability of being connected} again) that the second term is equivalent to $n^{-\frac16}$, so this lower bound is also asymptotically equivalent to $n^{\frac12}$, thus completing the proof of the statement concerning $\ella$.

The same reasoning can be applied to the random variables $\ellb$ and $\kb$ on the set of size $n$ $\PSL$-cyclically reduced subgroups, with $\E(\kb) = n^{\frac23} + o(n^{\frac23})$.

By Proposition~\ref{prop: compute isomorphism type}, the random variable $r$ is equal to $1 + \frac16(n - 2\kb - 3\ella - 4\ellb)$. It follows directly that
$$\E(r) = 1 + \frac16(n - 3\E(\ella) - 2\E(\kb) - 4\E(\ellb)) = \frac16n - \frac13 n^\frac23 + o(n^\frac23),$$
concluding the proof.
\eop

With the same reasoning, we show that the random variables $\ella$, $\ellb$, $\kb$ satisfy the same large deviations result on the set of labeled $\PSL$-cyclically reduced graphs as in Theorem~\ref{thm: large deviation ella}. Together with Proposition~\ref{prop: compute isomorphism type}, this translates into a large deviations statement for the random variable $r$.

\begin{corollary}\label{cor: large deviation rank}
For every $0 < \lambda < 1 < \mu$, there exists $0 < \gamma < 1$ such that
$$\P\left(r < \frac 16n - \frac\lambda3 n^{\frac23}\right) = \O\left(\gamma^{n^{\frac23}}\right) \enspace\textrm{and}\enspace \P\left(r > \frac 16n - \frac\mu3 n^{\frac23}\right) = \O\left(\gamma^{n^{\frac23}}\right).$$
\end{corollary}

\subsection{Random generation of finitely generated subgroups of $\PSL$}\label{sec: random generation}

Let $n \ge 2$. The algorithm to randomly generate a finitely generated size $n$ subgroup of $\PSL$, or more precisely, the Stallings graph of such a subgroup, is as follows:
\begin{itemize}
\item draw a labeled size $n$ proper $\PSL$-cyclically reduced graph $G$ uniformly at random;
\item if $G$ has $\ell$ loops, draw uniformly at random an integer $i \in [1,n+\ell]$;
\item if $i \le n$, return the graph $G$ with $i$ as the base vertex and forget the labeling of vertices; if $i > n$, delete from $G$ the $(i-n)$-th loop (in the order given by vertex labeling), say, at vertex $j$, and return the graph $G$ with $j$ as the base vertex and with all vertex labels forgotten.
\end{itemize}

We now discuss the first step, namely drawing a labeled size $n$ proper $\PSL$-cyclically reduced graph $G$ uniformly at random. This is done by drawing independently and uniformly at random, a $\tau_2$-structure and a $\tau_3$-structure on the set $[1,n]$. Of course, the labeled graph thus specified may not be connected, and therefore may not qualify as a $\PSL$-cyclically reduced graph, but Theorem~\ref{thm: probability of being connected} shows that this happens with vanishing probability. This justifies a \emph{rejection algorithm}:
\begin{itemize}
\item draw independently and uniformly at random a pair of a $\tau_2$- and a $\tau_3$-structure, both of size $n$;
\item return the corresponding labeled graph if it is connected, and repeat the previous step if it is not.
\end{itemize}
Theorem~\ref{thm: probability of being connected} states that with probability $1 - n^{-\frac16} + \O(n^{-\frac13})$, a single draw will be sufficient. In other words, the probability of needing at least $k+1$ draws is equivalent to $n^{-\frac k6}$.

The random generation of a size $n$ $\tau_2$- or $\tau_3$-structure was discussed at the end of Section~\ref{sec: series T2 and T3} and in Section~\ref{sec: exact counting and random generation}. In particular (Proposition~\ref{pro: complexity random}), it can be done in time $\O(n)$ in the unit-cost model, and in time $\O(n^2\log^2 n)$ in the bit-cost model.

\section{Finite index and free subgroups}\label{sec: finite and free}

The same techniques as in Sections~\ref{sec: counting} and~\ref{sec: asymptotics} can be used to count, randomly generate or asymptotically evaluate the number of finite index subgroups, free subgroups and free finite index subgroups of a given size. The counting and asymptotic results were already known in the finite index case, see Newman \cite{1976:Newman}. 

\subsection{Finite index subgroups}\label{sec: finite index subgroups}

\paragraph*{Asymptotic enumeration.} By Proposition~\ref{prop: kb=0}, a
subgroup of $\PSL$ has index $n$ if and only if its Stallings graph $\Gamma$ is $\PSL$-cyclically reduced, has $n$ vertices and has no isolated $b$-edge. That is, $\Gamma$ satisfies $\kb = 0$.

Let $H^\findex_n$ and $L^\findex_n$ be the numbers of index $n$ subgroups and labeled Stallings graphs of index $n$ subgroups, respectively. Let also $G^\findex$ be the EGS of labeled $\PSL$-cyclically reduced graphs with $\kb = 0$. The discussion in Sections~\ref{sec: counting strategy} and~\ref{sec: count cyclically reduced} can be reproduced identically, except for the following changes:
\begin{itemize}
\item since every finite index subgroup has a $\PSL$-cyclically reduced Stallings graph, the formula for $L_n^\findex$ is simply $L^\findex_n = n\ n!\ [z^n]G^\findex$ and $H^\findex_n = n[z^n]G^\findex$;
\item the EGS $T_3$ must be replaced by $T^\findex_3(z) = \exp(z + \frac{z^3}3)$, which counts the \emph{permutational} $\tau_3$-structures (that is, the permutations of order 3).
\end{itemize}
If $\widetilde G^\findex(z)$ is the EGS of pairs of a $\tau_2$-structure and a permutational $\tau_3$-structure over $n$ elements, then we have, as in the general case,
\begin{equation}\label{eq: tilde Gfi}
[z^n]\widetilde G^\findex = n!\ [z^n]T_2 \ [z^n]T_3^\findex\quad\text{and}\quad1 + \widetilde G^\findex = \exp(G^\findex).
\end{equation}

We can then turn to Section~\ref{sec: combinatorial digression 1} for exact computations (see Appendix~\ref{app:coef}) and for asymptotic estimates. In particular, we find the intermediary polynomials (with reference to Section~\ref{sec: asymptotics A-structures}) in the first column of Figure~\ref{table: RQT findex free}, which leads to Proposition~\ref{prop: asymptotics findex} below.

\bgroup
\def\arraystretch{1.5}
\begin{figure}[htp]
\begin{center}
\begin{tabular}{c|c|c|}
& $S_3^\findex = \frac13\,z^3 + z^2$ & $S_3^{(0)} = \frac13z^3+z^2$ \\
\hline
$R(z)$ & $-\frac13 z^2 $ & $-\frac{16}{81}z^3 + \frac49z^2 - \frac23z$ \\
\hline
$Q(z)$ & $ - \frac13 z$ & $-\frac23z^2 + \frac29z$ \\
\hline
$T(z)$ & $z$ & $z^2 - \frac23z + \frac49$ \\
\hline
\end{tabular}
\caption{Intermediary polynomials in the computation of asymptotics for $T_3^\findex$ and $T_3^{(0)}$}\label{table: RQT findex free}
\end{center}
\end{figure}
\egroup%

\begin{proposition}\label{prop: asymptotics findex}
The number $H_n^\findex$ of index $n$ subgroups of $\PSL$ satisfies the following:
\begin{equation}\label{eq: Hnfi}
H_n^\findex \sim \frac{e^{-\frac14}n^\frac12}{\sqrt{12\pi}}\ n^{\frac{n}6} \exp\left(-\frac n6 + n^\frac12 + n^\frac13\right).
\end{equation}
\end{proposition}

\proof
Proposition~\ref{pro:equiv exp polynom} applies to the EGS $T_3^\findex$ of permutational $\tau_3$-structures. In view of Table~\ref{table: RQT findex free}, it yields the following asymptotic equivalent:
$$[z^n]T_3^\findex(z) \sim \frac1{\sqrt{6\pi n}}\ n^{-\frac n3} \exp\left(\frac n3 + n^\frac13\right).$$
Equation~\eqref{eq: tilde Gfi} then yields an asymptotic equivalent of the coefficients of $\widetilde G^\findex$, namely
\begin{equation}\label{eq: eq tilde Gfi}
[z^n]\widetilde G^\findex(z) \sim \frac{e^{-\frac14}}{\sqrt{12\pi n}}\ n^{\frac n6} \exp\left(-\frac n6 + n^\frac12 + n^\frac13\right).
\end{equation}
We can then apply Proposition~\ref{prop: bender revisited} to $\widetilde G^\findex$, getting
\begin{equation}\label{eq: eq Gfi}
[z^n]G^\findex(z) \sim \frac{e^{-\frac14}}{\sqrt{12\pi n}}\ n^{\frac n6} \exp\left(-\frac n6 + n^\frac12 + n^\frac13\right).
\end{equation}
Equation~\eqref{eq: Hnfi} then follows.
\eop

\begin{remark}
Equation~\eqref{eq: Hnfi} coincides with Newman's result \cite[Thm 4]{1976:Newman}.
\end{remark}

Comparing with the asymptotic equivalent of $H_n$ (Theorem~\ref{thm: asymptotics for H} above), we get the following statement.

\begin{proposition}\label{prop: proba finite index}
The probability that a size $n$ $\PSL$-cyclically reduced subgroup of $\PSL$ has finite index is $\Theta\left(e^{-n^\frac23 -\frac13n^\frac13}\right)$.
\end{proposition}

\begin{remark}
\cite[Prop. 5.1]{2008:BassinoNicaudWeil} shows that the probability for a size $n$ subgroup of a rank 2 free group to have finite index is $\Theta\left(n^\frac12 e^{-4n^\frac12}\right)$.
\end{remark}

\paragraph*{Random generation.} The discussion of Section~\ref{sec: random generation} can also be adapted to the random generation of a finite index subgroup of index (size) $n$. As in the general case, we need to draw uniformly at random a $\PSL$-cyclically reduced graph of size $n$ but, in contrast with the general case, we do not need to decide between keeping this graph or deleting one of its loops.

Drawing uniformly at random a size $n$ $\PSL$-cyclically reduced graph is done here again by a rejection algorithm: one draws uniformly at random a size $n$ $\tau_2$-structure and a size $n$ permutational $\tau_3$-structure, using the procedure outlined in Section~\ref{sec: exact counting and random generation}; if they do not form a connected graph, they are rejected and the process is repeated. Proposition~\ref{prop: bender revisited} (or the comparison between Equations~\eqref{eq: eq tilde Gfi} and \eqref{eq: eq Gfi}) justifies the efficiency of the procedure.

\paragraph*{Isomorphism type.} Finally, as in Section~\ref{sec: asymptotics expected values}, we compute the expected value of the isomorphism type $(\ella^\findex, \ellb^\findex, r^\findex)$ of a finite index subgroup and we show a strong concentration phenomenon around these values.

Recall that $\ella^\findex$ is the number of $a$-loops in a $\tau_2$-structure, $\ellb^\findex$ is the number of $b$-loops in a permutational $\tau_3$-structure and $r^\findex = 1 + \frac16(n - 3\ella^\findex - 4\ellb^\findex)$ (Proposition~\ref{prop: compute isomorphism type}). We note that the random variables $\ella^\findex$ and $\ella$ coincide, so the expected value computed for $\ella$ in Theorem~\ref{thm: expected isomorphism type} also applies to $\ella^\findex$ --- on the set of $\tau_2$-structures as well as on the set of finite index subgroups, by the same reasoning as in Theorem~\ref{thm: expected isomorphism type}. Müller and Schlage-Puchta computed the expected value for $\ellb^\findex$ \cite[Corollary 12]{2004:MullerSchlage-Puchta} and showed that it is equivalent to $n^\frac13$ (it could also be computed as in Theorem~\ref{thm: expected isomorphism type}). Thus we get the following proposition.

\begin{proposition}\label{prop: expected loops findex}
The expected value of $\ella^\findex$ (resp. $\ellb^\findex$) is equivalent to $n^\frac12$ (resp. $n^\frac13$).
The expected value of $r^\findex$ is $\frac16n - \frac12n^\frac12 + o(n^\frac12)$.
\end{proposition}

The concentration result mentioned above, obtained directly from Proposition~\ref{pro:large dev}, is the following.

\begin{theorem}\label{thm: large deviations for finite index}
For every $0 < \lambda < 1 < \mu$, there exists $0 < \gamma < 1$ such that
\begin{itemize}
\item $\P(\ella^\findex \le \lambda n^\frac12) = \O(\gamma^{n^\frac12})$ and $\P(\ella^\findex \ge \mu n^\frac12) = \O(\gamma^{n^\frac12})$;
\item $\P(\ellb^\findex \le \lambda n^\frac13) = \O(\gamma^{n^\frac13})$ and $\P(\ellb^\findex \ge \mu n^\frac13) = \O(\gamma^{n^\frac13})$;
\item $\P(r^\findex \le 1 + \frac16(n - 3\mu n^\frac12 - 4\mu n^\frac13)) =  \O(\gamma^{n^\frac13})$;
\item and $\P(r^\findex \ge 1 + \frac16(n - 3\lambda n^\frac12 - 4\lambda n^\frac13)) =  \O(\gamma^{n^\frac13})$.
\end{itemize}
\end{theorem}

\proof
The reasoning is exactly the same as for Theorem~\ref{thm: large deviation ella}. Note that a common value for $\gamma$ can be chosen, by taking the maximum of the values given by Proposition~\ref{pro:large dev} for $\ella^\findex$, $\ellb^\findex$ and $r^\findex$, respectively.
\eop

\subsection{Free subgroups}\label{sec: free subgroups}

We now turn to free subgroups, that is, those whose Stallings graph $\Gamma$ satisfies $\ella = \ellb = 0$ (Proposition~\ref{prop: compute isomorphism type}). This happens if, either $\Gamma$ is $\PSL$-cyclically reduced and loop-free, or if $\Gamma$ is not $\PSL$-cyclically reduced and is obtained by deleting a loop from a $\PSL$-cyclically reduced graph with the same number of vertices and exactly 1 loop.

Let $\Gpr^{(0)}$ (resp. $\Gpr^{(1)}$) be the EGS of labeled proper $\PSL$-cyclically reduced graphs without any loops (resp. with a single loop). Let also $H_n^\ff$ and $H_n^\crff$ be the number of size $n$ free subgroups and cyclically free subgroups of $\PSL$, respectively. Then 
\begin{equation}\label{eq: Hnff}
H_n^\crff = n\ [z^n]\Gpr^{(0)}\quad\text{and}\quad H_n^\ff = n\ [z^n]\Gpr^{(0)} + [z^n]\Gpr^{(1)}.
\end{equation}

\paragraph*{Loop-free $\PSL$-cyclically reduced graphs.}
Let us start with labeled loop-free $\PSL$-cyclically reduced graphs. As we will see below, this is the most probable case. We proceed as in the previous sections, except that
\begin{itemize}
\item the EGS $T_2$ must be replaced by $T_2^{(0)}(z) = \exp\left(\frac{z^2}2\right)$, the EGS of loop-free $\tau_2$-structures;

\item and the EGS $T_3$ must be replaced by $T_3^{(0)}(z,v) = \exp\left(z^2v + \frac{z^3}3\right)$, the EGS of loop-free $\tau_3$-structures where variable $v$ counts the isolated $b$-edges.
\end{itemize}
Let $\widetilde\Gpr^{(0)}(z,v)$ be the EGS of pairs of a loop-free $\tau_2$-structure and a loop-free $\tau_3$-structure, and $\Gpr^{(0)}(z,v)$ the EGS of the connected objects enumerated by $\widetilde\Gpr^{(0)}$. Then, again, we have
\begin{equation}\label{eq: tilde G0}
[z^n]\widetilde\Gpr^{(0)}(z,1) = n! [z^n]T_2^{(0)}\ [z^n]T_3^{(0)}(z,1)\quad\text{and}\quad1 + \widetilde \Gpr^{(0)} = \exp(\Gpr^{(0)}),
\end{equation}
which allows exact computation (see  Appendix~\ref{app:coef}) and random generation as per Section~\ref{sec: random generation}.

We observe that the odd coefficients of $T_2^{(0)}$ are zero, and hence so are the odd coefficients of $[z^n]\widetilde\Gpr^{(0)}(z,1)$ and $[z^n]\Gpr^{(0)}(z,1)$. In order to estimate asymptotically the even coefficients of these series, we use the results in Section~\ref{sec: asymptotics A-structures}.

\begin{proposition}\label{prop: asymptotics free}
Asymptotic equivalents for $[z^{2n}]T_2^{(0)}$, $[z^{2n}]T_3^{(0)}(z,1)$, $[z^{2n}]\widetilde G^{(0)}(z,1)$ and $[z^{2n}] G^{(0)}(z,1)$ are as follows.
\begin{align}
[z^{2n}]T_2^{(0)}(z) &\sim \frac1{\sqrt{2\pi n}} \exp\left(-n\log n + (1-\log 2)n\right) \label{eq: eq T20} \\
[z^{2n}]T_3^{(0)}(z,1) & \sim \frac{e^{\frac{4}{9}}}{\sqrt{12\pi n}}\exp\left(-\frac23 n\log n + \frac23(1-\log 2)n + 2^\frac23 n^\frac23 -\frac{2^\frac43}3n^\frac13\right) \label{eq: eq T30} \\
[z^{2n}]\widetilde \Gpr^{(0)}(z,1) &\sim \frac{e^{\frac4{9}}}{\sqrt{6\pi}}\ \exp\left(\frac13n\log n - \frac13(1-\log2)n + 2^\frac23n^\frac23 - \frac{2^\frac43}3n^\frac13 -\frac12\log n\right) \label{eq: eq tilde G0}\\
[z^{2n}]\Gpr^{(0)}(z,1) &\sim \frac{e^{\frac4{9}}}{\sqrt{6\pi}}\ \exp\left(\frac13n\log n - \frac13(1-\log2)n + 2^\frac23n^\frac23 - \frac{2^\frac43}3n^\frac13 -\frac12\log n\right) \label{eq: eq G0}
\end{align}
\end{proposition}

\proof
Equation~\eqref{eq: eq T20} is a direct application of Proposition~\ref{pro:equiv exp polynom} to the series $\widehat T^{(0)}_2(z) = \frac z2$, since $[z^{2n}]T_2^{(0)} = [z^n]\widehat T^{(0)}_2$. It can also be obtained using Stirling's equivalent for the factorial function, since $[z^{2n}]T_2^{(0)} = \frac{2^{-n}}{n!}$.

Equation~\eqref{eq: eq T30} is also an application of Proposition~\ref{pro:equiv exp polynom}, in view of Table~\ref{table: RQT findex free}. Equation~\eqref{eq: eq tilde G0} then follows from \eqref{eq: tilde G0}.

Proposition~\ref{prop: bender revisited} applies to the series $\widehat G^{(0)}(z) = \sum_n [z^{2n}] \widetilde G^{(0)}(z,1)\ z^n$, and it yields
\begin{equation*}
\frac{[z^{2n}]\Gpr^{(0)}(z,1)}{[z^{2n}]\widetilde\Gpr^{(0)}(z,1)} = 1 - (2n)^{-\frac13} + o(n^{-\frac13}),
\end{equation*}
which justifies Equation~\eqref{eq: eq G0}.
\eop

We can also compute the expected value of the parameter $\kb$ and of the rank of a subgroup with a $\PSL$-cyclically reduced Stallings graph, on the set of loop-free $\tau_3$-structures and $\PSL$-cyclically reduced Stallings graphs (denoted $\kb^\ff$ and $r^\ff$) by the same method as in Theorem~\ref{thm: expected isomorphism type}, using Eq.~(\ref{eq: eq T30}) on the asymptotic behavior of the coefficients of the EGS of loop-free $\tau_3$-structures. This is then used to prove a large deviation statement, exactly as in Theorems~\ref{thm: large deviation ella} and~\ref{thm: large deviations for finite index}.

\begin{proposition}\label{prop: expected value kb cyclically reduced}
The expected value of $\kb^\ff$ on a size $2n$ loop-free $\PSL$-cyclically reduced graph is asymptotically equivalent to $(2n)^\frac23$.

The expected value of the rank $r^\ff$ of a size $2n$ free subgroup whose Stallings graph is $\PSL$-cyclically reduced, is $\frac16(2n-(2n)^\frac23 + o(n^\frac23))$.
\end{proposition}

\begin{theorem}\label{thm: large deviations for free}
For every $0 < \lambda < 1 < \mu$, there exists $0 < \gamma < 1$ such that
\begin{itemize}
\item $\P(\kb^\ff \le \lambda (2n)^\frac23) = \O(\gamma^{n^{2/3}})$ and $\P(\kb^\ff \ge \mu (2n)^\frac23) = \O(\gamma^{n^{2/3}})$;
\item $\P(r^\ff \le \frac16(n - 2\mu (2n)^\frac23)) =  \O(\gamma^{n^{2/3}})$;
\item and $\P(r^\ff \ge \frac16(n - 2\lambda (2n)^\frac23)) =  \O(\gamma^{n^{2/3}})$.
\end{itemize}
\end{theorem}

\paragraph*{Non $\PSL$-cyclically reduced graphs.}\label{par:one loop}
As we saw earlier, these graphs are obtained by removing the loop from a same size labeled $\PSL$-cyclically reduced graph $\Gamma$ with a single loop.
Either this loop is an $a$-loop at a vertex of a $b$-triangle, or it is an $a$- or $b$-loop sitting at the end of a path (of alternating isolated $a$- and $b$-edges), which ends in an isolated $a$-edge connected to a $b$-triangle. Moreover, removing that loop and the path connecting it to a $b$-triangle (including its last vertex, on the $b$-triangle) --- a path that we call the \emph{access path} --- yields a loop-free connected $\PSL$-cyclically reduced graph $\Delta$, called the \emph{cyclically reduced core} of $\Gamma$. In this process, the $b$-triangle in $\Gamma$ on which the access path is attached disappears, leaving a new isolated $b$-edge in the cyclically reduced core.

The set of labeled $\PSL$-cyclically reduced graphs with a single loop is therefore in bijection with the set of pairs of an access path (as described above), and a labeled loop-free $\PSL$-cyclically reduced graph with a designated isolated $b$-edge.

The EGS of  access paths is
\begin{equation}\label{eq: EGS P}
P(z,v) = \frac{z(z+1)}{1-2z^2v} = \sum_{n = 1}^\infty (2^{n-1}z^{2n-1}v^{n-1} + 2^{n-1}z^{2n}v^{n-1}).
\end{equation}
Since the process of building a 1-loop graph from an access path and a loop-free connected $\PSL$-cyclically reduced graph with a marked isolated $b$-edge, reduces the number of isolated $b$-edges by 1, we have
\begin{equation}\label{eq: EGS G1}
\Gpr^{(1)}(z,v) = P(z,v) \ \frac\partial{\partial v}\Gpr^{(0)}(z,v).
\end{equation}
This allows us to compute the coefficients of $\Gpr^{(1)}$, by the following formulas:
\begin{itemize}
\item $[z^0v^k]\Gpr^{(1)} = 0$ for all $k$;
\item for $n\ge 1$ and $k\le n-1$,
\begin{align*}
[z^{2n-1}v^k]\Gpr^{(1)} = [z^{2n}v^k]\Gpr^{(1)} &= \sum_{h=1}^{k+1} [z^{2h}v^{h-1}]P \ [z^{2(n-h)}v^{k-h+1}]\left(\frac\partial{\partial v}\Gpr^{(0)}\right) \\
&= \sum_{h=1}^{k+1} 2^{h-1} (k-h+2)\ [z^{2(n-h)}v^{k-h+2}]\Gpr^{(0)}.
\end{align*}
\end{itemize}
Then,  for $n\ge 1$, we have
\begin{align}
[z^{2n-1}]\Gpr^{(1)}(z,1) = [z^{2n}]\Gpr^{(1)}(z,1)&= \sum_{k=0}^{n-1} [z^{2n}v^k]\Gpr^{(1)} \nonumber \\
&= \sum_{k=0}^{n-1} \sum_{h=1}^{k+1} 2^{h-1} (k-h+2)\ [z^{2(n-h)}v^{k-h+2}]\Gpr^{(0)} \nonumber \\
&= \sum_{h=1}^{n} \sum_{k=h-1}^{n-1} 2^{h-1} (k-h+2)\ [z^{2(n-h)}v^{k-h+2}]\Gpr^{(0)} \nonumber \\
&= \sum_{h=1}^{n} 2^{h-1} \sum_{k=1}^{n-h+1} k\ [z^{2(n-h)}v^k]\Gpr^{(0)}\nonumber \\
&= \sum_{h=0}^{n-1} 2^{n-h-1} \sum_{k=1}^{h+1} k\ [z^{2h}v^k]\Gpr^{(0)}.\label{eq: sum for G1}
\end{align}
Using Eq.~(\ref{eq: Hnff}), we can now compute the exact values of the $H_n^\ff$.

\paragraph*{Random generation.} 
To randomly generate free subgroups, it is convenient to define a more atomic combinatorial description of 1-loop $\PSL$-cyclically reduced graphs. The set of such size $n$ graphs, whose only loop is labeled by $b$, is in bijection with the set of size $n-1$ $\PSL$-cyclically reduced graphs with exactly one loop labeled by $b$: the vertex $s$ carrying the $b$-loop must be at the end of an isolated $a$-edge joining, say, vertex $s$ to vertex $t$; removing $s$ and the adjacent $b$-loop and $a$-edge and adding a $b$-loop at $t$ yields a size $n-1$ $\PSL$-cyclically reduced graph with 1 $b$-loop. Consider now a graph with a single loop, labeled $a$, at vertex $s$. Then $s$ can be part of a $b$-triangle: removing it and the adjacent $b$-edges, as we did for access paths, yields a size $n-1$ graph with a designated isolated $b$-edge. If $s$ is not part of a $b$-triangle, then it is adjacent to an isolated $b$-edge, in which case we proceed as when the loop is labeled by $a$ (with two possibilities, depending on the orientation of the $b$-edge). Let $\Gpra(z)$ (resp.
$\Gprb(z)$) denote the EGS of $\PSL$-cyclically reduced graphs with exactly one loop labeled by $a$ (resp. $b$). We have $\Gpro(z)=\Gpra(z)+\Gprb(z)$ and
\begin{equation}\label{eq:decomp a b}
\Gprb(z) = z\Gpra(z)\text{ and }\Gpra(z) = z\frac{\partial}{\partial v}\Gprz(z,v)\Big|_{v=1}+
2z\Gprb(z).
\end{equation}
After precomputing the coefficients of all EGSs under consideration, we can use
Equations~\eqref{eq: Hnff} and~\eqref{eq:decomp a b} to build a random sampler. By Equation~\eqref{eq: Hnff}, the graph is not loop-free with probability $[z^n]\Gpro(z)/H_n^\text{fr}$, in which case its loop is labeled by $a$ with probability $[z^n]\Gpra(z)/[z^n]\Gpro(z)$ and by $b$ with probability $[z^n]\Gprb(z)/[z^n]\Gpro(z)$. If the loop is labeled by $a$ then we inductively build a $\PSL$-reduced graph with one loop, labeled by $b$ and change it accordingly. If the loop is labeled by $a$, then with probability $2[z^{n-1}]\Gprb(z)/\Gpra(z)$ the vertex with the $a$-loop is adjacent to an isolated $b$-edge, and we proceed inductively. Otherwise, we have to generate a loop-free $\PSL$-cyclically reduced graph of size $n$ with a distinguished isolated $b$-edge, and complete it into a $b$-triangle with an $a$-loop on the new vertex.

All in all, after some inductive steps, we reduced the process to generating a loop-free $\PSL$-cyclically reduced graph, possibly with a distinguished isolated $b$-edge. This is done as in the previous section, by generating pairs of loop-free $\tau_2$- and $\tau_3$-structures until the underlying graph is connected (weighting the graphs by their number of isolated $b$-edges if necessary).

\paragraph*{Asymptotic enumeration.} Before we discuss the asymptotic behavior of loop-free non-$\PSL$-cyclically reduced Stallings graphs, we record two technical statements.

\begin{lemma}\label{lm: about increasing sequence}
We have  $[z^2]\Gpr^{(0)}(z,1) = [z^4]\Gpr^{(0)}(z,1) = 1$ and, for $n\ge 2$, $2[z^{2n}]\Gpr^{(0)}(z,1) \le [z^{2(n+1)}]\Gpr^{(0)}(z,1)$.
\end{lemma}

\proof
We show that, for each $k\le n-1$,  $2[z^{2n}v^k]\Gpr^{(0)} \le [z^{2(n+1)}v^{k+1}]\Gpr^{(0)}$. Let $g^{(0)}_{n,k}$ be the number of labeled $\PSL$-cyclically reduced with $2n$ vertices and $k$ isolated $b$-edges. Then $[z^{2n}v^k]\Gpr^{(0)} = \frac{g^{(0)}_{n,k}}{(2n)!}$, and we want to show that $2(2n+1)(2n+2)g^{(0)}_{n,k} \le g^{(0)}_{n+1,k+1}$.

Let $\Gamma$ be a $\PSL$-cyclically reduced graph with $2n$ vertices (labeled 1 to $2n$) and $k$ isolated $b$-edges.

Let $(i_1,j_1)$ and $(i_2,j_2)$ be distinct isolated $a$-edges, say with $\min(i_1,j_1) < \min(i_2,j_2)$. There are $n$ isolated $a$-edges, and $\binom n2$ choices for these pairs. Let $\Delta$ be the graph obtained from $\Gamma$ by deleting these $a$-edges, adding an isolated $b$-edge between new vertices $2n+1$ and $2n+2$, and adding $a$-edges between $i_1$ or $i_2$ and one of $2n+1$ and $2n+2$, between $j_1$ or $j_2$ and the remaining new vertex, and between the remaining $x_1$ and $x_2$. Considering the orientation of the added $b$-edge, this gives rise to $8n(n-1)$ different labeled graphs $\Delta$, all of them $\PSL$-cyclically reduced, with $k+1$ isolated $b$-edges. Exchanging the labels of the vertex labeled $2n+1$ and its $a$-neighbor gives rise to $8n(n-1)$ more labeled graphs; and so does exchanging the labels of vertex $2n+2$ and its $a$-neighbor.

Moreover  if $(i,j)$ is an isolated $a$-edge in $\Gamma$ with, say, $i< j$, we can delete this edge, add an $a$-edge between new vertices $2n+1$ and $2n+2$, and $b$-edges between $i$ and $2n+1$ and between $j$ and $2n+2$, or between $i$ and $2n+2$ and between $j$ and $2n+2$. We can also exchange the label of vertex $2n+1$ with its $b$-neighor, of the label of vertex $2n+2$ and its $b$-neighbor. This gives rise, again, to $12n$ distinct $\PSL$-cyclically reduced graphs with $k+1$ isolated $b$-edges.

Thus we have
$(24n(n-1)+12n) g^{(0)}_{n,k}  \le g^{(0)}_{n+1,k+1},$
and the result follows since $2(2n+1)(2n+2) < 24n(n-1)+12n$ for all $n \ge 2$. 
\eop

The second technical result leads to an asymptotic estimate of the proportion of non $\PSL$-cyclically reduced graphs among Stallings graphs of free groups. 

\begin{lemma}\label{lm: bounds for G1}
$\Theta\left(n^{-\frac13(1+\log2)}\right) \enspace\le\enspace \frac{[z^{2n}]\Gpr^{(1)}(z,1)}{[z^{2n}]\Gpr^{(0)}(z,1)} \enspace\le\enspace \Theta(n^\frac23).$
\end{lemma}

\proof
Eq.~(\ref{eq: sum for G1}) shows that
\begin{align*}
\sum_{h=0}^{n-1} 2^{n-h-1} [z^{2h}]\Gpr^{(0)}(z,1) &\enspace\le\enspace  [z^{2n}]\Gpr^{(1)}(z,1) \enspace\le\enspace n\ \sum_{h=0}^{n-1} 2^{n-h-1} [z^{2h}]\Gpr^{(0)}(z,1).
\end{align*}
For each $h\ge 0$, let $S_h = 2^{n-h-1} [z^{2h}]\Gpr^{(0)}(z,1)$. Then
$\sum_{h = 0}^{n-1}S_h \le [z^{2n}]\Gpr^{(1)} \le n\,\sum_{h = 0}^{n-1}S_h.$
Let $S^- =  \sum_{h=0}^{\lfloor n-\log n\rfloor}S_h$ and $S^+ = \sum_{h = \lceil n-\log n\rceil}^{n-1} S_h$.
Since the sequence $(S_h)_h$ is increasing (by Lemma~\ref{lm: about increasing sequence}), we have
$$(n-\log n) 2^{n-1}  \le S^- \le (n-\log n) S_{\lfloor n-\log n\rfloor} \le (n-\log n) 2^{\log n} [z^{2\lfloor n-\log n\rfloor}]\Gpr^{(0)}(z,1),$$
and hence
$$\Theta\left(-\frac13n\log n + o(n\log n)\right)\enspace\le\enspace \frac{S^-}{[z^{2n}]\Gpr^{(0)}(z,1)} \enspace\le\enspace \Theta\left(-\frac13\log^2 n + o(\log^2 n)\right).$$

Let us now consider $n-\log n < h < n$. Then
\begin{equation}\label{eq: equiv quotient of Sh}
\frac{S_{h+1}}{S_h} = \frac{2^{n-(h+1)-1}\ [z^{2(h+1)}]\Gpr^{(0)}(z,1)}{2^{n-h-1}\ [z^{2h}]\Gpr^{(0)}(z,1)} \sim \frac12 (2h)^\frac13 = 2^{-\frac23}h^\frac13.
\end{equation}
It follows that
$$2\inv n^\frac13 \enspace\le\enspace \frac{S_{h+1}}{S_h} \enspace\le\enspace 2^{-\frac23}n^\frac13.$$
Then $S^+$ satisfies
\begin{align*}
S^+ &\enspace\le\enspace S_{\lceil n-\log n\rceil}\ \sum_{p=0}^{\log n-1} \left(2^{-\frac23}n^\frac13\right)^p \\
& \enspace\le\enspace  2^{\log n-1}\ [z^{2\lceil n-\log n\rceil}]\Gpr^{(0)}(z,1) \frac{\left(2^{-\frac23}n^\frac13\right)^{\log n}}{2^{-\frac23}n^\frac13-1}\\
&\enspace\le\enspace \frac{2^{\frac13\log n -1}n^{\frac13\log n}}{2^{-\frac23}n^\frac13-1}[z^{2\lceil n-\log n\rceil}]\Gpr^{(0)}(z,1).
\end{align*}
$S^+$ also satisfies
\begin{align*}
S^+ &\enspace\ge\enspace S_{\lceil n-\log n\rceil}\ \sum_{p=0}^{\log n-1} \left(2\inv n^\frac13\right)^p \\
&\enspace\ge\enspace 2^{\log n-2}[z^{2\lceil n-\log n\rceil}]\Gpr^{(0)}(z,1) \ \frac{\left(2\inv n^\frac13\right)^{\log n}-1}{2\inv n^\frac13} \\
&\enspace\ge\enspace  \frac13\,n^{\frac13(\log n-1)}\ [z^{2\lceil n-\log n\rceil}]\Gpr^{(0)}(z,1).
\end{align*}
It follows that
$$\Theta\left(n^{-\frac13(1+\log2)}\right) \enspace\le\enspace \frac{S^+}{[z^{2n}]\Gpr^{(0)}(z,1)} \enspace\le\enspace \Theta(n^{-\frac13}),$$
and this concludes the proof, since $S^-+S^+ \le [z^{2n}]\Gpr^{(1)} \le n(S^-+S^+)$.
\eop

Lemma~\ref{lm: bounds for G1} yields the following corollaries.

\begin{corollary}\label{cor: proba cyclically reduced for free}
The probability that a size $2n$ free subgroup of $\PSL$ is not $\PSL$-cyclically reduced is $\O\left(n^{-\frac13}\right)$.
\end{corollary}

\proof
The probability under study is the ratio
$$\frac{[z^{2n}]\Gpr^{(1)}(z,1)}{[z^{2n}]\Gpr^{(1)(z,1)}+ 2n[z^{2n}]\Gpr^{(0)}(z,1)} = \frac{\frac{[z^{2n}]\Gpr^{(1)}(z,1)}{[z^{2n}]\Gpr^{(0)}(z,1)}}{2n + \frac{[z^{2n}]\Gpr^{(1)}(z,1)}{[z^{2n}]\Gpr^{(0)}(z,1)}}.$$
An upper bound is deduced immediately from Lemma~\ref{lm: bounds for G1}.
\eop

We also get an asymptotic estimate for the number of free subgroups of $\PSL$.

\begin{corollary}\label{cor: asymptotics Hff}
The number $H_{2n}^\ff$ of size $2n$ free subgroups satisfies
\begin{align*}
H_{2n}^\ff &= \Theta\left(n\,[z^{2n}]\Gpr^{(0)}(z,1)\right) \\
&
= \Theta\left(\exp\left(\frac13n\log n -\frac13(1-\log 2)n + 2^\frac23\,n^\frac23 - \frac{2^\frac43}3\,n^\frac13 + \frac12\log n\right)\right),
\end{align*}
and the number $H_{2n-1}^\ff$ of size $2n-1$ free subgroups satisfies 
\begin{align*}
H_{2n-1}^\ff &\le \Theta\left(\exp\left(\frac13n\log n -\frac13(1-\log 2)n + 2^\frac23\,n^\frac23 - \frac{2^\frac43}3\,n^\frac13 + \frac16\log n\right)\right) \\
H_{2n-1}^\ff &\ge \Theta\left(\exp\left(\frac13n\log n -\frac13(1-\log 2)n + 2^\frac23\,n^\frac23 - \frac{2^\frac43}3\,n^\frac13 - \frac16(5+2\log 2)\log n\right)\right).
\end{align*}
\end{corollary}

Using Corollary~\ref{cor: asymptotics Hff} and our result on the number of subgroups of a given size (Theorem~\ref{thm: asymptotics for H}), we get the following result.

\begin{corollary}
The probability that a size $2n$ (resp. $2n-1$) subgroup of $\PSL$ is free is asymptotically
$\Theta\left(\exp\left(-2^\frac12 n^\frac12 - 2^\frac13 n^\frac13\right)\right)$ 
(resp. $\Theta\left(\exp\left(-2^\frac12 n^\frac12 - 2^\frac13 n^\frac13 -\frac16\log n\right)\right)$).
\end{corollary}

\paragraph*{Isomorphism type.} Let us now turn to the evaluation of the parameter $\kb$.
In a non-$\PSL$-cyclically reduced graph, isolated $b$-edges represent roughly half the edges of the access path. It is therefore important to estimate the length of the access path. It turns out to be very small, with high probability.

\begin{proposition}\label{prop: short access path}
Let $f\colon \N\to\N$ be a non-decreasing function such that $f(n) = o(\sqrt n)$, and $f(n) \ge 8$ for all $n$ (including the case of a constant function). Then the probability that a size $2n$ (resp. size $2n-1$) non-$\PSL$-cyclically reduced graph has an access path of length greater than or equal to $2f(n)$ is $\O\left(n^{-\frac{f(n)-7}4}\right)$.
\end{proposition}

\proof
Recall that $\Gpr^{(1)} = P\times\left(\frac\partial{\partial v}\Gpr^{(0)}\right)$, and that $[z^{2n-1}]P(z,1) = [z^{2n}]P(z,1) = 2^{n-1}$ for all $n\ge 1$ (Equations~(\ref{eq: EGS P}) and~(\ref{eq: EGS G1})). The probability $\P$ under investigation (in the even size case) is $\P = \frac X{[z^{2n}]\Gpr^{(1)}(z,1)}$, where
$$X = \sum_{k=0}^{n-1}\sum_{h = f(n)}^{k+1} [z^{2(h+1)}v^h]P\ [z^{2(n-h-1)}v^{k-h}]\left(\frac\partial{\partial v}\Gpr^{(0)}\right).$$
After the appropriate changes of variables (as in the proof of Equation~\eqref{eq: sum for G1}), we find that
\begin{align*}
X &= \sum_{h=0}^{n-f(n)} 2^{n-h-1} \sum_{k=1}^{h+1} k\, [z^{2h}v^k]\Gpr^{(0)} \\
&\le (n-f(n)+1)\sum_{h=0}^{n-f(n)} 2^{n-h-1} [z^{2h}]\Gpr^{(0)}(z,1) = (n-f(n)+1)\sum_{h=0}^{n-f(n)}S_h,
\end{align*}
with $S_h = 2^{n-h-1} [z^{2h}]\Gpr^{(0)}(z,1)$, as in the proof of Lemma~\ref{lm: bounds for G1}. Let $S^- = \sum_{h=0}^{n-f(n)}S_h$ and $S^+ = \sum_{h = n-f(n)+1}^{n-1}S_h$. Then $[z^{2n}]\Gpr^{(1)}(z,1) \ge S^-+S^+$ and
$$\P \le n\ \frac{S^-}{S^-+S^+} 
\le n\ \frac{S^-}{S^+}.$$
We saw (from Lemma~\ref{lm: about increasing sequence} and Eq.~(\ref{eq: equiv quotient of Sh})) that the sequence $(S_h)_h$ is increasing and that $\frac{S_{h+1}}{S_h} \sim 2^{-\frac23}h^\frac13$. It follows that, for $h > n-f(n)$, $\frac{S_{h+1}}{S_h} \ge 2\inv(n-f(n))^\frac13$, and
$$\P \enspace\le\enspace \frac{n\,(n-f(n)+1)\,S_{n-f(n)}}{S_{n-f(n)+1}\sum_{p = 0}^{f(n)-2}\left(2\inv(n-f(n))^\frac13\right)^p}$$
Observe that 
\[
\frac1{\sum_{p = 0}^{f(n)-2}\left(2\inv(n-f(n))^\frac13\right)^p} \leq 2^{f(n)-2}(n-f(n))^{\frac23-\frac13 f(n)}.
\]
Moreover, for $n$ sufficiently large,
\[
\frac{S_{n-f(n)}}{S_{n-f(n)+1}} \leq 2 (n-f(n))^{-\frac13}.
\]
This yields
\begin{align*}
\P & \le n^{2}2^{f(n)}\, (n-f(n))^{\frac13-\frac13f(n)}
= \exp\left( 2\log n + f(n)\log 2+ \left(\frac13-\frac13f(n)\right)\log(n-f(n))
\right)\\
& = \exp\left(\frac{7-f(n)}3\,\log n + f(n)\log 2 +o(1)\right) = \O\left(\exp\left(\frac{7-f(n)}4 \log n\right)\right),
\end{align*}
as announced.
\eop

\begin{corollary}
The expected value of $\kb$ for a size $n$ free subgroup is asymptotically equivalent to $n^\frac23$. The expected rank of a size $n$ free subgroup is $\frac16(n-n^\frac23 + o(n^\frac23))$.
\end{corollary}

\proof
Let us first consider size $2n$ (resp. $2n-1$), loop-free, non-$\PSL$-cyclically reduced graphs. 
Applying Proposition~\ref{prop: short access path} with $f(n)$ constantly equal to 9, we find that, with probability $1-\O(n^{-\frac12})$, the access path has length at most 18, and therefore has at most 8 isolated $b$-edges, and the cyclically reduced core has size at least $2(n-9)$. In the case where the access path is longer, we have $\kb\le n-1$. It follows that $\E(\kb) \sim (2n)^\frac23$ for these graphs.

Since this expected value is the same as that for $\PSL$-cyclically reduced graphs, we find that $\E(\kb) \sim (2n)^\frac23$ for free subgroups of $\PSL$, and that the expected rank of such a subgroup is $\frac16(2n-(2n)^\frac23 + o(n^\frac23))$.
\eop

\subsection{Finite index free subgroups}\label{sec: finite index free subgroups}

As we have seen before, a subgroup of $\PSL$ is free with finite index $n$, if and only if its Stallings graph is $\PSL$-cyclically reduced, has $n$ vertices, no loops and no isolated $b$-edges.

\begin{remark}
Here we retrieve results obtained by Stothers \cite{1978:Stothers-MathsComp}.
\end{remark}

Let $H^\fffi_n$ be the number of index $n$ free subgroups. Let also $G^\fffi(z)$ be the EGS of labeled proper $\PSL$-cyclically reduced graph with $\ella = \ellb = \kb = 0$. We follow again the broad lines of the discussion in Sections~\ref{sec: counting strategy} and~\ref{sec: count cyclically reduced}, with the following changes:
\begin{itemize}
\item since every finite index subgroup has a $\PSL$-cyclically reduced Stallings graph, we have $H_n^\fffi = n\ [z^n]G^\fffi$;
\item we consider loop-free, permutational $\tau_2$- and $\tau_3$-structures, whose EGSs are $T_2^\fffi(z) = T_2^{(0)}(z) = \exp(\frac{z^2}2)$ and $T_3^\fffi(z) = \exp(\frac{z^3}3)$.
\end{itemize}
Note in particular that $[z^n]T_2^\fffi = 0$ if $n$ is odd, and $[z^n]T_3^\fffi = 0$ if $n$ is not a multiple of 3. It follows that every free finite index subgroup of $\PSL$ has index a multiple of 6.
We already know (Proposition~\ref{prop: compute isomorphism type}) that a free, index $6n$ subgroup of $\PSL$ has rank exactly $n+1$.

 Let $\widetilde G^\fffi(z)$ be the EGS of pairs of a $\tau_2$-structure and a $\tau_3$-structure over $n$ elements, that are both loop-free and permutational, and $\Gpr^\fffi(z)$ the EGS of the connected objects enumerated by $\widetilde\Gpr^\fffi$. Then we have
\begin{equation}\label{eq: tilde Gfffi}
[z^n]\widetilde G^\fffi = n!\ [z^n]T_2^\ff \ [z^n]T_3^\fffi\quad\text{and}\quad1 + \widetilde G^\fffi = \exp(G^\fffi),
\end{equation}
which, again, allows exact computation (see Appendix~\ref{app:coef}) and random generation as per Section~\ref{sec: random generation}.

As for finite index and for free subgroups, we deduce asymptotic equivalents.

\begin{proposition}\label{prop: asymptotics free finite index}
Free finite index subgroups of $\PSL$ have index a multiple of 6. Moreover, the  number $H^\fffi_{6n}$ of index $6n$ free subgroups satisfies
$$H^\fffi_{6n} \sim \frac{\sqrt n}{\sqrt{2\pi}}\ \exp\left(n \log n - (1 - \log 6)n\right).$$
\end{proposition}

\proof
The formulas for $T_2^\fffi(z) = \exp(\frac{z^2}2)$ and $T_3^\fffi(z) = \exp(\frac{z^3}3)$ show that non-zero coefficients occur only for indices that are multiples of 2 and 3.

As noted above, $T_2^\fffi = T_2^{(0)}$ so we already know an asymptotic equivalent of its coefficients (Equation~\eqref{eq: eq T20}). An equivalent of the coefficients of $T_3^\fffi$ is obtained in the same fashion (see Proposition~\ref{prop: asymptotics free}), and we get
$$[z^{3n}]T_3^\fffi = \frac{3^{-n}}{n!} \sim \frac1{\sqrt{2\pi n}} \exp\left(-n\log n + (1-\log 3)n\right).$$
Next Equation~\eqref{eq: tilde Gfffi} yields an equivalent of the coefficients of $\widetilde G^\fffi$, namely
$$[z^{6n}]\widetilde G^\fffi \sim \frac1{\sqrt{2\pi n}}\ \exp\left(n \log n - (1 - \log 6)n\right).$$
Then Proposition~\ref{prop: bender revisited} applies to the series $\widehat G^\fffi(z) = \sum_n [z^{6n}] \widetilde G^\fffi\ z^n$, and it yields
\begin{equation*}
\frac{[z^{6n}]\Gpr^\fffi}{[z^{6n}]\widetilde\Gpr^\fffi} = 1 - \frac5{36} n\inv + o(n\inv).
\end{equation*}
Therefore $[z^{6n}] G^\fffi \sim [z^{6n}]\widetilde G^\fffi$, and the asymptotic equivalent of $H^\fffi_{6n}$ follows.
\eop

Comparing with the asymptotic equivalent of $H^\findex_n$ (Theorem~\ref{thm: asymptotics for H} above), we get the following statement.

\begin{proposition}\label{prop: proba free among finite index}
The probability that an index $6n$ subgroup of $\PSL$ is free is
$$\Theta\left(\exp\left(-6^\frac12n^\frac12 - 6^\frac13n^\frac13\right)\right).$$
\end{proposition}

\paragraph*{Beyond.}

The ideas in this paper, together with those in \cite{2008:BassinoNicaudWeil}, can be used to study the subgroups of $G$, where $G$ is the free product of a finite collection of (finite or infinite) cyclic groups: the program for such a study would start with a fine-grain analysis of the combinatorial conditions that characterize rooted $G$-reduced graphs (Stallings graphs, in one-to-one correspondence with subgroups) and unrooted $G$-cyclically reduced graphs (in one-to-one correspondence with conjugacy classes of subgroups). Further investigation would rely on the study of the EGS of $\tau_p$-structures for other integers than $2$ and $3$, along the lines developed above.

{\small\bibliographystyle{abbrv}
\bibliography{pwbiblio}}


\appendix

\newcommand{\tgpr}{\widetilde{\gpr}}

\section{Explicit computation}\label{app:coef}

Here we give more details on the exact computation of the number of the structures encountered in this paper for each size $n$: $\tau_2$- and $\tau_3$-structures (permutative, loop-free and otherwise), labeled $\PSL$-cyclically reduced graphs (rooted, permutative, loop-free or otherwise) and the numbers of subgroups, finite-index subgroups, free subgroups (cyclically reduced or not) and free finite index subgroups.

Figure~\ref{tab:counting} shows the first values of the numbers of size $n$ subgroups, finite index subgroups, free $\PSL$-cyclically reduced subgroups, free subgroups and free and finite index subgroups. The corresponding Python program is available on this \texttt{arXiv} post\footnote{\url{arXiv:2004.00437}.}.

The sequences of number of finite index subgroups and free finite index subgroups were already registered in OEIS, as sequences A005133\footnote{\url{http://oeis.org/A005133}.} and A062980\footnote{\url{http://oeis.org/A062980}.}, respectively. We added the other sequences to OEIS, where they are referred to as sequences \pw{xxx, xxx and xxx}.

\subsection{Number of finitely generated subgroups of $\PSL$}

Let us start with the total number of size $n$ subgroups of $\PSL$.
In view of the discussion in Section~\ref{sec: counting}, and especially Equations~\eqref{eq: Ln} and~\eqref{eq: Hn} on the number of size $n$ subgroups of $\PSL$, we need to compute the coefficients of the
exponential generating series of the $\tau_2$- and $\tau_3$-structures taking into account the number of $a$- and $b$-loops.

By Equation~\eqref{eq: T2}, the numbers $t_2(n,\ell)$ of size $n$ $\tau_2$-structures with $\ell$ loops satisfy
\[
T_2(z,u) = \sum_{n,\ell\geq 0} \frac{t_2(n,\ell)}{n!}z^nu^\ell = \exp\left(zu+\frac{z^2}2\right).
\]
Differentiating with respect to $z$, we get
\[
\frac\partial{\partial z}T_2(z,u) = (u+z)\ T_2(z,u).
\]
Then, considering the coefficient of $z^nu^\ell$ ($n\ge 1$, $0 \le \ell \le n$), we get
\[
\frac{t_2(n+1,\ell)}{n!} = \frac{t_2(n,\ell-1)}{n!} +  \frac{t_2(n-1,\ell)}{(n-1)!}. 
\]
With the convention that $t_2(n,\ell)=0$ if $\ell<0$ or $\ell>n$ or $n<0$, this yields
\begin{equation}\label{eq:rec tau2}
	\begin{cases}
	t_2(0,0) = 1 \\
	t_2(n,\ell) = t_2(n-1,\ell-1) + (n-1)t_2(n-2,\ell),&\text{for }n,\ \ell\geq0, n+\ell \ne 0 .
	\end{cases}
\end{equation}
The first values for $t_2(n,\ell)$ are given in Figure~\ref{fig:tau2 tau3 values}.

\begin{figure}[htb]
\begin{minipage}{.45\textwidth}
\[
\scriptscriptstyle
\begin{array}{c|ccccccc}
 n\setminus \ell& 0 & 1 & 2 & 3 & 4 & 5 & 6 \\
\hline
0 & 1 & 0 & 0 & 0 & 0 & 0 & 0 \\
\hline
1 & 0 & 1 & 0 & 0 & 0 & 0 & 0 \\
\hline
2& 1 & 0 & 1 & 0 & 0 & 0 & 0 \\
\hline
3& 0 & 3 & 0 & 1 & 0 & 0 & 0 \\
\hline
4& 3 & 0 & 6 & 0 & 1 & 0 & 0 \\
\hline
5& 0 & 15 & 0 & 10 & 0 & 1 & 0 \\
\hline
6& 15 & 0 & 45 & 0 & 15 & 0 & 1 \\
\hline
\end{array}
\]
\[
t_2(n,\ell)
\]
\end{minipage}
\begin{minipage}{.45\textwidth}
\[
\scriptscriptstyle
\begin{array}{c|ccccccc}
n\setminus \ell& 0 & 1 & 2 & 3 & 4 & 5 & 6 \\
\hline
0& 1 & 0 & 0 & 0 & 0 & 0 & 0 \\
\hline
1& 0 & 1 & 0 & 0 & 0 & 0 & 0 \\
\hline
2& 2 & 0 & 1 & 0 & 0 & 0 & 0 \\
\hline
3& 2 & 6 & 0 & 1 & 0 & 0 & 0 \\
\hline
4& 12 & 8 & 12 & 0 & 1 & 0 & 0 \\
\hline
5& 40 & 60 & 20 & 20 & 0 & 1 & 0 \\
\hline
6& 160 & 240 & 180 & 40 & 30 & 0 & 1 \\
\hline
\end{array}
\]
\[
t_3(n,\ell)
\]
\end{minipage}
\caption{The first values of $t_2(n,\ell)$ on the left and $t_3(n,\ell)$ on the right.\label{fig:tau2 tau3 values}}
\end{figure}
Similarly, let  $t_3(n,\ell)$ denote the number of size $n$ $\tau_3$-structures with $\ell$ loops, with the convention that $t_3(n,\ell)=0$ if $\ell<0$ or $\ell>n$ or $n<0$. We have (Equation~\eqref{eq: T3})
\[
T_3(z,u,1) = \sum_{n,\ell\geq 0} \frac{t_3(n,\ell)}{n!}z^nu^\ell = \exp\left(zu+z^2+\frac{z^3}3\right).
\]
Hence, we get
\[
\frac\partial{\partial z}T_3(z,u,1) = (u+2z+z^2)T_3(z,u,1).
\]
We then extract the coefficient of $z^nu^\ell$ for $n\geq 1$ and $\ell\geq 0$:
\[
\frac{t_3(n+1,\ell)}{n!} = \frac{t_3(n,\ell-1)}{n!} + \frac{2t_3(n-1,\ell)}{(n-1)!} 
+ \frac{t_3(n-2,\ell)}{(n-2)!}.
\]
This in turn yields the following inductive formula: for $n,\ell\geq 0$,
\begin{equation}\label{eq:rec tau3}
	\begin{cases}
	t_3(0,0) = 1 \\
	t_3(n,\ell) = t_3(n-1,\ell-1) + 2(n-1)t_3(n-2,\ell)+ (n-1)(n-2)t_3(n-3,\ell),&\text{if $n\ell > 0$.}
	\end{cases}
\end{equation}

\begin{remark}
We could also compute the number $t_3(n,\ell,k)$ of size $n$ $\tau_3$-structures with $\ell$ loops and $k$ isolated $b$-edges: the corresponding EGS is $T_3(z,u,v) = \exp(zu+z^2v+z^3/3)$, which can be handled using the same method by differentiating with respect to $z$.
\end{remark}

\begin{remark}
Using Equations~\eqref{eq:rec tau2} and~\eqref{eq:rec tau3}, one can compute the values of $t_2(m,k)$ and $t_3(m,k)$ for all $m\leq n$ and all $k\leq \ell$ in ${\cal O}(n\ell)$ time in the unit cost model.
\end{remark}

Let $\tgpr(n,\ell)$ denote the number of non-empty labeled graphs with $n$ vertices, whose set of $a$-edges (resp. $b$-edges) is determined by a $\tau_2$-structure (resp. a $\tau_3$-structure), having $\ell$ loops in total. We have, for $n\geq 1$ and $0\leq \ell\leq 2n$:
\begin{equation}\label{eq:tgpr}
\tgpr(n,\ell) = \sum_{i=0}^{\ell}t_2(n,i)\,t_3(n,\ell-i),
\end{equation}
which is the same as Equation~\eqref{eq: tilde G}. This directly provides a way to compute the values of $\tgpr$ from those of $t_2$ and $t_3$ (because of the summation, it takes time ${\cal O}(n\ell^2)$ in the unit-cost model to get all the values up to size $n$ and $\ell$ loops). The first values of $\gpr(n,\ell)$ are given in Figure~\ref{fig:tgpr values}.

\begin{figure}[htb]
\[
\scriptscriptstyle
\begin{array}{c|ccccccccccccc}
n\setminus \ell& 0 & 1 & 2 & 3 & 4 & 5 & 6 & 7 & 8 & 9 & 10 & 11 & 12 \\
\hline
1& 0 & 0 & 1 & 0 & 0 & 0 & 0 & 0 & 0 & 0 & 0 & 0 & 0 \\
\hline
2& 2 & 0 & 3 & 0 & 1 & 0 & 0 & 0 & 0 & 0 & 0 & 0 & 0 \\
\hline
3& 0 & 6 & 18 & 2 & 9 & 0 & 1 & 0 & 0 & 0 & 0 & 0 & 0 \\
\hline
4& 36 & 24 & 108 & 48 & 87 & 8 & 18 & 0 & 1 & 0 & 0 & 0 & 0 \\
\hline
5& 0 & 600 & 900 & 700 & 900 & 240 & 275 & 20 & 30 & 0 & 1 & 0 & 0 \\
\hline
6& 2400 & 3600 & 9900 & 11400 & 10950 & 5400 & 4225 & 840 & 675 & 40 & 45 & 0 & 1 \\
\hline
\end{array}
\]
\caption{The first values of $\tgpr(n,\ell)$.\label{fig:tgpr values}}
\end{figure}

Let $\gpr(n,\ell)$ be the number of those labeled graphs\footnote{We have $\gpr(n,\ell)=g_{n,\ell}$ of Section~\ref{sec: counting}. We use this notation for uniformity within the appendix.} counted by $\tgpr(n,\ell)$ which are connected.
By differentiating Equation~\eqref{eq: G} with respect to $z$, we have
\begin{equation}\label{eq:tilde g}
\frac\partial{\partial z}\widetilde \Gpr(z,u) = (1+\widetilde \Gpr(z,u))\ \frac\partial{\partial z} \Gpr(z,u).
\end{equation}
Since $\widetilde\Gpr(z,u)=\sum_{n,\ell} \frac1{n!}\tgpr(n,\ell)z^nu^\ell$ and
$\Gpr(z,u)=\sum_{n,\ell} \frac1{n!}\gpr(n,\ell)z^nu^\ell$, we can extract the coefficient of $z^{n}u^\ell$ in Equation~\eqref{eq:tilde g} to obtain the following (the sum on $m$ ends at $n-1$ as for all $\ell$, we have $\tgpr(0,\ell)=0$):
\[
\frac{\tgpr(n+1,\ell)}{n!} = \frac{\gpr(n+1,\ell)}{n!} + \sum_{m=0}^{n-1}\sum_{k=0}^\ell
\frac{\gpr(m+1,k)}{m!}  \frac{\tgpr(n-m,\ell-k)}{(n-m)!}. 
\]
There follows an inductive formula for the $\gpr(n,\ell)$ ($n\geq 1$, $\ell\geq 0$), with the convention that $\gpr(n,\ell)=0$ if $n\leq 0$, $\ell<0$ or $\ell\geq 2n$:
\begin{equation}\label{eq:gpr}
	\begin{cases}
	\gpr(1,2)  = 1\\
	\gpr(n,\ell) = \tgpr(n,\ell) - \sum_{m=1}^{n-1}\sum_{k=0}^\ell \binom{n-1}{m-1}
	\gpr(m,k)\,\tgpr(n-m,\ell-k),&\text{ otherwise.} \\
	\end{cases}
\end{equation}
Notice that, if $n > 1$, then  $\gpr(n,\ell)=0$ when $\ell>n$ since a graph with more than $n$ loops cannot be connected. The first values of $\gpr(n,\ell)$ are given in Figure~\ref{fig:gpr values}.

\begin{figure}[htb]
\[
\scriptscriptstyle
\begin{array}{c|ccccccccccccc}
n\setminus \ell& 0 & 1 & 2 & 3 & 4 & 5 & 6 \\
\hline
1& 0& 0& 1 & 0 & 0 & 0 & 0 \\
\hline
2& 2 & 0 & 3 & 0 & 0 & 0 & 0 \\
\hline
3& 0 & 6 & 12 & 2 & 0 & 0 & 0 \\
\hline
4& 24 & 24 & 72 & 24 & 0 & 0 & 0 \\
\hline
5& 0 & 480 & 480 & 360 & 0 & 0 & 0 \\
\hline
6& 1560 & 2880 & 5760 & 4560 & 360 & 0 & 0 \\
\hline
\end{array}
\]
\caption{The first values of $\gpr(n,\ell)$.\label{fig:gpr values}}
\end{figure}

Using Equations~\eqref{eq: Ln} and~\eqref{eq: Hn}, we can then compute the numbers $L_n$
of size $n$ labeled $\PSL$-reduced rooted graphs and $H_n$ of size $n$ subgroups of $\PSL$: for $n\geq 2$,
\[
L_n = \sum_{\ell=0}^n (n+\ell)\gpr(n,\ell)\quad\text{and}\quad H_n = \frac1{n!}L_n.
\]
The first values of $H_n$ are given in Figure~\ref{tab:counting} below.

\subsection{Number of  finite index subgroups}

As explained in Section~\ref{sec: finite index subgroups}, the number of finite index subgroups $H_n^\findex$ is $ng^\findex(n)$, where $g^\findex(n)$ is the number of pairs of $\tau_2$-structures and permutational $\tau_3$-structures that form a connected graph. In particular, we do not need to refine the counting in terms of the number of loops. Let $t_2(n)$ (resp. $t_3^\findex(n)$) denote the number of $\tau_2$-structures (resp. permutational $\tau_3$-structures) of size $n$. As the associated EGS of $\tau_2$-structures is $T_2(z)=\exp(z+z^2/2)$, by differentiating with respect to $z$ and extracting the coefficient in $z^n$ we obtain:

\begin{equation}\label{eq:tau2 n}
\begin{cases}
t_2(0) = t_2(1)  = 1 \\
t_2(n) = t_2(n-1)+(n-1)t_2(n-2) &\text{for }n\geq 2.
\end{cases}
\end{equation}
Similarly, the EGS of permutational $\tau_3$-structures is $T_2^\findex(z)=\exp(z+z^3/3)$, yielding

\begin{equation}\label{eq:tau3fi n}
\begin{cases}
t_3^\findex(0) = t_3^\findex(1) =  t_3^\findex(2) = 1 \\
t_3^\findex(n) = t_3^\findex(n-1)+(n-1)(n-2)t_3^\findex(n-3) &\text{for }n\geq 3.
\end{cases}
\end{equation}
The number $\widetilde{g}_\text{pr}^\findex(n)$ of pairs of a $\tau_2$-structure and a permutational $\tau_3$-structures of size $n$ is $t_2(n)t_3^\findex(n)$. And, from the equation $\widetilde{G}^\findex+1=\exp(G^\findex(z))$, the number
$\gpr^\findex(n)$ of such graphs that are also connected satisfies:
\begin{equation}\label{eq:gprfi}
	\begin{cases}
	\gpr^\findex(0)  = 0\\
	\gpr^\findex(n) = \widetilde{g}_\text{pr}^\findex(n) - \sum_{m=1}^{n-1}\binom{n-1}{m-1}
	\gpr^\findex(m)\,\widetilde{g}_\text{pr}^\findex(n-m)&\text{ for }n\geq 1. \\
	\end{cases}
\end{equation}
The number of finite index subgroups is therefore $H_n^\findex=\frac1{(n-1)!}\gpr^\findex(n)$, see Figure~\ref{tab:counting}.

\subsection{Number of  free subgroups}

Recall that the numbers $H_n^\ff$ and $H_n^\crff$ of size $n$ free and cyclically free subgroups of $\PSL$ satisfies $H_n^\crff = n\ [z^n]\Gpr^{(0)}$ and $H_n^\ff = n\ [z^n]\Gpr^{(0)} + [z^n]\Gpr^{(1)}$ (Equation~\eqref{eq: Hnff}), where $\Gpr^{(0)}$ (resp. $\Gpr^{(1)}$) is the EGS of labeled proper $\PSL$-cyclically reduced graphs without any loops (resp. with a single loop). We compute the coefficients of $\Gpr^{(0)}$ and $\Gpr^{(1)}$ separately.

\paragraph*{Loop-free $\PSL$-cyclically reduced graphs.}
Let $\gpr^{(0)}(n)$ denote the number of size $n$ loop-free labeled proper $\PSL$-cyclically reduced graphs. Let also $\widetilde\gpr^{(0)}(n)$ denote the number of pairs of
loop-free $\tau_2$-structures and loop-free $\tau_3$ structures. We have $\widetilde\gpr^{(0)}(0)=0$ and, for all $n\geq 1$,
\begin{equation}\label{eq:tgpr0}
\widetilde\gpr^{(0)}(n) = t_2^{(0)}(n)\,t_3^{(0)}(n),
\end{equation}
where the EGS of $t_2^{(0)}(n)$ (resp. $t_3^{(0)}(n)$) is $T_2^{(0)}(z)=\exp(z^2/2)$ (resp.  $T_3^{(0)}(z,1)=\exp(z^3/3+z^2)$), see Equation~\eqref{eq: tilde G0}. We already observed that $\widetilde\gpr^{0}(n)=0$ for odd values of $n$.

Since $\frac{d}{dz}T_2^{(0)}(z) = z T_2^{(0)}(z)$, we have $t_2^{(0)}(2n) = (2n-1)\,t_2^{(0)}(2n-2)$ for $n\geq 1$ and hence, for $n\geq 1$,
\[
t_2^{(0)}(2n) = (2n-1)(2n-3)\cdots 1
\]
and $t_2^{(0)}(2n+1)=0$. 

Similarly, $\frac{d}{dz}T_3^{(0)}(z,1) = (z^2+2z) T_3^{(0)}(z,1)$ and we get, by coefficient extraction,
\[
\begin{cases}
t_3^{(0)}(0)=1\\
t_3^{(0)}(1)=0\\
t_3^{(0)}(2)=2\\
t_3^{(0)}(n) = 2(n-1)\,t_3^{(0)}(n-2) + (n-1)(n-2)\,t_3^{(0)}(n-3),&\text{ for }n\geq 3. 
\end{cases}
\]
We can thus compute $t_2^{(0)}(n)$, $t_3^{(0)}(n)$ and
$\widetilde\gpr^{(0)}(n)$ by Equation~\ref{eq:tgpr0}. Proceeding as in the previous sections, we derive from
the equation $1+\widetilde\Gpr^{(0)}(z)=\exp(\Gpr^{(0)}(z))$ that
\begin{equation}\label{eq:gpr0}
\begin{cases}
\gpr^{(0)}(0)= 0,\\
\gpr^{(0)}(n)= \tgpr^{(0)}(n) - \sum_{m=1}^{n-1}\binom{n-1}{m-1} \gpr^{(0)}(m)\,\tgpr^{(0)}(n-m),&\text{ for }n\geq 1.
\end{cases}
\end{equation}

\begin{figure}[htb]
\[
\scriptscriptstyle
\begin{array}{c|ccccccc}
n &0&2&4&6&8&10&12\\
\hline
t_2^{(0)}(n) & 1&1&3&15&105&945&10395 \\
\hline
t_3^{(0)}(n) &1&2&12&160&3920&131040&5346880\\
\hline
\widetilde\gpr^{(0)}(n) & 0&2&36&2400&411600&123832800&55580817600 \\
\hline
\gpr^{(0)}(n) & 0&2&24&1560&282240&84188160&36883123200
\end{array}
\]
\caption{The first values of $t_2^{(0)}(2n)$, $t_3^{(0)}(2n)$, $\widetilde\gpr^{(0)}(2n)$ and $\gpr^{(0)}(2n)$.\label{fig:values g0}}
\end{figure}

\paragraph*{$\PSL$-cyclically reduced graphs with a single loop.}

We distinguish two kinds of $\PSL$-cyclically reduced graphs with a single loop, according to the loop label:
let $\Gpr^{(a)}(z)$ (resp. $\Gpr^{(b)}(z)$) denote the EGS of the size $n$ $\PSL$-cyclically reduced graphs with a single loop labeled by $a$ (resp. by $b$). We have:
\begin{equation}\label{eq:G1}
\Gpr^{(1)}(z) = \Gpr^{(a)}(z)+\Gpr^{(b)}(z).
\end{equation}
Let also $\gpr^{(a)}(n)$ and $\gpr^{(b)}(n)$ denote the number of such graphs.

If a size $n$ $\PSL$-cyclically reduced graphs has a single loop labeled by $b$, then it is necessarily at the end of the access path: if we remove it and add an $a$-loop on the previous vertex, we obtain a  size $n$ $\PSL$-cyclically reduced graphs has a single loop labeled by $a$; moreover this transformation is a bijection. Hence
\begin{equation}\label{eq:Gb}
\Gprb(z) = z\Gpra(z).
\end{equation}
If a size $n$ $\PSL$-cyclically reduced graph has a single loop labeled by $a$, then there are two possibilities: 
\begin{itemize}
\item The access path is empty and the $a$-loop is on a vertex that belongs to a $b$-triangle. If we remove this vertex and the adjacent $b$-edges, we obtain a loop-free $\PSL$-cyclically reduced graph. The transformation is not one-to-one, but the number of ways to go back from such a loop-free graph is exactly its number of isolated $b$-edges.
\item The  access path is not empty and if we remove its end vertex and add a $b$-loop on the previous vertex, we obtain a  size $n$ $\PSL$-cyclically reduced graphs that has a single loop labeled by $b$. There are two ways to perform the inverse transformation, depending on the orientation of the added isolated $b$-edge.
\end{itemize}
Translated into EGS, this yields the following equation, with $\Gprv(z) := \frac\partial{\partial v}\Gprz(z,v)\big|_{v=1}$
\begin{equation}\label{eq:Ga}
\Gpra(z) = z \Gprv(z) + 2z \Gprb(z).
\end{equation}
For all $n$, let $\gpr^{(\alpha)}(n)$ denote the number of graphs counted by $\Gpr^{(\alpha)}(z)$, for $\alpha\in\{0,1,a,b\}$. That is, $\gpr^{(\alpha)}(n):=n![z^n]\Gpr^{(\alpha)}(z)$. From 
Equations~\eqref{eq:G1} and \eqref{eq:Ga} we get:
\begin{equation}\label{eq:g1n and gbn}
\gprb(n)=n\gpra(n-1)\text{ and }\gpro(n)=\gpra(n)+n\gpra(n-1).
\end{equation}
We therefore focus on computing the values of $\gpra(n)$. From Equations~\eqref{eq:Gb} and
\eqref{eq:Ga}, we have
\[
\Gpra(z) = z\Gprv(z)+2z^2\Gpra(z).
\]
If we extract the coefficient in $z^n$, we obtain, for $n\geq2$,
\[
\frac{\gpra(n)}{n!} = \frac{\gprv(n-1)}{(n-1)!} + \frac{2\gpra(n-2)}{(n-2)!}. 
\]
Hence we have
\begin{equation}\label{eq:grpa(n)}
\begin{cases}
\gpra(0) = \gpra(1) = 0\\
\gpra(n) = n\,\gprv(n-1)+2n(n-1)\,\gpra(n-2),&\text{for }n\geq2.
\end{cases}
\end{equation}
The values of $\gpra(n)$, which are equal to $0$ for even $n$, can easily be computed from
Equation~\eqref{eq:grpa(n)}, provided the values of $\gprv(n)$ were already  computed.

To compute the values $\gprv(n)$, we need to refine the counting we did for  loop-free $\PSL$-cyclically reduced graphs, taking the number of isolated $b$-edges into account. Let
$t_3^{(0)}(n,k)$ be the number of size $n$ loop-free $\tau_3$-structures with $k$ pairs. We have:
\[
T_3^{(0)}(z,v) = \sum_{n,k} \frac{t_3^{(0)}(n,k)}{n!}z^n v^k = \exp\left(z^2v+\frac{z^3}3\right).
\]
Hence, by differentiating with respect to $z$ and extracting the coefficient in $z^{n-1}v^k$, we have
\[
\frac{t_3^{(0)}(n,k)}{(n-1)!} = \frac{2t_3^{(0)}(n-2,k-1)}{(n-2)!} + \frac{t_3^{(0)}(n-3,k)}{(n-3)!}. 
\]
Thus we have
\begin{equation}\label{eq:t30n}
\begin{cases}
t_3^{(0)}(0,0) = 1\\ 
t_3^{(0)}(n,k) = 2(n-1)t_3^{(0)}(n-2,k-1) +(n-1)(n-2) t_3^{(0)}(n-3,k),& \text{otherwise,}
\end{cases}
\end{equation}
with the convention that $t_3^{(0)}(n,k)=0$ if $n<0$ or $k<0$ or $2k>n$.

Let $\tgpr^{(0)}(n,k)$ denote the number of size $n$ pairs of loop-free $\tau_2$-structures and
loop-free $\tau_3$ structures with $k$ elements of size $2$. We have
\[
\tgpr^{(0)}(n,k) = t_2^{(0)}(n) t_3^{(0)}(n,k).
\]
And also, if $\gpr^{(0)}(n,k)$ denotes the number of size $n$ loop-free $\PSL$-cyclically reduced graphs with $k$ isolated $b$-edges, from $\widetilde\Gpr^{(0)}(z,v) + 1=\exp(\Gpr^{(0)}(z,v))$, we obtain
\begin{equation}\label{eq:gpr01}
	\begin{cases}
	\gpr^{(0)}(2,1)  = 2\\
	\gpr^{(0)}(n,k) = \tgpr(n,\ell) - \sum_{m=1}^{n-1}\sum_{j=0}^k \binom{n-1}{m-1}
	\gpr^{(0)}(m,j)\,\tgpr^{(0)}(n-m,k-j),&\text{ otherwise.} \\
	\end{cases}
\end{equation}
Since $\gprv(n)$ is $n![z^n]\frac{\partial}{\partial v}\Gprz(z,v)|_{v=1}$, we have
\[
\gprv(n) = \sum_{j=0}^{\lfloor n/2\rfloor} j\, \gprz(n,j),
\]
which can be easily computed from the values of $\gprz(n,j)$. Hence also the values of
$\gpra(n)$ by Equation~\ref{eq:grpa(n)}, then those of $\gpro(n)$ by Equation~\ref{eq:g1n and gbn}. Finally the values of $H_n^\crff$ and $H_n^\ff$ are obtain using Equation~\ref{eq: Hnff}:
\begin{align*}
H_n^\crff &= n\ [z^n]\Gpr^{(0)} = \frac{\gprz(n)}{(n-1)!} \\
H_n^\ff &= n\ [z^n]\Gpr^{(0)} + [z^n]\Gpr^{(1)} = \frac{\gprz(n)}{(n-1)!} + \frac{\gpra(n)}{n!}+\frac{\gpra(n-1)}{(n-1)!}.
\end{align*}
%
The first values of these numbers are given in Figure~\ref{tab:counting}.

\subsection{Number of  finite index free subgroups}

We proceed as in the previous sections. The number of size $n$ loop-free permutational $\tau_3$-structures $t^\text{fr-fi}_3(n)$ is $n![z^n]T_3^\text{fr-fi}(z)$, and since
$T_3^\text{fr-fi}(z)=\exp(z^3/3)$ we have, for $n$ multiple of $3$:
\[
t^\text{fr-fi}_3(n) = \frac{n!}{3^\frac{n}3 \left(\frac{n}3\right)!} =\prod_{\substack{i\in\{1\ldots n\}\\i\not\equiv 0\text{ mod }3}} i.
\]
Thus, the number $\widetilde\gpr^\text{fr-fi}(n)$ of pairs of size $n$ loop-free permutational
$\tau_2$- and $\tau_3$-structures is
\[
\widetilde\gpr^\text{fr-fi}(n) =t_2^{(0)}(n)\ t^\text{fr-fi}_3(n), 
\]
which is non-zero only for multiples of 6.
As before, from $\widetilde\Gpr^\text{fr-fi}(z) + 1 = \exp(\Gpr^\text{fr-fi}(z))$, we derive the number $\gpr^\text{fr-fi}(n)$ of labeled proper $\PSL$-cyclically reduced graphs, namely
\[
\begin{cases}
\gpr^\text{fr-fi}(0) = 0\\
\gpr^\text{fr-fi}(n) = \widetilde\gpr^\text{fr-fi}(n) - \sum_{m=1}^{n-1}\binom{n-1}{m-1}\gpr^\text{fr-fi}(m)\widetilde\gpr^\text{fr-fi}(n-m).  
\end{cases}
\]
Finally, the number of finite index free subgroups is $H_n^\text{fr-fi}=\frac{\gpr^\text{fr-fi}(n)}{(n-1)!}$, see Figure~\ref{tab:counting}.

\begin{figure}[htp]
\caption{Number of finitely generated subgroups of size up to 36}
\begin{center}
\begin{tabular}{|l||r|r|r|r|r|}
\hline
\small size & \small all subgroups & \small finite index & \small cycl. reduced free & \small free & \small free finite index \\
\hline
1 & 4 & 1 & 0 & 0 & 0 \\ 
2 & 8 & 1 & 2 & 2 & 0 \\ 
3 & 16 & 4 & 0 & 1 & 0 \\ 
4 & 34 & 8 & 4 & 5 & 0 \\ 
5 & 76 & 5 & 0 & 4 & 0 \\ 
6 & 167 & 22 & 13 & 17 & 5 \\ 
7 & 366 & 42 & 0 & 12 & 0 \\ 
8 & 846 & 40 & 56 & 68 & 0 \\ 
9 & 1870 & 120 & 0 & 37 & 0 \\ 
10 & 4353 & 265 & 232 & 269 & 0 \\ 
11 & 9900 & 286 & 0 & 130 & 0 \\ 
12 & 23054 & 764 & 924 & 1054 & 60 \\ 
13 & 53402 & 1729 & 0 & 492 & 0 \\ 
14 & 125379 & 2198 & 3768 & 4260 & 0 \\ 
15 & 293372 & 5168 & 0 & 1908 & 0 \\ 
16 & 694884 & 12144 & 15936 & 17844 & 0 \\ 
17 & 1641018 & 17034 & 0 & 7584 & 0 \\ 
18 & 3912272 & 37702 & 68817 & 76401 & 1105 \\ 
19 & 9319816 & 88958 & 0 & 31104 & 0 \\ 
20 & 22348358 & 136584 & 301524 & 332628 & 0 \\ 
21 & 53622232 & 288270 & 0 & 131025 & 0 \\ 
22 & 129319050 & 682572 & 1343388 & 1474413 & 0 \\ 
23 & 312184204 & 1118996 & 0 & 563574 & 0 \\ 
24 & 756855652 & 2306464 & 6087376 & 6650950 & 27120 \\ 
25 & 1837195988 & 5428800 & 0 & 2470536 & 0 \\ 
26 & 4475381885 & 9409517 & 27997712 & 30468248 & 0 \\ 
27 & 10918047864 & 19103988 & 0 & 11028448 & 0 \\ 
28 & 26714414272 & 44701696 & 130532224 & 141560672 & 0 \\ 
29 & 65467869902 & 80904113 & 0 & 50054608 & 0 \\ 
30 & 160853707175 & 163344502 & 616603418 & 666658026 & 828250 \\ 
31 & 395841123048 & 379249288 & 0 & 230641440 & 0 \\ 
32 & 976352297396 & 711598944 & 2949326656 & 3179968096 & 0 \\ 
33 & 2411988448210 & 1434840718 & 0 & 1077886298 & 0 \\ 
34 & 5970888317052 & 3308997062 & 14274174272 & 15352060570 & 0 \\ 
35 & 14803858849928 & 6391673638 & 0 & 5105099252 & 0 \\ 
36 & 36772848298022 & 12921383032 & 69861695744 & 74966794996 & 30220800 \\ 
\hline
\end{tabular}
\end{center}
\label{tab:counting}
\end{figure}

\end{document}